\documentclass[12pt,a4paper]{article}
\usepackage{latexsym}
\usepackage{mathrsfs}
\usepackage{amsmath}
\usepackage{amssymb}
\usepackage[usenames,dvipsnames,svgnames,table]{xcolor}
\usepackage{mathtools}

\newtheorem{theorem}{Theorem}[section]
\newtheorem{lemma}[theorem]{Lemma}

\def \Z {\mathbb{Z}}

\def \F {\overrightarrow F}
\def \S {{\cal{S}}}
\def \D {{\mathcal{D}}}

\def \deg {{\rm deg}}

\def \proof {\noindent{\bf Proof}\quad}

\newcommand{\qed}{\hfill$\Box$\vspace{0.2cm}}

\newcommand{\y}{\infty}

\title{\bf Hamilton decompositions of line graphs}

\author{
Darryn Bryant, Sara Herke,
Barbara Maenhaut\\ and Benjamin R. Smith
\\
{\small School of Mathematics and Physics,  The University of Queensland, QLD 4072, Australia.}\\
{\small \texttt{db@maths.uq.edu.au}, 
\, \texttt{bmm@maths.uq.edu.au},\,  \texttt{s.herke@uq.edu.au}\,} \\ {\small \texttt{bsmith.maths@gmail.com}}
}

\date{ }

\begin{document}
\maketitle\thispagestyle{empty}
\def\baselinestretch{1.25}\small\normalsize

\begin{abstract}
It is proved that if a graph is regular of even degree and contains a Hamilton cycle, 
or regular of odd degree and contains a Hamiltonian $3$-factor, 
then its line graph is Hamilton decomposable. This result partially extends Kotzig's result that a $3$-regular graph is Hamiltonian if and only if its line graph is Hamilton decomposable, 
and proves the conjecture of Bermond that 
the line graph of a Hamilton decomposable graph is Hamilton decomposable.
\end{abstract}

\section{Introduction}

Hamilton decomposability of line graphs has been studied extensively. The {\em line graph} of a graph $G$, denoted by $L(G)$, is the graph with a vertex corresponding to each edge of $G$, and in which two vertices are adjacent if and only if their corresponding edges are adjacent in $G$. A {\em Hamilton decomposition} of a graph $G$ is a set of Hamilton cycles in $G$ whose edge sets partition the edge set of $G$. A graph that has a Hamilton decomposition is said to be {\em Hamilton decomposable}.  A landmark result on the topic of Hamilton decomposability of line graphs, due to Kotzig \cite{Kot}, 
is that a 3-regular graph is Hamiltonian if and only if its line graph is Hamilton decomposable.
The goal of this paper is to prove the following theorem which addresses the extension of 
Kotzig's result to graphs of larger degree.

\begin{theorem}\label{mainthm}
If a graph is regular of even degree and contains a Hamilton cycle, or regular of odd degree and contains a Hamiltonian $3$-factor, 
then its line graph is Hamilton decomposable.
\end{theorem}

Theorem \ref{mainthm}, the proof of which follows immediately from Lemmas \ref{combineHamFrags}, \ref{alln} and \ref{alln_odd}, shows that the ``only if'' part of Kotzig's result holds for regular graphs of even degree, and comes close 
to showing that it holds for regular graphs of odd degree. It is possible that just the existence of a Hamiltonian 
cycle, rather than a Hamiltonian $3$-factor, in a regular graph $G$ of odd degree is also sufficient for Hamilton decomposability 
of the line graph of $G$, but we are unable to prove this using our methods.
On the other hand, it has recently been shown \cite{BryMaeSmi} that for all $k\geq 4$,
the existence of a Hamilton cycle in a $k$-regular graph $G$ is not necessary for Hamilton decomposability 
of the line graph of $G$.  

Theorem \ref{mainthm} proves, and considerably strengthens, a long-standing conjecture of Bermond \cite{Ber} which states that  
the line graph of a Hamilton decomposable graph is Hamilton decomposable. Bermond's conjecture was proved by
Jaeger \cite{Jae} in the case of regular graphs of degree $4$, and then
by Muthusamy and Paulraja \cite{MutPau} in the case of regular graphs of degree divisible by $4$. 
Also in \cite{MutPau}, it was shown that the line graph of a Hamiltonian regular graph of even degree can be decomposed into 
Hamilton cycles and a $2$-factor. This result was independently proved by Zahn \cite{Zah}. 
Other results relating to Hamilton decompositions of line graphs can be found in \cite{FleHilJac,HeiVer,Jac,JacWor,Pik1,Pik2,Pik3,Ver}
and in the survey on Hamilton decompositions \cite{AlsBerSot}.

A brief overview of the central elements of the construction used to prove 
Theorem \ref{mainthm} is as follows.
This description is for the case of regular graphs of even degree. Some additional complications are involved in 
the case of regular graphs of odd degree. 
In the line graph $L(G)$ of a $2n$-regular graph $G$, 
the $2n$ vertices of $L(G)$ that correspond to the $2n$ edges of $G$ that are incident with a vertex $v$ of $G$ induce
a complete subgraph of order $2n$ in $L(G)$, and we denote this complete subgraph by $L(G)_v$. 
By a well-known theorem of Petersen \cite{Pet}, if $G$ is Hamiltonian, then $G$ has a 
$2$-factorisation $\cal F$ in which one of the $2$-factors is a Hamilton cycle.

In Section \ref{Section3}, we define a Hamilton fragment to be a subgraph $H$ of a complete graph of order $2n$
such that for any given Hamiltonian $2n$-regular graph $G$, 
and any given $2$-factorisation $\cal F$ of 
$G$ containing a Hamilton cycle,
if a copy of $H$ is placed on the vertices of $L(G)_v$ for each vertex $v$ of 
$G$,
in a manner prescribed by $\cal F$, then the resulting subgraph of $L(G)$ 
has a Hamilton decomposition. We prove various conditions under which $H$ is a Hamilton fragment, and then 
show that the complete graph of order $2n$ can be decomposed into Hamilton fragments. 
The union of the resulting Hamilton decompositions is thus a Hamilton decomposition of $L(G)$. 

The paper is structured as follows. In Section \ref{Section2} we prove several technical lemmas which are used later in the paper. Section \ref{Section3} is divided into two subsections for the two cases of regular graphs of even and odd degree. The main goal of
Section \ref{Section3} is to prove conditions under which a subgraph $H$ is a Hamilton fragment. In Sections \ref{Section4} and \ref{Section5}, the required 
decompositions of complete graphs into Hamilton fragments are given.
If we are after a Hamilton decomposition of $L(G)$, then the order of the complete 
graph to be decomposed into Hamilton fragments is equal to the degree of $G$.
Section \ref{Section4} gives the required 
decompositions of complete graphs into Hamilton fragments for orders 12, 16, 18, and for all orders greater than 19. 
Decompositions for the other small orders 
require different methods than those used in the general case 
and are given in Section \ref{Section5}.

\section{From $2$-factorisations to Hamilton decompositions}\label{Section2}

Let $V$ be a set of vertices and let 
$E_1,E_2,\ldots,E_r$ be pairwise disjoint sets of edges.
The set $\{E_1,E_2,\ldots,E_r\}$ is said to be a {\em $V$-connector}
if for any $2$-factorisation $\{F_1,F_2,\ldots,F_r\}$ of any $2r$-regular graph $G$ 
with $E_i\subseteq E(F_i)$ for $i=1,2,\ldots,r$, 
there exists a $2$-factorisation $\{F'_1,F'_2,\ldots,F'_r\}$ of $G$ such that 
for $i=1,2,\ldots,r$
\begin{itemize}
\item $E(F_i)\setminus E_i\subseteq E(F'_i)$;
\item if $u$ and $v$ are vertices in the same component of $F_i$, 
then $u$ and $v$ are in the same component of $F'_i$; and
\item each of the vertices in $V$ belongs to the same component in $F'_i$.
\end{itemize}
We say that $\{F'_1,F'_2,\ldots,F'_r\}$ is a {\em $2$-factorisation of $G$ obtained from $\{F_1,F_2,\ldots,F_r\}$ by applying the $V$-connector $\{E_1,E_2,\ldots,E_r\}$}.

Given a $2$-factorisation $\mathcal{F}=\{F_1,F_2,\ldots,F_r\}$ of a graph $G$, we say a subgraph $H$ of $G$ {\em induces a $V$-connector in $\mathcal{F}$} if 
$$\{E(F_1)\cap E(H),E(F_2)\cap E(H),\ldots,E(F_r)\cap E(H)\}$$ is a $V$-connector, and in this case we call $\{E(F_1)\cap E(H),E(F_2)\cap E(H),\ldots,E(F_r)\cap E(H)\}$ the {\em $V$-connector induced by $H$ in $\mathcal{F}$}.

\begin{lemma}\label{edgedisjointconnectors}
Suppose $\mathcal{F}$ is a $2$-factorisation of a graph $G$, and $H$ and $H'$ are edge-disjoint subgraphs of $G$ such that 
\begin{itemize}
\item [$(1)$] $H$ induces a $V$-connector in $\mathcal{F}$; and 
\item [$(2)$] $H'$ induces a $V'$-connector in $\mathcal{F}$.
\end{itemize}
If $\mathcal{F}'$ is any $2$-factorisation of $G$ obtained from $\mathcal{F}$ by applying the $V$-connector induced by $H$, then $H'$ induces a $V'$-connector in $\mathcal{F}'$.
\end{lemma}

\proof
Let $\mathcal{F}=\{F_1,F_2,\ldots,F_r\}$ and let 
$\mathcal{F'}=\{F'_1,F'_2,\ldots,F'_r\}$ be any $2$-factorisation of $G$ obtained from $\mathcal{F}$ by applying the $V$-connector $\{E_1,E_2,\ldots,E_r\}$ induced by $H$. 
Since $H$ and $H'$ are edge-disjoint, it follows from 
$E(F_i)\setminus E_i\subseteq E(F'_i)$ that $E(F_i)\cap E(H')=E(F'_i)\cap E(H')$ for $i=1,2,\ldots,r$. Thus, since $H'$ induces a $V'$-connector in $\mathcal{F}$, $H'$ also induces a 
$V'$-connector in $\mathcal{F'}$.
\qed

Note that Lemma \ref{edgedisjointconnectors} does not require $V$ and $V'$ to be disjoint.

\begin{lemma}\label{connectors}
Suppose $\alpha$, $\beta$, $u$, $v$, $w$, $x$, $u'$, $v'$ and $w'$ are distinct vertices. Then the following sets are $\{\alpha,\beta\}$-connectors:
\begin{itemize}
\item[$(0)$] $\{\{\alpha u, uw, \beta w\}, \{\alpha w, \beta u, uv\}\}$; 
\item[$(1)$] $\{\{\alpha u,uv,\beta w\},\{\alpha w,wv,\beta u\}\}$;
\item[$(2)$] $\{\{\alpha u,uv,\beta w\},\{\alpha w,wv,\beta u\},\{\alpha\beta\}\}$;
\item[$(3)$] $\{\{\alpha u,uv,\beta w, wx\},\{\alpha w,\beta u\},\{\alpha\beta,vw,wu,ux\}\}$;
\item[$(4)$] $\{\{\alpha u,uv,vw,w\beta\},\{\alpha v,\beta u, uw\},\{\alpha w, \beta v\}, \{\alpha x,x\beta\}\}$; 
\item[$(5)$] $\{\{\alpha x,vw,u\beta\},\{\alpha \beta, uv, wx\},\{\alpha u, ux, x\beta\}, \{\alpha v,uw, w\beta\}, \{\alpha w,xv, v\beta\}\}$; and
\item[$(6)$] $\{\{\alpha u,ux,xv,vw,\beta u',u'v',v'w'\},\{\alpha u',u'v,vw',\beta u,uv',v'w\}\}.$ 
\end{itemize}
Furthermore, the set 
\begin{itemize}
\item[$(7)$] $\{\{\alpha u,ux,xv,vw,\beta u',u'v',v'w'\},\{\alpha u',u'v,vw',\beta u,uv',v'w\},\{\alpha x, x\beta, uv\}\}$
\end{itemize}
 is an $\{\alpha,\beta,u\}$-connector.
\end{lemma}

\proof
For (0), let $E_1=\{\alpha u, uw, \beta w\}$, let $E_2=\{\alpha w, \beta u, uv\}$ and suppose $\{F_1,F_2\}$ is a $2$-factorisation of 
some $4$-regular graph $G$ such that $E_1\subseteq E(F_1)$ and $E_2\subseteq E(F_2)$.
We need to allocate the edges of $E_1\cup E_2$ to $E'_1$ and $E'_2$ such that if $F'_1=F_1-E_1+E'_1$ and $F'_2=F_2-E_2+E'_2$,
then $\{F'_1,F'_2\}$ is a $2$-factorisation of $G$ with the required properties.
Note that $\alpha$ and $\beta$ are currently in the same component of $F_1$. If $\alpha$ and $\beta$ are also in the same component of $F_2$ then we let $E'_1 = E_1$ and $E'_2 = E_2$. Otherwise, we let $E'_1= \{\alpha w, uw, \beta u\}$ and $E'_2=\{\alpha u, uv, \beta w\}$. This completes the proof for (0).
  
For (1), let $E_1=\{\alpha u,uv,\beta w\}$, let $E_2=\{\alpha w,wv,\beta u\}$ and suppose $\{F_1,F_2\}$ is a $2$-factorisation of 
some $4$-regular graph $G$ such that $E_1\subseteq E(F_1)$ and $E_2\subseteq E(F_2)$.
We need to allocate the edges of $E_1\cup E_2$ to $E'_1$ and $E'_2$ such that if $F'_1=F_1-E_1+E'_1$ and $F'_2=F_2-E_2+E'_2$,
then $\{F'_1,F'_2\}$ is a $2$-factorisation of $G$ with the required properties. 

Each of $E_1$ and $E_2$ induces a union of disjoint paths, and the allocation of edges to $E'_1$ and $E'_2$ depends only on 
how these paths are connected up (into cycles) by the paths in $F_1-E_1$ and $F_2-E_2$. There are three distinct ways that the two paths induced by $E_1$ can be connected up (namely, $\alpha$ to $v$ and $\beta$ to $w$, 
$\alpha$ to $\beta$ and $v$ to $w$, or $\alpha$ to $w$ and $v$ to $\beta$). Similarly, there are three distinct ways that the two paths 
induced by $E_2$ can be connected up. Thus, there are nine possibilities for which we need to find a suitable 
allocation of edges to $E'_1$ and $E'_2$. As listed below, one of the following three
allocations works for each of the nine possibilities. 
\begin{itemize} 
\item [(a)] $E'_1=E_1$ and $E'_2=E_2$; 
\item [(b)] $E'_1=\{\alpha w, \beta u, uv\}$ and $E'_2=\{\alpha u, \beta w, wv\}$;
\item [(c)] $E'_1=\{\alpha u, u \beta, vw\}$ and $E'_2=\{\alpha w, w \beta, uv\}$.
\end{itemize}

\noindent (a) is used when 
\begin{itemize}
\item in $F_1-E_1$ there are paths from $\alpha$ to $\beta$ and from $v$ to $w$, and in $F_2-E_2$ there are paths from 
$\alpha$ to $u$ and $\beta$ to $v$;
\item in $F_1-E_1$ there are paths from $\alpha$ to $\beta$ and from $v$ to $w$, and in $F_2-E_2$ there are paths from 
$\alpha$ to $\beta$ and $u$ to $v$;
\item in $F_1-E_1$ there are paths from $\alpha$ to $w$ and from $\beta$ to $v$, and in $F_2-E_2$ there are paths from 
$\alpha$ to $\beta$ and $u$ to $v$; and when
\item in $F_1-E_1$ there are paths from $\alpha$ to $w$ and from $\beta$ to $v$, and in $F_2-E_2$ there are paths from 
$\alpha$ to $u$ and $\beta$ to $v$. 
\end{itemize}

\noindent (b) is used when 
\begin{itemize}
\item in $F_1-E_1$ there are paths from $\alpha$ to $v$ and from $\beta$ to $w$, and in $F_2-E_2$ there are paths from 
$\alpha$ to $v$ and $\beta$ to $u$;
\item in $F_1-E_1$ there are paths from $\alpha$ to $v$ and from $\beta$ to $w$, and in $F_2-E_2$ there are paths from 
$\alpha$ to $\beta$ and $u$ to $v$; and when
\item in $F_1-E_1$ there are paths from $\alpha$ to $\beta$ and from $v$ to $w$, and in $F_2-E_2$ there are paths from 
$\alpha$ to $v$ and $\beta$ to $u$.
\end{itemize}

\noindent (c) is used when 
\begin{itemize}
\item in $F_1-E_1$ there are paths from $\alpha$ to $v$ and from $\beta$ to $w$, and in $F_2-E_2$ there are paths from 
$\alpha$ to $u$ and $\beta$ to $v$; and when
\item in $F_1-E_1$ there are paths from $\alpha$ to $w$ and from $\beta$ to $v$, and in $F_2-E_2$ there are paths from 
$\alpha$ to $v$ and $\beta$ to $u$. 
\end{itemize}
This completes the proof for (1).

Case (2) is an immediate consequence of (1). Using the same method as for (1), 
it is a routine (and somewhat tedious) exercise to check that each of sets listed in 
(3)-(6) is an $\{\alpha,\beta\}$-connector, and that the set listed in (7) is an $\{\alpha,\beta,u\}$-connector.
\qed

Suppose $H$ is a graph. We say a subset $U$ of $V(H)$ {\em links} $H$ if it contains at least one vertex from each connected component of $H$. Similarly, we say $U$ {\em links} a set $\mathcal{H}$ of graphs if it links each graph $H\in \mathcal{H}$.
Observe that a $2$-factorisation $\mathcal{F}$ of a graph $G$ is a Hamilton decomposition of $G$ if and only if $\{v\}$ links $\mathcal{F}$ for every $v\in V(G)$.

\begin{lemma}\label{intersects}
Suppose $\mathcal{F}$ is a factorisation of a graph $G$ in which each factor has no isolated vertices. If at most two of the factors in $\mathcal{F}$ are not connected, then there is a partition $\{U,U'\}$ of $V(G)$ such that each of $U$ and $U'$ links $\mathcal{F}$.
\end{lemma}

\proof
If fewer than two of the factors in $\mathcal{F}$ are not connected then the result is obvious. Suppose then that $\mathcal{F}$ contains precisely two non-connected factors, say $F$ and $F'$, and that $F$ has connected components $C_1,C_2,\ldots,C_r$ and $F'$ has connected components $C'_1,C'_2,\ldots,C'_t$.   
Define a bipartite (multi)graph $B$ with parts $X=\{C_1,C_2,\ldots,C_r\}$ and $Y=\{ C'_1,C'_2,\ldots,C'_t\}$ and edge set defined 
as follows. For each vertex $v\in V(G)$ join $C_i$ to $C_j'$ where $C_i$ is the connected component of $F$ that contains $v$ and $C'_j$ is the connected component of $F'$ that contains $v$. Since each connected component of $F$ or $F'$ has at least two vertices, $B$ has minimum degree at least $2$. 

It is well-known that the edges of any multigraph can be oriented so that the indegree of each vertex differs from its outdegree by at most 1.
Give the edges of $B$ such an orientation. Thus, since $B$ has minimum degree at least $2$, each vertex of $B$ has indegree at least 1, and outdegree at least 1. For each $v\in G$, if the edge of $B$ corresponding to $v$ is oriented from $X$ to $Y$, then we place $v$ in $U$. Otherwise, the edge of $B$ corresponding to $v$ is oriented from $Y$ to $X$ and we place $v$ in $U'$. It follows that each of $U$ and $U'$ links $\mathcal{F}$.
\qed

The following result allows us to obtain a Hamilton decomposition of a graph from a $2$-factorisation with some additional properties.

\begin{lemma}\label{manyrepairs}
Suppose there are pairwise edge-disjoint subgraphs $H_1,H_2,\ldots,H_m$ of $G$, subsets $V_1,V_2,\ldots,V_m$ of $V(G)$, and a $2$-factorisation $\mathcal{F}$ of $G$ such that
\begin{itemize}
\item [$(1)$] $\bigcup_{i=1}^m V_i$ links $\mathcal{F}$;
\item [$(2)$] $V_i\cap V_{i+1}\ne \emptyset$ for $i=1,2,\ldots,m-1$; and
\item [$(3)$] $H_i$ induces a $V_i$-connector in $\mathcal{F}$ for $i=1,2,\ldots,m$.
\end{itemize}
Then $G$ is Hamilton decomposable.
\end{lemma}

\proof Let $\mathcal{F}_0=\mathcal{F}$. By repeated application of Lemma \ref{edgedisjointconnectors}, 
we can inductively obtain a sequence $\mathcal{F}_0,\mathcal{F}_1,\ldots,\mathcal{F}_m$ of $2$-factorisations of $G$ by letting
$\mathcal{F}_i$ be a $2$-factorisation obtained from $\mathcal{F}_{i-1}$ by applying the $V_i$-connector induced by $H_i$ in $\mathcal{F}_{i-1}$. 
Observe that $\bigcup_{i=j+1}^m V_i$ links $\mathcal{F}_i$ for each $j=0,1,\ldots,m-1$, 
and in particular that $V_m$ links $\mathcal{F}_{m-1}$. It follows that $\mathcal{F}_m$ is a Hamilton decomposition of $G$.\qed

\section{Hamilton decomposable subgraphs of line graphs}\label{Section3}

We begin this section by introducing some notation that we will be using.
Let $G$ be a given regular graph such that $G$ contains a Hamilton cycle if $G$ has even degree, 
and $G$ contains a Hamiltonian $3$-factor if $G$ has odd degree.
Let $n=\lfloor\frac{\deg(G)}{2}\rfloor$ so that $G$ has degree $2n$ or $2n+1$, and when $G$ has degree $2n+1$ let $F$ be a $1$-factor in $G$ such that $G-F$ is Hamiltonian. 
Let $\{F_0,F_1,\ldots,F_{n-2},F_\infty\}$ 
be a $2$-factorisation of $G$ (if $G$ has degree $2n$) or $G-F$ (if $G$ has degree $2n+1$) such that $F_\infty$ is a Hamilton cycle
(such a $2$-factorisation exists by a well-known theorem first proved by Petersen \cite{Pet}, see \cite{Wes}).

Let $N_n$ denote the set $\{0,1,\ldots,n-2,\infty\}$ 
and for each $i\in N_n$, define $\S_i$ to be the set of vertices of $L(G)$ that correspond to edges in 
$F_i$. Further, when $G$ has degree $2n+1$ define $\S$ to be the set of vertices of $L(G)$ that correspond to edges in $F$. 
Thus, when $G$ has degree $2n$, $\{\S_i:i\in N_n\}$ partitions the vertex set of $L(G)$,
and when $G$ has degree $2n+1$, $\{\S_i:i\in N_n\}\cup\{\S\}$ partitions the vertex set of $L(G)$.

For each 
$i\in N_n$ let $\F_i$ be a directed graph obtained from $F_i$ by (arbitrarily) orienting its edges to form directed cycles. 
We call $\{\F_0,\F_1,\ldots,\F_{n-2},\F_\infty\}$ a {\em directed $2$-factorisation}.
The following definitions are made in the 
context of an existing directed $2$-factorisation $\{\F_0,\F_1,\ldots,\F_{n-2},\F_\infty\}$
of $G$ or $G-F$ for some given graph $G$.

For each $v\in V(G)$ and each $i\in N_n$, 
assign the label $a^v_i$ to the vertex of $L(G)$
whose corresponding edge in $G$ is the edge directed into 
$v$ in $\F_i$, and 
assign the label $b^v_i$ to the vertex of $L(G)$
whose corresponding edge in $G$ is the edge directed out of  
$v$ in $\F_i$.
Further, if $G$ has degree $2n+1$, then 
assign the label $c^v$ to the vertex of $L(G)$ whose corresponding edge in $G$ is the edge of $F$ incident with $v$.
Thus, each vertex of $L(G)$ is assigned two labels. Each vertex of $L(G)$ corresponding to an edge $uv\in F_i$ is assigned 
labels $a^v_i$ and $b^u_i$ where $uv$ is oriented from $u$ to $v$ in $\F_i$, 
and, in the case $G$ has degree $2n+1$, each vertex of $L(G)$ corresponding to an edge $uv$ of $F$ is assigned labels $c^u$ and $c^v$. 

If $G$ has degree $2n$, then for each $v\in V(G)$, 
the $2n$ vertices with labels in the set 
$$\{a^v_0,a^v_1,\ldots,a^v_{n-2},a^v_\infty,b^v_0,b^v_1,\ldots,b^v_{n-2},b^v_\infty\}$$
induce a complete subgraph in $L(G)$.
Similarly, if $G$ has degree $2n+1$, then for each $v\in V(G)$, 
the $2n+1$ vertices with labels in the set 
$$\{a^v_0,a^v_1,\ldots,a^v_{n-2},a^v_\infty,b^v_0,b^v_1,\ldots,b^v_{n-2},b^v_\infty,c^v\}$$
induce a complete subgraph in $L(G)$.
In either case, we denote this complete subgraph by $L(G)_v$, and note that 
$\{L(G)_v:v\in V(G)\}$ is a decomposition of
$L(G)$ into complete subgraphs. 

Let $A_n=\{a_0,a_1,\ldots,a_{n-2},a_\infty\}$ and let
$B_n=\{b_0,b_1,\ldots,b_{n-2},b_\infty\}$. 
For each $v\in V(G)$, there is an obvious bijection 
$$
\begin{array}{ll}
\sigma_v:A_n\cup B_n\rightarrow V(L(G)_v)&\mbox{ if $G$ has degree $2n$};\\
\sigma_v:A_n\cup B_n\cup\{c\}\rightarrow V(L(G)_v)&\mbox{ if $G$ has degree $2n+1$;}
\end{array}
$$ 
given by $\sigma_v(a_i)=a^v_i$ and $\sigma_v(b_i)=b^v_i$ for each $i\in N_n$, and $\sigma_v(c)=c^v$.
For any subgraph $H$ of $K_{A_n\cup B_n}$ (if $G$ has degree $2n$) or $K_{A_n\cup B_n\cup\{c\}}$ (if $G$ has degree $2n+1$), and 
for each $v\in V(G)$, we define $\sigma_v(H)$ to be the graph with $V(\sigma_v(H))=V(L(G)_v)$ and 
$\sigma_v(x)\sigma_v(y)\in E(\sigma_v(H))$ if and only if $xy\in E(H)$. 

We call a subgraph $H$ of $K_{A_n\cup B_n}$ 
a {\em Hamilton fragment} if for any given Hamiltonian
$2n$-regular graph $G$ and any directed $2$-factorisation 
$\{\F_0,\F_1,\ldots,\F_{n-2},\F_\infty\}$
of $G$ where $\F_\infty$ is a Hamilton cycle, the subgraph
$$\bigcup_{v\in V(G)}\sigma_v(H)$$
of $L(G)$ has a Hamilton decomposition. 
Similarly, we 
call a subgraph $H$ of $K_{A_n\cup B_n\cup\{c\}}$ 
a {\em Hamilton fragment} if for any given 
$2n+1$-regular graph $G$ and any directed $2$-factorisation 
$\{\F_0,\F_1,\ldots,\F_{n-2},\F_\infty\}$
of $G-F$ where $F$ is a $1$-factor and $\F_\infty$ is a Hamilton cycle, the subgraph
$$\bigcup_{v\in V(G)}\sigma_v(H)$$
of $L(G)$ has a Hamilton decomposition.

Lemma \ref{combineHamFrags} below follows immediately from the definition of Hamilton fragment, because if $\{H_1, H_2, \dots, H_t\}$ is a decomposition of $K_{A_n\cup B_n}$ into Hamilton fragments, then $\{\cup_{v\in V(G)}\sigma_v(H_1), \cup_{v\in V(G)}\sigma_v(H_2), \dots, \cup_{v\in V(G)}\sigma_v(H_t)\}$ is a decomposition of $L(G)$ into Hamilton decomposable factors whenever $G$ is a Hamiltonian $2n$-regular graph; and similarly if $K_{A_n\cup B_n\cup \{c\}}$ has a decomposition into Hamilton fragments and $G$ is any $(2n+1)$-regular graph with a Hamiltonian 3-factor.
 
\begin{lemma}\label{combineHamFrags}
If $K_{A_n\cup B_n}$ can be decomposed into Hamilton fragments,
then the line graph of any Hamiltonian $2n$-regular graph has a Hamilton decomposition, and if  
$K_{A_n\cup B_n\cup \{c\}}$ 
can be decomposed into Hamilton fragments,
then the line graph of any 
$(2n+1)$-regular graph with a Hamiltonian 3-factor 
has a Hamilton decomposition.
\end{lemma}

For any subgraph $H$ of $K_{A_n\cup B_n}$ (if $G$ has degree $2n$) or $K_{A_n\cup B_n\cup\{c\}}$ (if $G$ has degree $2n+1$), 
we define $^{A}H$ to be the (multi)graph obtained from $H$ by amalgamating vertices $a_i$ and $b_i$ into a single vertex labelled $i$ for each $i\in N_n$. Thus, $^AH$ has vertex set $N_n$ (if $G$ has degree $2n$) or $N_n\cup\{c\}$ (if $G$ has degree $2n+1$), 
and has
an edge with endpoints $i$ and $j$ for each edge of $H$ whose endpoints have subscripts $i$ and $j$. Further, when $G$ has degree $2n+1$, 
$^AH$ has 
an edge with endpoints $c$ and $i$ for each edge of $H$ whose endpoints are $c$ and a vertex with subscript $i$. 
Note that an edge joining $a_i$ to $b_i$ in $H$ results in a loop on vertex $i$ in $^AH$.

\subsection{Regular graphs of even degree}\label{SectionRegulargraphsofevendegree}

\begin{lemma}\label{linking2factors}
Suppose $G$ is a $2n$-regular graph, $\{\F_0,\F_1,\ldots,\F_{n-2},\F_\infty\}$ is a directed $2$-factorisation of $G$ and $X$ is a subgraph of $K_{A_n\cup B_n}$ such that 
$^{A}X$ is an $n$-cycle.
Then the graph $$J=\bigcup_{v\in V(G)}\sigma_v(X)$$
is a $2$-factor of $L(G)$ and $\S_i$ links $J$ for each $i\in N_n$.
\end{lemma}

\proof
Let $uv\in E(G)$, let $F_i$ be the $2$-factor containing $uv$, and let $uv$ be oriented from $u$ to $v$ in $\F_i$. 
Then $\deg_J(uv)=\deg_{\sigma_v(X)}(a_i)+\deg_{\sigma_u(X)}(b_i)=\deg_X(a_i)+\deg_X(b_i)=\deg_{^AX}(i)=2$. Thus, $J$ is a $2$-factor of $L(G)$.
We now show that $\S_i$ links $J$ for each $i\in N_n$. 
Each edge of $J$ if of the form $\sigma_v(x)\sigma_v(y)$ where $xy$ is an edge of $X$, and so if an edge of $J$ has its endpoints 
in $\S_i$ and $\S_j$, then $ij$ is an edge of $^AX$. It thus follows from the fact that $^AX$ is an $n$-cycle that each component of 
$J$ contains at least one vertex from $\S_i$ for each $i\in N_n$. That is, $\S_i$ links $J$ for each $i\in N_n$. 
\qed

\begin{lemma}\label{linking2factors1}
Suppose $G$ is a $2n$-regular graph, $\{\F_0,\F_1,\ldots,\F_{n-2},\F_\infty\}$ is a directed $2$-factorisation of $G$, $T$ is a non-empty subset of $N_n$, $U$ is a subset of $V(G)$ that links $\{\F_j : j\in T\}$, and $X$ and $X'$ are subgraphs of $K_{A_n\cup B_n}$ such that 
\begin{itemize}
\item [$(1)$] $\deg_{X}(a_i)=\deg_{X'}(a_i)$ and $\deg_{X}(b_i)=\deg_{X'}(b_i)$ for each $i\in N_n$; 
\item [$(2)$] $\deg_X(a_i)=1$ if and only if $i\in T$ and $\deg_X(b_i)=1$ if and only if $i\in T$;
\item [$(3)$] $^{A}X$ is an $n$-cycle; and
\item [$(4)$] both $X'$ and $^{A}X'$ have $|T|$ connected components.
\end{itemize}
Then the graph $$J=\left(\bigcup_{v\in U}\sigma_v(X)\right)\cup \left(\bigcup_{v\in V(G)\setminus U}\sigma_v(X')\right)$$
is a $2$-factor of $L(G)$ and $\S_i$ links $J$ for each $i\in N_n$.
\end{lemma}

\proof
Let $uv\in E(G)$, let $F_i$ be the $2$-factor containing $uv$, and let $uv$ be oriented from $u$ to $v$ in $\F_i$. 
Then $\deg_J(uv)=\deg_{Y_u}(b_i)+\deg_{Y_v}(a_i)$ where $Y_u=X$ if $u\in U$, $Y_u=X'$ if $u\notin U$, $Y_v=X$ if $v\in U$ and $Y_v=X'$ if $v\notin U$.
But $\deg_X(b_i)=\deg_{X'}(b_i)$ and $\deg_X(a_i)=\deg_{X'}(a_i)$. So $\deg_J(uv)=\deg_X(b_i)+\deg_{X}(a_i)=\deg_{^AX}(i)=2$.
Thus, $J$ is a $2$-factor of $L(G)$.

We now show that $\S_i$ links $J$ for each $i\in N_n$.
Since $^AX$ is an $n$-cycle, it follows from $(1)$ that $^AX'$ is a $2$-regular graph with vertex set $N_n$. 
Thus, each connected component of $X$ is a path and each connected component of $X'$ is a path or a cycle. 
But by $(4)$, $X'$ and $^{A}X'$ have the same number of components and it follows from this that any maximal path in 
$X'$ has endpoints $a_k$ and $b_k$ for some $k$. Thus, by $(2)$ (and $(1)$) $X'$ contains $|T|$ vertex disjoint paths; one from 
$a_k$ to $b_k$ for each $k\in T$. Also, by $(4)$, there are no other connected components in $X'$. In particular, there are no cycles
in $X'$.

Let $Z$ be a connected component of $J$.
Then $Z$ is a cycle and consists of a sequence $P_1,P_2,\ldots,P_r$ of paths 
where for $j=1,2,\ldots,r$ we have 
$P_j$ is a maximal path in $\sigma_{v_j}(Y_j)$, 
$v_j\in V(G)$, 
$Y_j=X$ if $v_j\in U$ and $Y_j=X'$ if $v_j\in V(G)\setminus U$.

Since any path in $X'$ has endpoints $a_k$ and $b_k$ for some $k\in T$, 
if $v_j\in V(G)\setminus U$ for $j=1,2,\ldots,r$, then $Y_j=X'$ for $j=1,2,\ldots,r$ and 
$(v_1,v_2,\ldots,v_r)$ is a connected component of $\F_k$ for some $k\in T$.
This contradicts the fact that $U$ links $\{\F_j : j\in T\}$. Thus, there exists $j^*\in\{1,2,\ldots,r\}$ such that 
$v_{j^*}\in U$ and $Y_{j^*}=X$.
Since any path in $X'$ has endpoints $a_k$ and $b_k$ for some $k\in T$, and since $^AX$ is an $n$-cycle, it thus 
follows that $Z$ contains at least one vertex of $\S_i$ for each $i\in N_n$. That is, $\S_i$ links $J$ for each $i\in N_n$.
\qed

\begin{lemma}\label{linking2factors2}
Suppose $G$ is a $2n$-regular graph, $\{\F_0,\F_1,\ldots,\F_{n-2},\F_\infty\}$ is a directed $2$-factorisation of $G$, $T$ is a non-empty proper subset of $N_n$, 
$U$ is a subset of $V(G)$ such that $U$ links $\{\F_j : j\in T\}$, and $X$, $C$ and $C'$ are subgraphs of $K_{A_n\cup B_n}$ such that 
\begin{itemize}
\item [$(1)$] $\deg_{C}(a_i)=\deg_{C'}(a_i)$ and $\deg_{C}(b_i)=\deg_{C'}(b_i)$ for each $i\in N_n$; 
\item [$(2)$] $C\cup C'$ is a $2(|T|+1)$-cycle with $a_jb_j\in E(C')$ for each $j\in T$; and
\item [$(3)$] $^{A}(X\cup C)$ is an $n$-cycle.
\end{itemize}
Then the graph $$J=\left(\bigcup_{v\in U}\sigma_v(X\cup C)\right)\cup \left(\bigcup_{v\in V(G)\setminus U}\sigma_v(X\cup C')\right)$$
is a $2$-factor of $L(G)$ and $\S_i$ links $J$ for each $i\in N_n\setminus T$.
\end{lemma}

\proof
Let $uv\in E(G)$, let $F_i$ be the $2$-factor containing $uv$, and let $uv$ be oriented from $u$ to $v$ in $\F_i$. 
Then $\deg_J(uv)=\deg_{Y_u}(b_i)+\deg_{Y_v}(a_i)$ where $Y_u=X\cup C$ if $u\in U$, $Y_u=X\cup C'$ if $u\notin U$, $Y_v=X\cup C$ if $v\in U$ and $Y_v=X\cup C'$ if $v\notin U$.
But by $(1)$ we have $\deg_{X\cup C}(b_i)=\deg_{X\cup C'}(b_i)$ and $\deg_{X\cup C}(a_i)=\deg_{X\cup C'}(a_i)$. 
So $\deg_J(uv)=\deg_{X\cup C}(b_i)+\deg_{X\cup C}(a_i)=\deg_{^A(X\cup C)}(i)=2$.
Thus, $J$ is a $2$-factor of $L(G)$.

We now show that $\S_i$ links $J$ for each $i\in N_n\setminus T$.
It is a consequence of $(1)$ and the fact that $C\cup C'$ is a $2(|T|+1)$-cycle that each of $C$ and $C'$ is a matching with $|T|+1$ edges. Also, since $a_jb_j\in E(C')$ for each $j\in T$, the edge set of $C'$ consists of these $|T|$ edges and one other edge $e$. 
Let the subscripts of the endpoints of $e$ be $j'$ and $j''$ (so without loss of generality 
$e=a_{j'}a_{j''}$, $a_{j'}b_{j''}$ or $b_{j'}b_{j''}$).

The special case where $|T|=n-1$ will be dealt with later. So for now assume that $|T|\leq n-2$. 
If $j'=j''$, then $^AC$ is a $2$-regular graph with vertex set $T\cup\{j'\}$. But this contradicts the fact that 
$^A(X\cup C)$ is an $n$-cycle, and so we have $j'\neq j''$. This means that $^AC$ is a path with vertex set 
$T\cup\{j',j''\}$ and endpoints $j'$ and $j''$, 
and hence that $^AX$ is a path with vertex set $N_n\setminus T$ and
endpoints $j'$ and $j''$. 
Let $^AX$ be the path $[j',j_1,j_2,\ldots,j_s,j'']$ where $s=|n-(|T|+2)|$ (if $|T|=n-2$, then $^AX$ is the path $[j',j''$]). So $\{j_1,j_2,\ldots,j_s\}=N_n\setminus (T\cup\{j',j''\})$.

Let $Z$ be a connected component of $J$. 
It follows from the observations made in the preceding paragraph that if $Z$ has a vertex in $\S_i$ where $i\in N_n\setminus T$, 
then 
$Z$ contains a path $x_{j'},x_1,\ldots,x_s,x_{j''}$ 
where $x_{j'}\in \S_{j'}$, $x_{j''}\in \S_{j''}$ and $x_i\in \S_{j_i}$ for $i=1,2,\ldots,s$. 
Thus, $Z$ has at least one vertex in $\S_i$ for each $i\in N_n\setminus T$. 
On the other hand, if every vertex of $Z$ is in $\cup_{i\in T}\S_i$, then it follows (again from the observations made in the preceding paragraph) that $V(Z)\subseteq \S_i$ for some $i\in T$. Moreover, every edge of $Z$ is of the form $a_i^vb_i^v$ where $v\in V(G)\setminus U$. 
But this implies that there is a cycle in $\F_i$ that contains no vertex of $U$, contradicting the fact that $U$ links $\{\F_j:j\in T\}$.

Thus, we are left with the special case where $|T|=n-1$ which was mentioned earlier. In this case we have $E(C')=\{a_ib_i:i\in N_n\}$
and by the same argument as in the previous paragraph we reach a contradiction if we assume every vertex of a component of $J$ 
is in $\cup_{i\in T}\S_i$. Thus, any component of $J$ has a vertex in $\S_i$ where $i$ is the unique element of $N_n\setminus T$.  
This completes the proof that $\S_i$ links $J$ for each $i\in N_n\setminus T$. 
\qed

The next few lemmas establish some sufficient conditions for $H$ to be a Hamilton fragment.

\begin{lemma}\label{4regular_withHam}
Suppose $r\geq 2$ is a positive integer, $s$ and $t$ are distinct elements of $\{0,1,\ldots,n-2\}$, 
$H$ is a subgraph of $K_{A_n\cup B_n}$, and $H$ admits a decomposition $\{X_1,X_2,\ldots,X_r\}$ such that:
\begin{itemize}
\item [$(1)$] $\deg_{X_1}(a_i)=\deg_{X_2}(a_i)$ and $\deg_{X_1}(b_i)=\deg_{X_2}(b_i)$ for each $i\in N_n$;
\item [$(2)$] $\deg_{X_1}(u)=1$ if and only if $u\in\{a_s,b_s,a_t,b_t,a_\infty,b_\infty\}$;
\item [$(3)$] $^{A}X_1$ is an $n$-cycle;
\item [$(4)$] both $X_2$ and $^{A}X_2$ have three connected components;
\item [$(5)$] if $r\geq 3$ then $^{A}X_k$ is an $n$-cycle for $k=3,4,\ldots,r$; and 
\item [$(6)$] $\{E(X_1),E(X_2),\ldots,E(X_r)\}$ is an $\{a_{\infty},b_{\infty}\}$-connector.
\end{itemize} 
Then $H$ is a Hamilton fragment.
\end{lemma}

\proof
Let $r$, $s$, $t$, $H$ and $X_1,X_2,\ldots,X_r$ satisfy the conditions of the lemma, let $G$ be a $2n$-regular graph,  and let $\{\F_0,\F_1,\ldots,\F_{n-2},\F_\infty\}$ be a directed $2$-factorisation of $G$ where $\F_\infty$ is a Hamilton cycle, say $\F_\infty
=(v_1,v_2,\ldots,v_m)$. Our aim is to show that the subgraph $$L=\bigcup_{v\in V(G)}\sigma_v(H)$$ of $L(G)$ decomposes into Hamilton cycles,
and we do this by applying Lemma \ref{manyrepairs} to $L$ with $H_i=\sigma_{v_i}(H)$ and $V_i=\{a_\infty^{v_i},b_\infty^{v_i}\}$ for each $i=1,2,\ldots,m$.
To this end, observe that $\sigma_{v_1}(H),\sigma_{v_2}(H),\ldots,\sigma_{v_m}(H)$ are edge-disjoint subgraphs of $L$,
and since $b_\infty^{v_i}=a_\infty^{v_{i+1}}$ for $i=1,2,\ldots,m-1$ it follows that $V_i\cap V_{i+1}\ne \emptyset$ for $i=1,2,\ldots,m-1$ as required. It remains to show there is a $2$-factorisation $\{J_1,J_2,\ldots,J_r\}$ of $L$ such that 
\begin{itemize}
\item $\{a_\infty^{v_1},a_\infty^{v_2},\ldots,a_\infty^{v_m}\}$ links $\{J_1,J_2,\ldots,J_r\}$; and
\item $\sigma_{v_i}(H)$ induces an $\{a_\infty^{v_{i}},b_\infty^{v_i}\}$-connector in $\{J_1,J_2,\ldots,J_r\}$, for $i=1,2,\ldots,m$.
\end{itemize} 
Let $\{U,U'\}$ be a partition of $V(G)$ such that both $U$ and $U'$ link  $\{\F_s,\F_t,\F_\infty\}$ (such a partition exists by Lemma \ref{intersects}), and let $\{J_1,J_2,\ldots,J_r\}$ be the decomposition of $L$ defined by 
\begin{itemize}
\item $J_1=(\bigcup_{v\in U}\sigma_v(X_1))\cup (\bigcup_{v\in U'}\sigma_v(X_2))$;
\item $J_2= (\bigcup_{v\in U}\sigma_v(X_2))\cup (\bigcup_{v\in U'}\sigma_v(X_1))$; and
\item if $r\geq 3$, then $J_k=(\bigcup_{v\in V(G)}\sigma_v(X_k))$ for $k=3,4,\ldots,r$.
\end{itemize}
It follows easily from Lemma \ref{linking2factors1} that both $J_1$ and $J_2$ are $2$-factors of $L$, and that $\{a_\infty^{v_1},a_\infty^{v_2},\ldots,a_\infty^{v_m}\}$ links $\{J_1,J_2\}$ as required. Similarly, if $r\geq 3$, it follows easily from Lemma \ref{linking2factors} that each of $J_3,J_4,\ldots,J_r$ is in fact a $2$-factor of $L$, and that $\{a_\infty^{v_1},a_\infty^{v_2},\ldots,a_\infty^{v_m}\}$ links $\{J_3,J_4,\ldots,J_r\}$ as required.

Finally, observe that 
$$\{E(J_k)\cap E(\sigma_{v_i}(H)): k=1,2,\ldots,r\}= \{E(\sigma_{v_i}(X_1)),E(\sigma_{v_i}(X_2)),\ldots,E(\sigma_{v_i}(X_r))\}$$
for $i=1,2,\ldots,m$. Then, since $\{E(X_1),E(X_2),\ldots,E(X_r)\}$ is an $\{a_{\infty},b_{\infty}\}$-connector, 
it follows that $\sigma_{v_i}(H)$ induces an $\{a_\infty^{v_{i}},b_\infty^{v_i}\}$-connector in $\{J_1,J_2,\ldots,J_r\}$, for $i=1,2,\ldots,m$ as required.
\qed

\begin{lemma}\label{6regular_2}
Suppose $s$ and $t$ are distinct elements of $\{0,1,\ldots,n-2\}$, $H$ is a subgraph of $K_{A_n\cup B_n}$, and $H$ admits decompositions $\{X_1,X_2,X_3\}$ and $\{X_1',X_2',X_3\}$ such that:
\begin{itemize}
\item [$(1)$] $\deg_{X_1}(u)=\deg_{X'_2}(u)$ for each $u\in A_n\cup B_n$;
\item [$(2)$] $\deg_{X_1}(u)=1$ if and only if $u\in\{a_s,b_s,a_\infty,b_\infty\}$;
\item [$(3)$] $\deg_{X_1'}(u)=\deg_{X_2}(u)$ for each $u\in A_n\cup B_n$;
\item [$(4)$] $\deg_{X_1'}(u)=1$ if and only if $u\in\{a_t,b_t,a_\infty,b_\infty\}$;
\item [$(5)$] both $^{A}X_1$ and $^{A}X_1'$ are $n$-cycles;
\item [$(6)$] each of $X_2$, $X_2'$, $^{A}X_2$ and $^{A}X_2'$ has two connected components;
\item [$(7)$] $^{A}X_3$ is an $n$-cycle; and 
\item [$(8)$] both $\{E(X_1),E(X_2),E(X_3)\}$ and $\{E(X_1'),E(X_2'),E(X_3)\}$ are $\{a_{\infty},b_{\infty}\}$-connectors.
\end{itemize} 
Then $H$ is a Hamilton fragment.
\end{lemma}

\proof
Let $s$, $t$, $H$, $X_1$, $X_2$, $X_3$, $X_1'$ and $X_2'$ satisfy the conditions of the lemma, let $G$ be a $2n$-regular graph,  and let $\{\F_0,\F_1,\ldots,\F_{n-2},\F_\infty\}$ be a directed $2$-factorisation of $G$ where $\F_\infty$ is a Hamilton cycle, say $\F_\infty=(v_1,v_2,\ldots,v_m)$. Our aim is to show that the subgraph $$L=\bigcup_{v\in V(G)}\sigma_v(H)$$ of $L(G)$ decomposes into Hamilton cycles, and we do this by applying Lemma \ref{manyrepairs} with $H_i=\sigma_{v_i}(H)$ and $V_i=\{a_\infty^{v_i},b_\infty^{v_i}\}$ for each $i=1,2,\ldots,m$.
To this end, observe that
\linebreak $\sigma_{v_1}(H),\sigma_{v_2}(H),\ldots,\sigma_{v_m}(H)$ are edge-disjoint subgraphs of $L$,
and since $b_\infty^{v_i}=a_\infty^{v_{i+1}}$ for $i=1,2,\ldots,m-1$ it follows that $V_i\cap V_{i+1}\ne \emptyset$ for $i=1,2,\ldots,m-1$ as required. It remains to show there is a $2$-factorisation $\{J_1,J_2,J_3\}$ of $L$ such that 
\begin{itemize}
\item $\{a_\infty^{v_0},a_\infty^{v_1},\ldots,a_\infty^{v_m}\}$ links $\{J_1,J_2,J_3\}$; and
\item $\sigma_{v_i}(H)$ induces an $\{a_\infty^{v_{i-1}},a_\infty^{v_i}\}$-connector in $\{J_1,J_2,J_3\}$, for $i=1,2,\ldots,m$.
\end{itemize} 
Let $\{U,U'\}$ be a partition of $V(G)$ such that both $U$ and $U'$ link  $\{\F_s,\F_t,\F_\infty\}$ (such a partition exists by Lemma \ref{intersects}), and let $\{J_1,J_2,J_3\}$ be the decomposition of $L$ defined by 
\begin{itemize}
\item $J_1=(\bigcup_{v\in U}\sigma_v(X_1))\cup (\bigcup_{v\in U'}\sigma_v(X_2'))$;
\item $J_2= (\bigcup_{v\in U}\sigma_v(X_2))\cup (\bigcup_{v\in U'}\sigma_v(X_1'))$; and
\item $J_3= (\bigcup_{v\in V(G)}\sigma_v(X_3))$.
\end{itemize}
It follows easily from Lemma \ref{linking2factors1} that both $J_1$ and $J_2$ are $2$-factors of $L$, and that $\{a_\infty^{v_1},a_\infty^{v_2},\ldots,a_\infty^{v_m}\}$ links $\{J_1,J_2\}$ as required. Similarly, it follows easily from Lemma \ref{linking2factors} that $J_3$ is a $2$-factor of $L$, and that $\{a_\infty^{v_1},a_\infty^{v_2},\ldots,a_\infty^{v_m}\}$ links $\{J_3\}$ as required.

Finally, observe that for each $i=1,2,\ldots,m$ we have that $$\{E(J_k)\cap E(\sigma_{v_i}(H)) : k=1,2,3\}=\{E(\sigma_{v_i}(X_1)),E(\sigma_{v_i}(X_2)),E(\sigma_{v_i}(X_3))\}$$ when $v_i\in U$, and 
$$\{E(J_k)\cap E(\sigma_{v_i}(H)) : k=1,2,3\}=\{E(\sigma_{v_i}(X_1')),E(\sigma_{v_i}(X_2')),E(\sigma_{v_i}(X_3))\}$$ when $v_i\in U'$.
 Then, since both $\{E(X_1),E(X_2),E(X_3)\}$ and $\{E(X_1'),E(X_2'),E(X_3)\}$ are $\{a_{\infty},b_{\infty}\}$-connectors, 
it follows that $\sigma_{v_i}(H)$ induces an $\{a_\infty^{v_{i}},b_\infty^{v_i}\}$-connector in $\{J_1,J_2,J_3\}$, for $i=1,2,\ldots,m$ as required.
\qed

A subgraph $H$ of $K_{A_n\cup B_n}$ is said to be an {\em $R$-adjustable Hamilton fragment} if $R$ is a subgraph of $H$ such that $H-R+Q$ is a Hamilton fragment for any subgraph $Q$ of $K_{A_n\cup B_n}-(H-R)$ satisfying $^{A}Q=\prescript{A}{}{R}$.
The following lemmas makes use of this concept while establishing some sufficient conditions for $H$ to be a Hamilton fragment in cases where $\bigcup_{v\in V(G)}\sigma_v(H)$ is $4$-regular. 

\begin{lemma}\label{4regular_easy}
Suppose $H$ is a subgraph of $K_{A_n\cup B_n}$ and $H$ admits a decomposition $\{S,S',R,R'\}$ such that:
\begin{itemize}
\item [$(1)$] $^{A}(S\cup R)$ is an $n$-cycle;
\item [$(2)$] $^{A}(S' \cup R')$ is an $n$-cycle; and
\item [$(3)$] $\{E(S),E(S')\}$ is an $\{a_{\infty},b_{\infty}\}$-connector.
\end{itemize} 
Then $H$ is an $(R\cup R')$-adjustable Hamilton fragment.
\end{lemma}

\proof
Suppose $\{S,S',R,R'\}$ is a decomposition of $H$ that satisfies the conditions of the lemma. It is easy to see that if $Q$ and $Q'$ are edge-disjoint subgraphs of $K_{A_n\cup B_n}-(S\cup S')$ satisfying $^{A}Q=\prescript{A}{}R$ and $^{A}Q'=\prescript{A}{}R'$, then $H-(R\cup R')+(Q\cup Q')$ has a decomposition, namely $\{S,S',Q,Q'\}$, which also satisfies the conditions of the lemma. Thus we need only show $H$ is a Hamilton fragment and it follows immediately that it is $(R\cup R')$-adjustable. 

Let $G$ be a $2n$-regular graph and let $\{\F_0,\F_1,\ldots,\F_{n-2},\F_\infty\}$ be a directed $2$-factorisation of $G$ where $\F_\infty$ is a Hamilton cycle, say $\F_\infty=(v_1,v_2,\ldots,v_m)$. Our aim is to show that the subgraph $$L=\bigcup_{v\in V(G)}\sigma_v(H)$$
of $L(G)$ decomposes into Hamilton cycles, and we do this by applying Lemma \ref{manyrepairs} with 
with $H_i=\sigma_{v_i}(H)$ and $V_i=\{a_\infty^{v_i},b_\infty^{v_i}\}$ for each $i=1,2,\ldots,m$.
To this end, observe that $\sigma_{v_1}(H),\sigma_{v_2}(H),\ldots,\sigma_{v_m}(H)$ are edge-disjoint subgraphs of $L$,
and since $b_\infty^{v_i}=a_\infty^{v_{i+1}}$ for $i=1,2,\ldots,m-1$ it follows that $V_i\cap V_{i+1}\ne \emptyset$ for $i=1,2,\ldots,m-1$ as required. It remains to show there is a $2$-factorisation $\{J_1,J_2\}$ of $L$ such that 
\begin{itemize}
\item $\{a_\infty^{v_1},a_\infty^{v_2},\ldots,a_\infty^{v_m}\}$ links $\{J_1,J_2\}$; and
\item $\sigma_{v_i}(H)$ induces an $\{a_\infty^{v_{i}},b_\infty^{v_i}\}$-connector in $\{J_1,J_2\}$, for $i=1,2,\ldots,m$.
\end{itemize} 

Let $\{J_1,J_2\}$ be the decomposition of $L$ defined by 
\begin{itemize}
\item $J_1=\bigcup_{v\in V(G)}\sigma_v(S\cup R)$; and
\item $J_2= \bigcup_{v\in V(G)}\sigma_v(S'\cup R')$.
\end{itemize}
Since both $^{A}(S\cup R)$ and $^{A}(S'\cup R')$ are $n$-cycles, it follows from Lemma \ref{linking2factors} that both $J_1$ and $J_2$ are $2$-factors of $L$ and that  
$\{a_\infty^{v_0},a_\infty^{v_1},\ldots,a_\infty^{v_m}\}$ links $\{J_1,J_2\}$ as required.

Finally, observe that 
$E(\sigma_{v_i}(S))\subseteq E(J_1)\cap E(\sigma_{v_i}(H))$ and $E(\sigma_{v_i}(S'))\subseteq E(J_2)\cap E(\sigma_{v_i}(H))$ for $i=1,2,\ldots,m$. Then, since $\{E(S),E(S')\}$ is an $\{a_{\infty},b_{\infty}\}$-connector, 
it follows that $\sigma_{v_i}(H)$ induces an $\{a_\infty^{v_{i}},b_\infty^{v_i}\}$-connector in $\{J_1,J_2\}$, for $i=1,2,\ldots,m$ as required.
\qed

\begin{lemma}\label{4regular}
Suppose $T$ is a subset of $\{0,1,\ldots,n-2\}$ with $|T|\in\{1,2\}$, $H$ is a subgraph of $K_{A_n\cup B_n}$, and $H$ admits a decomposition $\{S,S',R,R',C,C'\}$ such that:
\begin{itemize}
\item [$(1)$] $\deg_C(u)=\deg_{C'}(u)$ for each $u\in A_n\cup B_n$;
\item [$(2)$] $C\cup C'$ is a $2(|T|+1)$-cycle with $a_j b_j\in E(C')$ for each $j\in T$; 
\item [$(3)$] $^{A}(S\cup R\cup C)$ is an $n$-cycle;
\item [$(4)$] $^{A}(S' \cup R'\cup C')$ is the vertex disjoint union of an $(n-|T|)$-cycle and a loop on each vertex $j\in T$; and
\item [$(5)$] $\{E(S),E(S')\}$ is an $\{a_{\infty},b_{\infty}\}$-connector.
\end{itemize} 
Then $H$ is an $(R\cup R')$-adjustable Hamilton fragment.
\end{lemma}

\proof
Suppose $T$ is a subset of $\{0,1,\ldots,n-2\}$ with $|T|\in\{1,2\}$, and \linebreak $\{S,S',R,R',C,C'\}$ is a decomposition of $H$ that satisfies the conditions of the lemma. It is easy to see that if $Q$ and $Q'$ are edge-disjoint subgraphs of $K_{A_n\cup B_n}-(S\cup S'\cup C\cup C')$ satisfying $^{A}Q=\prescript{A}{}R$ and $^{A}Q'=\prescript{A}{}R'$, then $H-(R\cup R')+(Q\cup Q')$ has a decomposition, namely $\{S,S',Q,Q',C,C'\}$, which also satisfies the conditions of the lemma. Thus we need only show $H$ is a Hamilton fragment and it follows immediately that it is $(R\cup R')$-adjustable. 

Let $G$ be a $2n$-regular graph and let $\{\F_0,\F_1,\ldots,\F_{n-2},\F_\infty\}$ be a directed $2$-factorisation of $G$ where $\F_\infty$ is a Hamilton cycle, say $\F_\infty=(v_1,v_2,\ldots,v_m)$. Our aim is to show that the subgraph $$L=\bigcup_{v\in V(G)}\sigma_v(H)$$
of $L(G)$ decomposes into Hamilton cycles, and we do this by applying Lemma \ref{manyrepairs} 
with $H_i=\sigma_{v_i}(H)$ and $V_i=\{a_\infty^{v_i},b_\infty^{v_i}\}$ for $i=1,2,\ldots,m$.
To this end, observe that \linebreak $\sigma_{v_1}(H),\sigma_{v_2}(H),\ldots,\sigma_{v_m}(H)$ are edge-disjoint subgraphs of $L$,
and since $b_\infty^{v_i}=a_\infty^{v_{i+1}}$ for $i=1,2,\ldots,m-1$ it follows that $V_i\cap V_{i+1}\ne \emptyset$ for $i=1,2,\ldots,m-1$ as required. It remains to show there is a $2$-factorisation
$\{J_1,J_2\}$ of $L$ such that 
\begin{itemize}
\item $\{a_\infty^{v_1},a_\infty^{v_2},\ldots,a_\infty^{v_m}\}$ links $\{J_1,J_2\}$; and
\item $\sigma_{v_i}(H)$ induces an $\{a_\infty^{v_{i}},b_\infty^{v_i}\}$-connector in $\{J_1,J_2\}$, for $i=1,2,\ldots,m$.
\end{itemize} 
Let $\{U,U'\}$ be a partition of $V(G)$ such that both $U$ and $U'$ link  $\{\F_j : j\in T\}$ (such a partition exists by Lemma \ref{intersects}), 
and let $\{J_1,J_2\}$ be the decomposition of $L$ defined by 
\begin{itemize}
\item $J_1=(\bigcup_{v\in U}\sigma_v(S\cup R\cup C))\cup (\bigcup_{v\in U'}\sigma_v(S\cup R\cup C'))$;
\item $J_2= (\bigcup_{v\in U}\sigma_v(S'\cup R'\cup C'))\cup (\bigcup_{v\in U'}\sigma_v(S'\cup R'\cup C))$.
\end{itemize}
It follows easily from Lemma \ref{linking2factors2} that both $J_1$ and $J_2$ are $2$-factors of $L$, and that $\{a_\infty^{v_1},a_\infty^{v_2},\ldots,a_\infty^{v_m}\}$ links $\{J_1,J_2\}$ as required.

Finally, observe that 
$E(\sigma_{v_i}(S))\subseteq E(J_1)\cap E(\sigma_{v_i}(H))$ and $E(\sigma_{v_i}(S'))\subseteq E(J_2)\cap E(\sigma_{v_i}(H))$ for  $i=1,2,\ldots,m$. Then, since $\{E(S),E(S')\}$ is an $\{a_{\infty},b_{\infty}\}$-connector, 
it follows that $\sigma_{v_i}(H)$ induces an $\{a_\infty^{v_{i}},b_\infty^{v_i}\}$-connector in $\{J_1,J_2\}$ for $i=1,2,\ldots,m$, as required.
\qed

\subsection{Regular graphs of odd degree}

\begin{lemma}\label{linking2factors2_odd}
Suppose $G$ is a $(2n+1)$-regular graph, $F$ is a $1$-factor of $G$, \linebreak $\{\F_0,\F_1,\ldots,\F_{n-2},\F_\infty\}$ is a directed $2$-factorisation of $G-F$, $t\in N_n$, $\{U,U'\}$ is a partition of $V(G)$ such that both $U$ and $U'$ link $\{F,\F_t\}$, and $X$, $X'$, $C$ and $C'$ are subgraphs of $K_{A_n\cup B_n\cup \{c\}}$ such that 
\begin{itemize}
\item [$(1)$] $\deg_{X}(u)=\deg_{X'}(u)$ for each $u\in A_n\cup B_n$;
\item [$(2)$] $\{\deg_{X}(c),\deg_{X'}(c)\}=\{0,2\}$ and $\{|E(X\cup C)|, |E(X'\cup C')|\}=\{n,n+1\}$;
\item [$(3)$] if $u$ and $v$ are vertices such that $\deg_X(u)=\deg_X(v)=1$, then $u$ and $v$ belong to the same component of $X'$ if and only if they belong to the same component of $X$;
\item [$(4)$] $\deg_{C}(u)=\deg_{C'}(u)$ for each $u\in A_n\cup B_n\cup\{c\}$; 
\item [$(5)$] $C\cup C'$ is a $4$-cycle or a $6$-cycle with $a_t b_t\in E(C')$; 
\item [$(6)$] $^{A}C$ is a path;
\item [$(7)$] $^{A}C'$ consists of a path and a loop on vertex $t$;
\item [$(8)$] $^{A}(X\cup C)$ is a cycle; and
\item [$(9)$] $^{A}(X'\cup C')$ is the vertex disjoint union of a cycle and a loop on vertex $t$.
\end{itemize}
Then the graph $$J=\left(\bigcup_{v\in U}\sigma_v(X\cup C)\right)\cup \left(\bigcup_{v\in U'}\sigma_v(X'\cup C')\right)$$
is a $2$-factor of $L(G)$ and $\S_i$ links $J$ for each $i\in N_n\setminus \{t\}$ such that $\deg_X(a_i)$ or $\deg_X(b_i)=1$.
\end{lemma}

\proof
Let $uv\in E(G)$. First suppose $uv\notin F$, let $F_i$ be the $2$-factor containing $uv$, and let $uv$ be oriented from $u$ to $v$ in $\F_i$. 
Then $\deg_J(uv)=\deg_{Y_u}(b_i)+\deg_{Y_v}(a_i)$ where $Y_u=X\cup C$ if $u\in U$, $Y_u=X'\cup C'$ if $u\in U'$, $Y_v=X\cup C$ if $v\in U$ and $Y_v=X'\cup C'$ if $v\notin U$.
By $(1)$ and $(4)$ we have $\deg_{X\cup C}(b_i)=\deg_{X'\cup C'}(b_i)$ and $\deg_{X\cup C}(a_i)=\deg_{X'\cup C'}(a_i)$. 
So $\deg_J(uv)=\deg_{X\cup C}(b_i)+\deg_{X\cup C}(a_i)=\deg_{^A(X\cup C)}(i)$.
However, it follows from $(1)$, $(2)$, $(4)$ and $(8)$ that either $^A(X\cup C)$ is an $n$-cycle and $\deg_X(c)=0$ or 
$^A(X\cup C)$ is an $(n+1)$-cycle and $\deg_X(c)=2$. In either case, we have $\deg_{^A(X\cup C)}(i)=2$, and hence $\deg_J(uv)=2$.

Now suppose $uv\in F$. Then $\deg_J(uv)=\deg_{Y_u}(c)+\deg_{Y_v}(c)$ where $Y_u=X\cup C$ if $u\in U$, $Y_u=X'\cup C'$ if $u\in U'$, $Y_v=X\cup C$ if $v\in U$ and $Y_v=X'\cup C'$ if $v\notin U$. But since $U$ and $U'$ each link $F$, exactly one of $u$ and $v$ is in $U$ and the other is in $U'$. So we have $\deg_J(uv)=\deg_{Y_u}(c)+\deg_{Y_v}(c)=\deg_{X\cup C}(c)+\deg_{X'\cup C'}(c)=0+2=2$ (because it follows from $(2)$, $(8)$ and $(9)$ that $\deg_C(c)=\deg_{C'}(c)=0$).
This completes the proof that $J$ is a $2$-factor of $L(G)$.

We now show that $\S_i$ links $J$ for each $i\in N_n\setminus\{t\}$ such that $\deg_X(a_i)$ or $\deg_X(b_i)=1$.
Since $^A(X\cup C)$ is a cycle (by $(8)$), $X$ is a union of vertex-disjoint paths, and since $^A(X'\cup C')$ 
is the vertex-disjoint union of a cycle and a loop on vertex $t$ (by $(9)$), $X'$ is a union of vertex-disjoint paths 
(the loop in $^A(X'\cup C')$ arises from the edge $a_tb_t$ of $C'$). Also, $^AX$ is a path and $^AX'$ is a path. 
Moreover, by $(1)$ and $(3)$, 
for each maximal path $P$ of $X$ there is corresponding maximal path $P'$ in $X'$ such that $P$ and $P'$ have the same endpoints. Let $P_1,P_2,\ldots,P_k$ the vertex disjoint paths that comprise $X$ and let $P'_1,P'_2,\ldots,P'_k$ the vertex disjoint paths that comprise $X'$
where $P_i$ has the same endpoints as $P'_i$ for $i=1,2,\ldots,k$. 
These endpoints are precisely the vertices where $\deg_X(a_i)=1$ or $\deg_X(b_i)=1$ (and the values of $i\in N_n\setminus \{t\}$ where $\deg_X(a_i)=1$ or $\deg_X(b_i)=1$ are the values of $i$ for which we need to show $\S_i$ links $J$). 

Let $Z$ be a connected component of $J$. 
It follows from the observations made in the 
preceding paragraph that if $Z$ contains any vertex of $\sigma_v(X)$ for any $v\in U$ or $\sigma_v(X')$ for any $v\in U'$, 
then for $i=1,2,\ldots,k$,
$Z$ contains either $\sigma_v(P_i)$ for some $v\in U$ or $\sigma_v(P'_i)$ for some $v\in U'$. Thus, $Z$ contains a vertex of $\S_i$ for each 
$i\in N_n\setminus \{t\}$ such that $\deg_X(a_i)$ or $\deg_X(b_i)=1$.

Thus, we can assume that $Z$ contains no vertex of $\sigma_v(X)$ for any $v\in U$ and no vertex of $\sigma_v(X')$ for any $v\in U'$. 
It is a consequence of $(4)$ and the fact that $C\cup C'$ is a $4$-cycle or a $6$-cycle that each of $C$ and $C'$ is a matching with $2$ or $3$ edges. However, since $^AC$ is a path (by $(6)$) and $^AC'$ consists of a path and a loop on vertex $t$ (by $(7)$), it can be seen 
that every edge of $Z$ is of the form $a^v_tb^v_t$ where $v\in U'$ (otherwise we must have a vertex of $\sigma_v(X)$ for some $v\in U$ 
or a vertex of $\sigma_v(X')$ for some $v\in U'$, and we have assumed that this is not the case). 
This means that that there is a cycle in $\F_t$ that contains no vertex of $U$, contradicting the fact that $U$ links $\F_t$. 
We conclude that $Z$ contains a vertex of $\S_i$ for each 
$i\in N_n\setminus \{t\}$ such that $\deg_X(a_i)$ or $\deg_X(b_i)=1$, and this completes the proof.
\qed

\begin{lemma}\label{4regular_odd}
Suppose $t\in\{0,1,\ldots,n-2\}$, and $H$ is a subgraph of $K_{A_n\cup B_n\cup\{c\}}$ that admits decompositions $\{S,S',R,R',C,C'\}$ and $\{T,T',R,R',C,C'\}$ such that:
\begin{itemize}
\item [$(1)$] $\deg_{S}(u)=\deg_{T}(u)$ and $\deg_{S'}(u)=\deg_{T'}(u)$ for each $u\in A_n\cup B_n$;

\item [$(2)$] $\deg_S(a_\infty)=\deg_S(b_\infty)=1$ and $\deg_{S'}(a_\infty)=deg_{S'}(b_\infty)=1$;

\item [$(3)$] $\{\deg_{S}(c),\deg_{T}(c)\}=\{\deg_{S'}(c),\deg_{T'}(c)\}=\{0,2\}$ and $\{|E(S\cup R\cup C)|,|E(T\cup R\cup C')|\}=\{|E(S'\cup R'\cup C')|,|E(T'\cup R'\cup C)|\}=\{n,n+1\}$;
\item [$(4)$] if $u$ and $v$ are distinct vertices such that $\deg_S(u)=\deg_S(v)=1$, then $u$ and $v$ belong to the same component of $T$ if and only if they belong to the same component of $S$;
\item [$(5)$] if $u$ and $v$ are distinct vertices such that $\deg_{S'}(u)=\deg_{S'}(v)=1$, then $u$ and $v$ belong to the same component of $T'$ if and only if they belong to the same component of $S'$;
\item [$(6)$] $\deg_{C}(u)=\deg_{C'}(u)$ for each $u\in A_n\cup B_n\cup \{c\}$;
\item [$(7)$] $C\cup C'$ is a $4$-cycle or a $6$-cycle with $a_t b_t\in E(C')$; 
\item [$(8)$] $^{A}C$ is a path;
\item [$(9)$] $^{A}C'$ is the vertex disjoint union of a path and a loop on vertex $t$;
\item [$(10)$] each of $^{A}(S\cup R\cup C)$ and $^{A}(T'\cup R'\cup C)$ is a cycle;
\item [$(11)$] each of $^{A}(T \cup R\cup C')$ and $^{A}(S' \cup R'\cup C')$ is the vertex disjoint union of a cycle and a loop on vertex $t$;
 and
\item [$(12)$] both $\{E(S),E(S')\}$ and $\{E(T),E(T')\}$ are $\{a_{\infty},b_{\infty}\}$-connectors.
\end{itemize} 
Then $H$ is an $(R\cup R')$-adjustable Hamilton fragment.
\end{lemma}

\proof
Suppose $t\in\{0,1,\ldots,n-2\}$ and $\{S,S',R,R',C,C'\}$ and  $\{T,T',R,R',C,C'\}$ are decompositions of $H$ that satisfy the conditions of the lemma. It is easy to see that if $Q$ and $Q'$ are edge-disjoint subgraphs of $K_{A_n\cup B_n\cup\{c\}}-(S\cup S'\cup C\cup C')$ satisfying $^{A}Q=\prescript{A}{}R$ and $^{A}Q'=\prescript{A}{}R'$, then $H-(R\cup R')+(Q\cup Q')$ has decompositions, namely $\{S,S',Q,Q',C,C'\}$ and $\{T,T',Q,Q',C,C'\}$ which also satisfy the conditions of the lemma. Thus we need only show $H$ is a Hamilton fragment and it follows immediately that it is $(R\cup R')$-adjustable. 

Let $G$ be a $2n+1$-regular graph, let $F$ be a $1-$factor of $G$ and let \linebreak $\{\F_0,\F_1,\ldots,\F_{n-2},\F_\infty\}$ be a directed $2$-factorisation of $G-F$ where $\F_\infty$ is a Hamilton cycle, say $\F_\infty=(v_1,v_2,\ldots,v_m)$. Our aim is to show that the subgraph $$L=\bigcup_{v\in V(G)}\sigma_v(H)$$
of $L(G)$ decomposes into Hamilton cycles, and we do this by applying Lemma \ref{manyrepairs} 
with $H_i=\sigma_{v_i}(H)$ and $V_i=\{a_\infty^{v_i},b_\infty^{v_i}\}$ for each $i=1,2,\ldots,m$.
To this end, observe that \linebreak $\sigma_{v_1}(H),\sigma_{v_2}(H),\ldots,\sigma_{v_m}(H)$ are edge-disjoint subgraphs of $L$,
and since $b_\infty^{v_i}=a_\infty^{v_{i+1}}$ for $i=1,2,\ldots,m-1$ it follows that $V_i\cap V_{i+1}\ne \emptyset$ for $i=1,2,\ldots,m-1$ as required. It remains to show there is a $2$-factorisation
$\{J_1,J_2\}$ of $L$ such that 
\begin{itemize}
\item $\{a_\infty^{v_1},a_\infty^{v_2},\ldots,a_\infty^{v_m}\}$ links $\{J_1,J_2\}$; and
\item $\sigma_{v_i}(H)$ induces an $\{a_\infty^{v_{i}},b_\infty^{v_i}\}$-connector in $\{J_1,J_2\}$, for $i=1,2,\ldots,m$.
\end{itemize} 
Let $\{U,U'\}$ be a partition of $V(G)$ such that both $U$ and $U'$ link $\{F,\F_t\}$ (such a partition exists by Lemma \ref{intersects}),   
and let $\{J_1,J_2\}$ be the decomposition of $L$ defined by 
\begin{itemize}
\item $J_1=(\bigcup_{v\in U}\sigma_v(S\cup R\cup C))\cup (\bigcup_{v\in U'}\sigma_v(T\cup R\cup C'))$;
\item $J_2= (\bigcup_{v\in U}\sigma_v(S'\cup R'\cup C'))\cup (\bigcup_{v\in U'}\sigma_v(T'\cup R'\cup C))$.
\end{itemize}
It follows easily from Lemma \ref{linking2factors2_odd} that both $J_1$ and $J_2$ are in fact $2$-factors of $L$ and that  
$\{a_\infty^{v_1},a_\infty^{v_2},\ldots,a_\infty^{v_m}\}$ links $\{J_1,J_2\}$ as required.

Finally, observe that 
$E(\sigma_{v_i}(S))\subseteq E(J_1)\cap E(\sigma_{v_i}(H))$ and $E(\sigma_{v_i}(S'))\subseteq E(J_2)\cap E(\sigma_{v_i}(H))$ whenever $v_i\in U$, and 
$E(\sigma_{v_i}(T))\subseteq E(J_1)\cap E(\sigma_{v_i}(H))$ and $E(\sigma_{v_i}(T'))\subseteq E(J_2)\cap E(\sigma_{v_i}(H))$ whenever $v_i\in U'$. Then, since both $\{E(S),E(S')\}$ and $\{E(T),E(T')\}$ are  $\{a_{\infty},b_{\infty}\}$-connectors, 
it follows that $\sigma_{v_i}(H)$ induces an $\{a_\infty^{v_{i}},b_\infty^{v_i}\}$-connector in $\{J_1,J_2\}$, for $i=1,2,\ldots,m$ as required.
\qed

\begin{lemma}\label{12regular_odd}
Suppose $t\in\{0,1,\ldots,n-2\}$ and $H$ is a subgraph of $K_{A_n\cup B_n\cup\{c\}}$ that admits a decomposition $\{X,X',P_1, P_2,H_1,H_2\}$ such that $X=[a_\infty,c,b_\infty]$, $X'=[a_\infty,b_\infty]$, and
for each $i\in\{1,2\}$: 
\begin{itemize}
\item [$(1)$] $P_i$ is a path with end vertices $a_t$ and $b_t$;
\item [$(2)$] $^{A}P_i$ is an $(n-1)$-cycle on $\{0,1,\ldots,n-2\}$;
\item [$(3)$] there is a $t_i\in\{0,1,\ldots,n-2\}$ and decompositions $\{S_i,S_i',R_i,R_i',C_i,C_i'\}$ and $\{T_i,T_i',R_i,R_i',C_i,C_i'\}$ of $H_i$ that satisfy the conditions of Lemma \ref{4regular_odd};
\item [$(4)$] both $\{E(S_i),E(S_i'),E(X\cup P_i)\}$ and $\{E(T_i),E(T_i'),E(X\cup P_i)\}$ are $\{a_{\infty},b_{\infty},u_i\}$-connectors, for some $u_i\in V(P_i)$.
\end{itemize} 
Then $H$ is a Hamilton fragment.
\end{lemma}

\proof
Let $G$ be a $(2n+1)$-regular graph, let $F$ be a $1-$factor of $G$, let \linebreak $\{\F_0,\F_1,\ldots,\F_{n-2},\F_\infty\}$ be a directed $2$-factorisation of $G-F$ where $\F_\infty$ is a Hamilton cycle, say $\F_\infty=(v_1,v_2,\ldots,v_m)$, and let $\{U,U'\}$ be a partition of $V(G)$ such that both $U$ and $U'$ link $\{F,\F_{t}\}$ (such a partition exists by Lemma \ref{intersects}). Our aim is to show that the subgraph $$L=(\bigcup_{v\in U}\sigma_v(H_1\cup X\cup P_1))\cup (\bigcup_{v\in U'}\sigma_v(H_1\cup X'\cup P_1))$$ of $L(G)$ 
decomposes into Hamilton cycles, and we do this by applying Lemma \ref{manyrepairs} with
\begin{itemize}
\item $H_i=\sigma_{v_i}(H_1\cup X\cup P_1)$ and $V_i=\{a_\infty^{v_i},b_\infty^{v_i},u_1^{v_i}\}$ for each $v_i\in U$; and
\item $H_i=\sigma_{v_i}(H_1\cup X'\cup P_1)$ and $V_i=\{a_\infty^{v_i},b_\infty^{v_i}\}$ for each $v_i\in U'$.
\end{itemize}
It then follows, by symmetry, that the subgraph
$$L'=(\bigcup_{v\in V(G)}\sigma_v(H))\setminus L =(\bigcup_{v\in U}\sigma_v(H_2\cup X'\cup P_2))\cup (\bigcup_{v\in U'}\sigma_v(H_2\cup X\cup P_2))$$ of $L(G)$ decomposes into Hamilton cycles, and the result then follows.

To this end, observe that $H_1,H_2,\ldots,H_m$ are edge-disjoint subgraphs of $L$,
and since $b_\infty^{v_i}=a_\infty^{v_{i+1}}$ for $i=1,2,\ldots,m-1$ it follows that $V_i\cap V_{i+1}\ne \emptyset$ for $i=1,2,\ldots,m-1$ as required. It remains to show there is a $2$-factorisation
$\{J_1,J_2,J_3\}$ of $L$ such that
\begin{itemize}
\item $\{a_\infty^{v_1},a_\infty^{v_2},\ldots,a_\infty^{v_m}\}\cup \{u_1^{v}\mid v\in U\}$ links $\{J_1,J_2,J_3\}$; and 
\item $H_i$ induces a $V_i$-connector in $\{J_1,J_2,J_3\}$ for each $i=1,2,\ldots,m$; that is
\begin{itemize}
\item[$\bullet$] $\sigma_{v}(H_1\cup X\cup P_1)$ induces an $\{a_\infty^{v},b_\infty^{v},u_1^{v}\}$-connector in $\{J_1,J_2,J_3\}$ whenever $v\in U$; and 
\item[$\bullet$] $\sigma_{v}(H_1\cup X'\cup P_1)$ induces an $\{a_\infty^{v},b_\infty^{v}\}$-connector in $\{J_1,J_2,J_3\}$ whenever $v\in U'$.
\end{itemize}
\end{itemize} 

Let $\{W,W'\}$ be a partition of $V(G)$ such that both $W$ and $W'$ link $\{F,\F_{t_1}\}$ (such partitions exist by Lemma \ref{intersects}), and let $\{J_1,J_2,J_3\}$ be the decomposition of $L$ defined by 
\begin{itemize}
\item $J_{1}=(\bigcup_{v\in W}\sigma_v(S_1\cup R_1\cup C_1))\cup (\bigcup_{v\in W'}\sigma_v(T_1\cup R_1\cup C'_1))$;
\item $J_{2}= (\bigcup_{v\in W}\sigma_v(S'_1\cup R'_1\cup C'_1))\cup (\bigcup_{v\in W'}\sigma_v(T'_1\cup R'_1\cup C_1))$; and
\item $J_3=(\bigcup_{v\in U}\sigma_v(X\cup P_1))\cup (\bigcup_{v\in U'}\sigma_v(X'\cup P_1))$.
\end{itemize}
As in the proof of Lemma \ref{4regular_odd}, it follows easily from Lemma \ref{linking2factors2_odd} that both of $J_1$ and $J_2$ are $2$-factors of $L$, and that  
$\{a_\infty^{v_1},a_\infty^{v_2},\ldots,a_\infty^{v_m}\}$ links $\{J_1,J_2\}$ as required.
Furthermore, it is easy to see that $J_3$ is a $2$-factor of $L$ and, since $U$ links $\F_{t}$ and 
$u_1\in V(P_1)$, it follows easily from the conditions of the lemma that $\{a_\infty^{v_1},a_\infty^{v_2},\ldots,a_\infty^{v_m}\}\cup \{u_1^{v}\mid v\in U\}$ links $J_3$ as required.

Finally, observe that 
\begin{itemize}
\item $E(\sigma_{v_i}(S_1))\subseteq (E(J_1)\cap H_i)$ and $E(\sigma_{v_i}(S'_1))\subseteq (E(J_2)\cap H_i)$ whenever $v_i\in W$;
\item $E(\sigma_{v_i}(T_1))\subseteq (E(J_1)\cap H_i)$ and $E(\sigma_{v_i}(T'_1))\subseteq (E(J_2)\cap H_i)$ whenever $v_i\in W'$; 
\item $E(\sigma_{v_i}(X\cup P_1))\subseteq (E(J_3)\cap H_i)$ whenever $v_i\in U$; and 
\item $E(\sigma_{v_i}(X'))\subseteq (E(J_3)\cap H_i)$ whenever $v_i\in U'$.
\end{itemize}
Then, since both $\{E(S_1),E(S_1'),E(X\cup P_1)\}$ and $\{E(T_1),E(T_1'),E(X\cup P_1)\}$ are \linebreak $\{a_{\infty},b_{\infty},u_1\}$-connectors (by assumption), and both $\{E(S_1),E(S_1'),E(X')\}$ and \linebreak $\{E(T_1),E(T_1'),E(X')\}$ are  $\{a_{\infty},b_{\infty}\}$-connectors (by the properties of Lemma \ref{4regular_odd} and the fact that $E(X')=\{a_\infty b_\infty\}$), 
it follows that $H_i$ induces a $V_i$-connector in $\{J_1,J_2,J_3\}$, for each $i=1,2,\ldots,m$ as required.\qed

\section{Decompositions into Hamilton fragments}\label{Section4}

\begin{lemma}\label{SuffCond2}
Let $n$ be a positive integer, let $K\in \{K_{A_{n}\cup B_{n}},K_{A_{n}\cup B_{n}\cup\{c\}}\}$, let 
$\rho$ be a permutation, of order $q$ say, on $\Z_{n-1}$, and let $\rho$ act on $V(K)$ by mapping $a_i$ to $a_{\rho(i)}$ and $b_i$ to $b_{\rho(i)}$ for each $i\in\Z_{n-1}$ (keeping all other vertices fixed).
If there are pairwise edge-disjoint subgraphs $I$, $Z$ and $Q$ of $K$ such that:
\begin{itemize}
\item[$(1)$] $Z\cup Q$ is a $Q$-adjustable Hamilton fragment;
\item[$(2)$] the orbit of $^{A}(Z\cup Q)$ under $\rho$ decomposes $^{A}(K-I)$; and
\item[$(3)$] the graphs $I,\rho^0(Z),\rho^1(Z),\ldots,\rho^{q-1}(Z)$ are pairwise edge-disjoint;
\end{itemize}
then there is a decomposition $\D=\{\rho^i(Z)\cup Q_i: i\in \Z_q\}$ of $K-I$ into Hamilton fragments in which $^{A}Q_i=\prescript{A}{}(\rho^i(Q))=\rho^i(^{A}Q)$, for $i=0,1,\ldots,q-1$.
Furthermore, if
\begin{itemize}
\item[$(P1)$] $Q$ is edge-disjoint from each of $I, \rho^0(Z), \rho^1(Z),\ldots, \rho^{q-1}(Z)$;
\end{itemize}
 then there is such a decomposition $\D$ in which $Q_0=Q$; and if
\begin{itemize}
\item [$(P2)$] $I\cup \left( \bigcup_{H\in \D'} H\right)$ is a Hamilton fragment for some subset $\D'$ of $\D$,
\end{itemize}
 then $K$ admits a decomposition into Hamilton fragments.
\end{lemma}

\proof 
It follows directly from properties $(2)$ and $(3)$ that there is a decomposition $\D=\{\rho^i(Z)\cup Q_i: i\in\Z_q\}$
of $K-I$, where $Q_0,Q_1,\ldots,Q_{q-1}$ are graphs satisfying $^{A}Q_i=\prescript{A}{}(\rho^i(Q))=\rho^i(^{A}Q)$, for $i=0,1,\ldots,q-1$.
Then, since $\rho$ fixes $\infty$, property $(1)$ ensures that $\rho^i(Z)\cup \rho^i(Q)$ is a $\rho^i(Q)$-adjustable Hamilton fragment for each $i\in\Z_q$.  
Thus $\rho^i(Z)\cup Q_i$ is a Hamilton fragment for each $i\in\Z_q$ and each element of $\D$ is a Hamilton fragment as required.

Furthermore, if $(P1)$ holds then $Q$ is a subgraph of $Q_0\cup Q_1\cup \cdots \cup Q_{q-1}$ and, since $^A Q=\prescript{A}{}Q_0$, we are free to set $Q_0=Q$ in our above definition of $\D$. 

Finally, if $(P2)$ holds, say $H'=I\cup \left( \bigcup_{H\in \D'} H\right)$ is a Hamilton fragment for some subset $\D'$ of $\D$, then $\{H'\}\cup\{H:H\in \D\setminus\D'\}$ is the required decomposition of $K$ into Hamilton fragments.\qed

\begin{lemma}\label{alln}
There is a decomposition of $K_{A_{n}\cup B_{n}}$ into Hamilton fragments for each positive integer $n$.
\end{lemma}

\proof When $n=1$ the result is obvious. Furthermore, for each $n\in\{2,3,4,5,7\}$ a suitable decomposition is given in Section \ref{Section5}. Suppose then that $n\notin\{1,2,3,4,5,7\}$.

Let $K=K_{A_{n}\cup B_{n}}$ and let $\rho$ be the permutation on $\Z_{n-1}$ which maps $v$ to $v+2$ for each $v\in\Z_{n-1}$. Observe that $\rho$ has order $n-1$ when $n$ is even, and order $(n-1)/2$ when $n$ is odd.
Our aim is to show, for each value $n$, that there are edge-disjoint subgraphs $I$, $Z$ and $Q$ of $K$ which, together with $\rho$, satisfy properties $(2)-(3)$ and $(P1), (P2)$ (with $\D'=\{Z\cup Q\}$) of Lemma \ref{SuffCond2}. The problem now splits according to the parity of $n$.\\

\noindent{\bf Case 1: $n\geq 6$ is even}

Let $n=2m$ and let $I$ be the $1$-factor of $K$ with $$E(I)=\{a_\infty b_\infty,a_0b_1,a_1b_2,\ldots,a_{2m-3}b_{2m-2},a_{2m-2}b_0\}.$$
Observe that both $I$ and $^{A}I$ are fixed under the permutation $\rho$. 
Define 
\begin{itemize}
\item $O_0=\{a_ib_i\mid i\in \Z_{2m-1}\}$; 
\item $O_1=\{a_ia_{i+1},b_ib_{i+1},a_{i+1}b_{i}\mid i\in \Z_{2m-1}\}$;
\item $O_j=\{a_ia_{i+j},b_ib_{i+j},a_{i+j}b_{i},a_ib_{i+j}\mid i\in \Z_{2m-1}\}$ for each $j=2,3,\ldots,m-1$; and
\item $O_\infty = \{a_ia_{\infty},b_ib_{\infty},a_{\infty}b_{i},a_ib_{\infty}\mid i\in \Z_{2m-1}\}$.
\end{itemize}
Observe that $\{O_0,O_1,\ldots,O_{m-1},O_\infty\}$ partitions $E(K-I)$, and that each $O_i$ is the union of one or more edge orbits of $K-I$ under $\rho$. Furthermore, for any subgraph $H$ of $K-I$ we define  
$$\mathcal{O}(H)=(|E(H)\cap O_0|,|E(H)\cap O_1|,\ldots,|E(H)\cap O_{m-1}|,|E(H)\cap O_\infty|).$$
Observe that 
$$\mathcal{O}(K-I)=(2m-1,3(2m-1),4(2m-1),4(2m-1),\ldots,4(2m-1)),$$
and that the orbit of $^{A}H$ under $\rho$ decomposes $^{A}(K-I)$, whenever $H$ is a subgraph of $K-I$ such that $\mathcal{O}(H)=(1,3,4,4,\ldots,4).$\\

Suppose firstly that $n=6$.  
Let $Z=C\cup C' \cup S\cup S'$ and $Q=R\cup R'$ where 
\begin{itemize}
\item $C=[b_{4},b_{2}]\cup [a_{3},a_{4}]$;
\item $C'=[a_{4},b_{4}]\cup [b_{2},a_{3}]$;
\item $S=[a_\infty,b_{1}]\cup [b_\infty,a_{0},b_2]$; 
\item $S'=[a_\infty,a_{0}]\cup [b_\infty,b_{1},b_{2}]$;
\item $R=[a_{1},a_{3}]$; and
\item $R'=[a_{0},a_{3}]$.
\end{itemize}

\noindent Property $(1)$ follows from Lemma \ref{4regular}, setting $T=\{4\}$ and $H=Z\cup Q$, and noting in particular that $\{E(S),E(S')\}$ is an $\{a_\infty,b_\infty\}$-connector by Lemma \ref{connectors} $(1)$, with $(\alpha,\beta,u,v,w)=(b_\infty,a_\infty,a_0,b_2,b_1)$.\\

\noindent 
Property $(2)$ follows by noting that 
\begin{itemize}
\item $\mathcal{O}(C)=(0,1,1,0)$;
\item $\mathcal{O}(C')=(1,1,0,0)$;
\item $\mathcal{O}(S)=(0,0,1,2)$;
\item $\mathcal{O}(S')=(0,1,0,2)$;
\item $\mathcal{O}(R)=(0,0,1,0)$;
\item $\mathcal{O}(R')=(0,0,1,0)$;
\end{itemize}
and hence $\mathcal{O}(Z\cup Q)=(1,3,4,4)$ as required.\\

\noindent Property $(3)$ follows by noting that each of the ten edges in $E(Z)$ belong to distinct edge orbits of $K-I$ under $\rho$.

\noindent Property $(P1)$ follows by noting that $\mathcal{O}(Z)=(1,3,2,4)$, $\mathcal{O}(Q)=(0,0,2,0)$ and the edge orbits of the two edges in $E(Z)\cap O_2$ are distinct from the edge orbits of the two edges in $E(Q)\cap O_2$.\\

\noindent Property $(P2)$ follows by seting $D'=\{Z\cup Q\}$ and applying Lemma \ref{4regular_withHam} with  
\begin{itemize}
\item $(r,s,t)=(3,1,4)$;
\item $H=I \cup Z \cup Q$;
\item $X_1 = C\cup S\cup R$;
\item $X_2= C'\cup R' + \{a_\infty b_\infty,a_{0}b_{1},a_1b_2\}$;  
\item $X_3 = I \cup S' -\{a_\infty b_\infty,a_{0}b_{1},a_1b_2\}$; 
\end{itemize}
 and noting in particular that $\{E(X_1),E(X_3),E(X_2)\}$ is an $\{a_\infty,b_\infty\}$-connector by Lemma \ref{connectors} $(2)$, with $(\alpha,\beta,u,v,w)=(b_\infty,a_\infty,a_0,b_2,b_1)$.\\

Suppose then $n\geq 8$. 
Let $Z=C\cup C' \cup S\cup S'$ and $Q=R\cup R'$ where
\begin{itemize}
\item $C=[a_m,a_{m+2}]\cup [a_{m+1},b_m]$; 
\item $C'=[a_m,b_m]\cup [a_{m+1},a_{m+2}]$;
\item $S=[a_\infty,b_{1},b_{2}] \cup [b_\infty,a_{0},b_{3}]$; 
\item $S'=[a_\infty,a_{0},b_{2}]\cup [b_\infty,b_{1},b_{3}]$;
\item $R=[b_{2},b_{2m-2},b_{4},b_{2m-4},\ldots,b_{m-2},b_{m+2}]$\\
$\cup$ $[b_{3},b_{2m-3},b_{5},b_{2m-5},\ldots,b_{m-1},a_{m+1}]$; and
\item $R'=[b_{2},b_{2m-3},b_{4},b_{2m-5},\ldots,b_{m-2},a_{m+1}]$ \\
$\cup$ $[b_{3},b_{2m-2},b_{5},b_{2m-4},\ldots,b_{m-1},b_{m+2}]$;
\end{itemize} 
when $m$ is even, and 
\begin{itemize}
\item $C=[b_{m+1},b_{m-1}]\cup [a_{m},a_{m+1}]$;
\item $C'=[a_{m+1},b_{m+1}]\cup [b_{m-1},a_{m}]$;
\item $S=[a_\infty,b_{1},b_{2}]\cup [b_\infty,a_{0},a_{3}]$; 
\item $S'=[a_\infty,a_{0},b_{2}]\cup [b_\infty,b_{1},a_{3}]$;
\item $R=[b_{2},b_{2m-2},b_{4},b_{2m-4},\ldots,b_{m-3},b_{m+3},b_{m-1}]$ \\
$\cup$ $[a_{3},b_{2m-3},b_{5},b_{2m-5},\ldots,b_{m-2},b_{m+2}]\cup [a_{m+2},a_m]$; and
\item $R'=[b_{2},b_{2m-3},b_{4},b_{2m-5},\ldots,b_{m-3},b_{m+2}]\cup [a_{m+2},b_{m-1}]$ \\
$\cup$ $[a_{3},b_{2m-2},b_{5},b_{2m-4},\ldots,b_{m-2},b_{m+3},a_m]$;
\end{itemize}
when $m$ is odd.

\noindent Property $(1)$ follows from Lemma \ref{4regular}, setting $T=\{m\}$ and $H=Z\cup Q$ when $m$ is even and $T=\{m+1\}$ and $H=Z\cup Q$ when $m$ is odd, and noting in particular that $\{E(S),E(S')\}$ is an $\{a_\infty,b_\infty\}$-connector by Lemma \ref{connectors} $(1)$, with $(\alpha,\beta,u,v,w)=(a_\infty,b_\infty,b_1,b_2,a_0)$.\\

\noindent 
Property $(2)$ follows by noting that 
\begin{itemize}
\item $\mathcal{O}(C)=(0,1,1,0,0,\ldots,0)$;
\item $\mathcal{O}(C')=(1,1,0,0,\ldots,0)$;
\item $\mathcal{O}(S)=(0,1,0,1,0,0,\ldots,0,2)$;
\item $\mathcal{O}(S')=(0,0,2,0,0\ldots,0,2)$;
\item $\mathcal{O}(R)=(0,0,1,1,2,2,\ldots,2,0)$;
\item $\mathcal{O}(R')=(0,0,0,2,2,\ldots,2,0)$;
\end{itemize}
and hence $\mathcal{O}(Z\cup Q)=(1,3,4,4,\ldots,4)$ as required.\\

\noindent Property $(3)$ follows by noting that each of the twelve edges in $E(Z)$ belong to distinct edge orbits of $K-I$ under $\rho$.\\

\noindent Property $(P1)$ follows by noting that $\mathcal{O}(Z)=(1,3,3,1,0,0,\ldots,0,4)$,  \linebreak $\mathcal{O}(Q)=(0,0,1,3,4,4,\ldots,4,0)$ and the edge orbits of the four edges in $E(Z)\cap (O_2 \cup O_3)$ are distinct from the edge orbits of the four edges in $E(Q)\cap (O_2 \cup O_3)$.\\ 

\noindent Property $(P2)$ follows by setting $\D'=\{Z\cup Q\}$ and applying Lemma \ref{4regular_withHam} with  
\begin{itemize}
\item $(r,s,t)=(3,m,m+2)$ when $m$ is even, and $(r,s,t)=(3,m+1,m+2)$ when $m$ is odd;
\item $H=I \cup Z\cup Q$;
\item $X_1 = C\cup S\cup R$;
\item $X_2= C'\cup S'\cup R'-\{a_\infty a_{0},b_\infty b_{1}\} + \{a_\infty b_\infty,a_{0}b_{1}\}$; 
\item $X_3 = I-\{a_\infty b_\infty,a_{0}b_{1}\}+\{a_\infty a_{0},b_\infty b_{1}\}$;
\end{itemize}
and noting in particular that $\{E(X_1),E(X_3),E(X_2)\}$ is an $\{a_\infty,b_\infty\}$-connector by Lemma \ref{connectors} $(3)$, with $(\alpha,\beta,u,v,w,x)=(b_\infty,a_\infty,a_0,b_3,b_1,b_2)$ when $m$ is even and $(\alpha,\beta,u,v,w,x)=(b_\infty,a_\infty,a_0,a_3,b_1,b_2)$ when $m$ is odd.\\

\noindent{\bf Case 2: $n\geq 9$ is odd}

Let $n=2m+1$ and let $I$ be the $1$-factor of $K$ with 
$$E(I)=\{a_\infty b_\infty,a_0b_1,a_2b_3,\ldots,a_{2m-2}b_{2m-1},a_1b_4,a_3b_6,\ldots,a_{2m-1}b_2\}.$$
Observe that both $I$ and $^{A}I$ are fixed under the permutation $\rho$. 
Define 
\begin{itemize}
\item $O_0=\{a_ib_i\mid i\text{ even, } i\in \Z_{2m}\}$; 
\item $O_1=\{a_ia_{i+1},b_ib_{i+1},a_{i+1}b_{i}\mid i\text{ even, } i\in \Z_{2m}\}$;
\item $O_j=\{a_ia_{i+j},b_ib_{i+j},a_{i+j}b_{i},a_ib_{i+j}\mid i\text{ even, } i\in \Z_{2m}\}$ for each $j=2,3,\ldots,m$;
\item $O_\infty = \{a_ia_{\infty},b_ib_{\infty},a_{\infty}b_{i},a_ib_{\infty}\mid i\text{ even, } i\in \Z_{2m}\}$;
\item $O'_0=\{a_ib_i\mid i\text{ odd, } i\in \Z_{2m}\}$; 
\item $O'_3=\{a_ia_{i+3},b_ib_{i+3},a_{i+3}b_{i}\mid i\text{ odd, } i\in \Z_{2m}\}$;
\item $O'_j=\{a_ia_{i+j},b_ib_{i+j},a_{i+j}b_{i},a_ib_{i+j}\mid i\text{ odd, } i\in \Z_{2m}\}$ for each $j=1,2,4,5,\ldots,m$; and
\item $O'_\infty = \{a_ia_{\infty},b_ib_{\infty},a_{\infty}b_{i},a_ib_{\infty}i\text{ odd, } i\in \Z_{2m}\}$.
\end{itemize}
Observe that $\{O_0,O_1,\ldots,O_{m},O_\infty,O'_0,O'_1,\ldots,O'_{m},O'_\infty\}$ partitions $E(K-I)$ (with $O_m=O'_m$ when $m$ is odd), and that each $O_i$ and $O'_i$ is the union of one or more edge orbits of $K-I$ under $\rho$. Furthermore, for any subgraph $H$ of $K-I$ we define  
$$\mathcal{O}(H)=(|E(H)\cap O_0|,|E(H)\cap O_1|,\ldots,|E(H)\cap O_{m}|,|E(H)\cap O_\infty|)$$ 
and 
$$\mathcal{O}'(H)=(|E(H)\cap O'_0|,|E(H)\cap O'_1|,\ldots,|E(H)\cap O'_{m}|,|E(H)\cap O'_\infty|).$$ 
Observe that
$$\mathcal{O}(K-I)=(m,3m,4m,4m,\ldots,4m,2m,4m)$$
and 
$$\mathcal{O}'(K-I)=(m,4m,4m,3m,4m,4m,\ldots,4m,2m,4m)$$
when $m$ is even, and that 
$$\mathcal{O}(K-I)=(m,3m,4m,4m,\ldots,4m,4m,4m)$$
and 
$$\mathcal{O}'(K-I)=(m,4m,4m,3m,4m,4m,\ldots,4m,4m,4m)$$
when $m$ is odd.
It follows that the orbit of $^{A}H$ under $\rho$ decomposes $^{A}(K-I)$, whenever $H$ is a subgraph of $K-I$ such that 
\begin{itemize}
\item $\mathcal{O}(H)=(1,3,4,4,\ldots,4,2,4)$ and $\mathcal{O}'(H)=(1,4,4,3,4,4,\ldots,4,2,4)$ when $m$ is even; and
\item $\mathcal{O}(H)=(1,3,4,4,\ldots,4,4,4)$ and $\mathcal{O}'(H)=(1,4,4,3,4,4,\ldots,4,4,4)$ when $m$ is odd.
\end{itemize}

Let $Z=C\cup C' \cup S\cup S'\cup S_1\cup S_1'$ and $Q=R\cup R'\cup R_1\cup R_1'$ where
\begin{itemize}
\item $C=[a_m,a_{m+1}]\cup [b_{m+1},b_{m+2}]\cup [a_{m+3},b_{m}]$;
\item $C'=[b_m,a_m]\cup [a_{m+1},b_{m+1}]\cup [b_{m+2},a_{m+3}]$;
\item $S=[a_\infty,b_{0},b_{3}]\cup [b_\infty,a_{1},b_{2}]$; 
\item $S'=[a_\infty,a_{1},b_{3}]\cup [b_\infty,b_{0},b_{2}]$;
\item $S_1=[a_\infty,a_{0},a_{3}]\cup [b_\infty,b_{1},a_{2}]$; 
\item $S_1'=[a_\infty,b_{1},a_{3}]\cup [b_\infty,a_{0},a_{2}]$;
\item $R=[b_{2},b_{2m-1},b_{4},b_{2m-3},\ldots,b_{m+5},b_{m-2},b_{m+3}]\\ 
\cup [b_{3},b_{2m-2},b_{5},b_{2m-4},\ldots,b_{m+4},b_{m-1},b_{m+2}]$; 
\item $R'=[b_{2},b_{2m-2},b_{4},b_{2m-4},\ldots,b_{m+4},b_{m-2},b_{m+2}]\\
\cup [b_{3},b_{2m-1},b_{5},b_{2m-3},\ldots,b_{m+5},b_{m-1},b_{m+3}]$;
\item $R_1=[a_{2},a_{2m-1},a_{4},a_{2m-3},\ldots,a_{m+5},a_{m-2},b_{m+3},a_m,b_{m+2}]\\ 
\cup [a_{3},a_{2m-2},a_{5},a_{2m-4},\ldots,a_{m+4},a_{m-1},a_{m+1},a_{m+2}]$; and
\item $R_1'=[a_{2},a_{2m-2},a_{4},a_{2m-4},\ldots,a_{m+4},a_{m-2},a_{m+2},b_m,b_{m+1}]\\
\cup [a_{3},a_{2m-1},a_{5},a_{2m-3},\ldots,a_{m+5},a_{m-1},b_{m+3},b_{m+1}]$.
\end{itemize}

We begin by proving the following useful properties:
\begin{itemize}
\item[(I)] $I\cup S_1\cup S_1'\cup R_1\cup R_1'$ is a Hamilton fragment;
\item[(II)] $C\cup C' \cup S\cup S'\cup R\cup R'$ is an $(R\cup R')$-adjustable Hamilton fragment; and
\item[(III)] $S_1\cup S_1'\cup R_1\cup R_1'$ is an $(R_1\cup R_1')$-adjustable Hamilton fragment.
\end{itemize}

Property (I) follows from Lemma \ref{6regular_2}, setting
\begin{itemize}
\item $(s,t)=(m+2,m+1)$;
\item $H=I\cup S_1\cup S_1'\cup R_1\cup R_1'$; 
\item $X_1 = S_1\cup R_1$;
\item $X_2= S_1'\cup R_1'-\{a_\infty b_{1},b_\infty a_{0}\} + \{a_\infty b_\infty,a_{0}b_{1}\}$; 
\item $X_1'= X_2-\{a_\infty b_\infty, a_{0}b_{1}\} + \{a_\infty a_{0},b_\infty b_{1}\}$;
\item $X_2' = X_1-\{a_\infty a_{0},b_\infty b_{1}\} + \{a_\infty b_\infty,a_{0}b_{1}\}$;
\item $X_3 = I-\{a_\infty b_\infty,a_{0}b_{1}\}+\{a_\infty b_{1},b_\infty a_{0}\}$; 
\end{itemize}
and noting in particular that $\{E(X_1),E(X_3),E(X_2)\}$ is an $\{a_\infty,b_\infty\}$-connector by Lemma \ref{connectors} $(3)$ with $(\alpha,\beta,u,v,w,x)=(a_\infty,b_\infty,a_0,a_3,b_1,a_2)$, 
and that \linebreak $\{E(X_1'),E(X_3),E(X_2')\}$ is an $\{a_\infty,b_\infty\}$-connector by Lemma \ref{connectors} $(3)$ with \linebreak $(\alpha,\beta,u,v,w,x)=(a_\infty,b_\infty,a_0,a_2,b_1,a_3)$.\\

Property (II) follows from Lemma \ref{4regular}, setting
\begin{itemize}
\item $T=\{1,2\}$;
\item $H=C\cup C' \cup S\cup S'\cup R\cup R'$; 
\end{itemize}
and noting in particular that $\{E(S),E(S')\}$ is an $\{a_\infty,b_\infty\}$-connector by Lemma \ref{connectors} $(1)$ with $(\alpha,\beta,u,v,w)=(a_\infty,b_\infty,b_0,b_3,a_1)$.\\

Property (III) follows from Lemma \ref{4regular_easy}, setting
\begin{itemize}
\item $H=S_1\cup S_1'\cup R_1\cup R_1'$;
\item $\{S,S',R,R'\}=\{S_1,S_1',R_1,R_1'\}$; 
\end{itemize}
and noting in particular that $\{E(S_1),E(S'_1)\}$ is an $\{a_\infty,b_\infty\}$-connector by Lemma \ref{connectors} $(1)$ with $(\alpha,\beta,u,v,w)=(a_\infty,b_\infty,a_0,a_3,b_1)$.\\

\noindent Property $(1)$ then follows from (II), (III) and the fact that $\{C\cup C' \cup S\cup S'\cup R\cup R',S_1\cup S_1'\cup R_1\cup R_1',\}$ is a decomposition of $Z\cup Q$ and $\{R,R',R_1,R_2\}$ is a decomposition of $Q$. \\

\noindent Property $(2)$ then follows by noting that 
\begin{itemize}
\item $\mathcal{O}(C)=(0,1,0,1,0,0,\ldots,0)$ and $\mathcal{O}'(C)=(0,1,0,0,\ldots,0)$;
\item $\mathcal{O}(C')=(1,1,0,0,\ldots,0)$ and $\mathcal{O}'(C')=(1,0,0,\ldots,0)$;
\item $\mathcal{O}(S)=(0,0,0,1,0,0,\ldots,0,1)$ and $\mathcal{O}'(S)=(0,1,0,0\ldots,0,1)$;
\item $\mathcal{O}(S')=(0,0,1,0,0,\ldots,0,1)$ and $\mathcal{O}'(S')=(0,0,1,0,0\ldots,0,1)$;
\item $\mathcal{O}(S_1)=(0,0,0,1,0,0,\ldots,0,1)$ and $\mathcal{O}'(S_1)=(0,1,0,0\ldots,0,1)$;
\item $\mathcal{O}(S'_1)=(0,0,1,0,0,\ldots,0,1)$ and $\mathcal{O}'(S'_1)=(0,0,1,0,0\ldots,0,1)$;
\item $\mathcal{O}(R)=(0,0,0,0,0,2,0,2,0,\ldots,2,0,0)$ and $\mathcal{O}'(R)=(0,0,0,2,0,2,0,\ldots,2,0,0)$;
\item $\mathcal{O}(R')=(0,0,0,0,2,0,2,0,\ldots,2,0,1,0)$ and $\mathcal{O}'(R')=(0,0,0,0,2,0,2,0,\ldots,2,0,1,0)$ when $m$ is even;
\item $\mathcal{O}(R')=(0,0,0,0,2,0,2,0,\ldots,2,0,2,0)$ and $\mathcal{O}'(R')=(0,0,0,0,2,0,2,0,\ldots,2,0,2,0)$ when $m$ is odd;
\item $\mathcal{O}(R_1)=(0,0,1,1,0,2,0,2,0,\ldots,2,0,0)$ and $\mathcal{O}'(R_1)=(0,1,1,1,0,2,0,2,0,\ldots,2,0,0)$;
\item $\mathcal{O}(R'_1)=(0,1,1,0,2,0,2,0,\ldots,2,0,1,0)$ and $\mathcal{O}'(R'_1)=(0,0,1,0,2,0,2,0,\ldots,2,0,1,0)$ when $m$ is even;
\item $\mathcal{O}(R'_1)=(0,1,1,0,2,0,2,0,\ldots,2,0,2,0)$ and $\mathcal{O}'(R'_1)=(0,0,1,0,2,0,2,0,\ldots,2,0,2,0)$ when $m$ is odd; 
\end{itemize}
and hence 
\begin{itemize}
\item $\mathcal{O}(Z\cup Q)=(1,3,4,4,\ldots,4,2,4)$ and $\mathcal{O}'(Z\cup Q)=(1,4,4,3,4,4,\ldots,4,2,4)$ when $m$ is even; and \item $\mathcal{O}(Z\cup Q)=(1,3,4,4,\ldots,4,4,4)$ and $\mathcal{O}'(Z\cup Q)=(1,4,4,3,4,4,\ldots,4,4,4)$ when $m$ is odd;
\end{itemize}
as required.\\

\noindent Property $(3)$ follows by noting that each of the twenty-two edges in $E(Z)$ belong to distinct edge orbits of $K-I$ under $\rho$.\\

\noindent Property $(P1)$ follows by noting that 
\begin{itemize}
\item $\mathcal{O}(Z)=(1,2,2,3,0,0,\ldots,0,4)$; 
\item $\mathcal{O}'(Z)=(1,3,2,0,0,\ldots,0,4)$; 
\item $\mathcal{O}(Q)=(0,1,2,1,4,4,\ldots,4,2,0)$ when $m$ is even and $\mathcal{O}(Q)=(0,1,2,1,4,4,\ldots,4,4,0)$ when $m$ is odd; 
\item $\mathcal{O}'(Q)=(0,1,2,3,4,4,\ldots,4,2,0)$ when $m$ is even and $\mathcal{O}'(Q)=(0,1,2,3,4,4,\ldots,4,4,0)$ when $m$ is odd; and 
\end{itemize}
the edge orbits of the twelve edges in $E(Z)\cap(O_1\cup O_2\cup O_3\cup O'_1\cup O'_2)$ are distinct from the edge orbits of the seven edges in $E(Q)\cap(O_1\cup O_2\cup O_3\cup O'_1\cup O'_2)$.\\

\noindent Property $(P2)$ (with $\D'=\{Z\cup Q\}$) then follows from (I), (II) and the fact that $\{I\cup S_1\cup S_1'\cup R_1\cup R_1', C\cup C'\cup S\cup S'\cup R\cup R'\}$ is a decomposition of $I\cup Z\cup Q$.\qed

\begin{lemma}\label{alln_odd}
There is a decomposition of $K_{A_{n}\cup B_{n}\cup\{c\}}$ into Hamilton fragments for each positive integer $n$.
\end{lemma}

\proof 
For each $n\in\{1,2,\ldots,9\}$ a suitable decomposition is given in Section \ref{Section5}. Suppose then that $n\geq 10$.

Let $K=K_{A_{n}\cup B_{n}\cup\{c\}}$ and let $\rho$ be the permutation on $\Z_{n-1}$ which maps $v$ to $v+2$ for each $v\in\Z_{n-1}$. Observe that $\rho$ has order $n-1$ when $n$ is even, and order $(n-1)/2$ when $n$ is odd. Our aim is to show, for each value $n$, that there are edge-disjoint subgraphs $I$, $Z$ and $Q$ of $K$ which satisfy properties $(1)-(3)$ and $(P2)$ of Lemma \ref{SuffCond2}.
The problem now splits according to the parity of $n$.\\

\noindent{\bf Case 1: $n\geq 10$ is even}

Let $n=2m$ and let $I$ be the $2$-factor of $K$ with 
$$I=(a_\infty, b_\infty,c)\cup (a_0,b_1,a_2,b_3,\ldots,a_{2m-2},b_0,a_1,b_2,a_3,\ldots,b_{2m-2}).$$
Define 
\begin{itemize}
\item $O_c=\{a_ic,b_ic\mid i\in \Z_{2m-1}\}$;
\item $O_0=\{a_ib_i\mid i\in \Z_{2m-1}\}$; 
\item $O_1=\{a_ia_{i+1},b_ib_{i+1}\mid i\in \Z_{2m-1}\}$;
\item $O_j=\{a_ia_{i+j},b_ib_{i+j},a_{i+j}b_{i},a_ib_{i+j}\mid i\in \Z_{2m-1}\}$ for each $j=2,3,\ldots,m-1$; and
\item $O_\infty = \{a_ia_{\infty},b_ib_{\infty},a_{\infty}b_{i},a_ib_{\infty}\mid i\in \Z_{2m-1}\}$.
\end{itemize}
Observe that $\{O_c,O_0,O_1,\ldots,O_{m-1},O_\infty\}$ partitions $E(K-I)$, and that each $O_i$ is the union of one or more edge orbits of $K-I$ under $\rho$. Furthermore, for any subgraph $H$ of $K-I$ we define  
$$\mathcal{O}(H)=(|E(H)\cap O_c|,|E(H)\cap O_0|,|E(H)\cap O_1|,\ldots,|E(H)\cap O_{m-1}|,|E(H)\cap O_\infty|).$$
Observe that 
$$\mathcal{O}(K-I)=(2(2m-1),2m-1,2(2m-1),4(2m-1),4(2m-1),\ldots,4(2m-1)),$$
and that the orbit of $^{A}H$ under $\rho$ decomposes $^{A}(K-I)$, whenever $H$ is a subgraph of $K-I$ such that $\mathcal{O}(H)=(2,1,2,4,4,\ldots,4).$

Let $Z=C\cup C' \cup S\cup S'=C\cup C'\cup T\cup T'$ and $Q=R\cup R'$ where
\begin{itemize}
\item $C=[a_{m-1},a_{m}]\cup [b_{m},a_{m+2}]$; 
\item $C'=[a_m,b_m]\cup [a_{m+2},a_{m-1}]$;
\item $S=[a_\infty,a_0,c,b_1,a_{2m-5}]\cup [b_\infty,b_{2m-3},a_2,b_{2m-2}]$; 
\item $S'=[a_\infty,b_{2m-3},b_1,b_{2m-2}]\cup [b_\infty,a_0,a_2,a_{2m-5}]$;
\item $T=[a_\infty,b_{2m-3},b_1,a_{2m-5}]\cup [b_\infty,a_0,a_2,b_{2m-2}]$;
\item $T'=[a_\infty,a_0,c,b_1,b_{2m-2}]\cup [b_\infty,b_{2m-3},a_2,a_{2m-5}]$; 
\item $R=[a_{2m-2},a_4,a_{2m-4},a_6,\ldots, a_{m+4},a_{m-2},a_{m+2}]$\\
$\cup$ $[a_{2m-5},a_3,a_{2m-7},a_5,\ldots, a_{m+1},a_{m-3},a_{m-1}]$; and
\item $R'=[a_{2m-2},a_3,a_{2m-4},a_5,\ldots, a_{m+4},a_{m-3},a_{m+2}]$ \\
$\cup$ $[a_{2m-5},a_4,a_{2m-7},a_6,\ldots, a_{m+1},a_{m-2},a_{m-1}]$;
\end{itemize} 
when $m$ is even, and 
\begin{itemize}
\item $C=[b_{m-1},b_{m+1}]\cup [a_{m+1},b_{m-2}]$; 
\item $C'=[b_{m+1},a_{m+1}]\cup [b_{m-2},b_{m-1}]$;
\item $S=[a_\infty,a_0,c,b_1,b_{2m-5}]\cup [b_\infty,b_{2m-3},a_2,a_{2m-2}]$; 
\item $S'=[a_\infty,b_{2m-3},b_1,a_{2m-2}]\cup [b_\infty,a_0,a_2,b_{2m-5}]$;
\item $T=[a_\infty,b_{2m-3},b_1,b_{2m-5}]\cup [b_\infty,a_0,a_2,a_{2m-2}]$;
\item $T'=[a_\infty,a_0,c,b_1,a_{2m-2}]\cup [b_\infty,b_{2m-3},a_2,b_{2m-5}]$; 
\item $R=[a_{2m-2},a_4,a_{2m-4},a_6,\ldots,a_{m+3},a_{m-1}]$\\
$\cup$ $[a_{2m-5},a_3,a_{2m-7},a_5,\ldots,a_{m},a_{m-2}]$; and
\item $R'=[a_{2m-2},a_3,a_{2m-4},a_5,\ldots,a_{m+3},a_{m-2}]$ \\
$\cup$ $[a_{2m-5},a_4,a_{2m-7},a_6,\ldots, a_{m},a_{m-1}]$;
\end{itemize}
when $m$ is odd.

\noindent Property $(1)$ follows from Lemma \ref{4regular_odd}, setting $t=m$ and $H=Z\cup Q$ when $m$ is even and $t=m+1$ and $H=Z\cup Q$ when $m$ is odd, and noting in particular that $\{E(S),E(S')\}$ and $\{E(T),E(T')\}$ are $\{a_\infty,b_\infty\}$-connectors by Lemma \ref{connectors} $(6)$.\\

\noindent 
Property $(2)$ follows by noting that 
\begin{itemize}
\item $\mathcal{O}(C)= (0,0,0,1,1,0,0,0,0,\ldots,0,0)$;
\item $\mathcal{O}(C')=(0,1,1,0,0,0,0,0,0,\ldots,0,0)$;
\item $\mathcal{O}(S)= (2,0,0,0,1,1,1,0,0,\ldots,0,2)$;
\item $\mathcal{O}(S')=(0,0,0,2,1,0,0,1,0,\ldots,0,2)$;
\item $\mathcal{O}(R)= (0,0,0,1,0,2,1,2,2,\ldots,2,0)$;
\item $\mathcal{O}(R')=(0,0,1,0,1,1,2,1,2,\ldots,2,0)$;
\end{itemize}
and hence $\mathcal{O}(Z\cup Q)=
                       (2,1,2,4,4,\ldots,4)$ as required.\\

\noindent Property $(3)$ follows by noting that each of the twelve edges in $E(Z)$ belong to distinct edge orbits of $K-I$ under $\rho$.\\ 

\noindent Property $(P2)$
follows from Lemma \ref{12regular_odd}, setting
$H=I\cup (Z\cup Q_0)\cup (\rho^m(Z)\cup Q_m)$, 
$t=0$,  
$P_1=[a_0,b_1,a_2,b_3,\ldots,a_{2m-2},b_0]$,
$P_2=[b_0,a_1,b_2,a_3,\ldots,b_{2m-2},a_0]$,
$H_1=Z\cup Q_0$, and
$H_2=\rho^m(Z)\cup Q_m$, noting that
\begin{itemize}
\item from the proof of property $(1)$ above, for each $i\in\{1,2\}$ there is a \linebreak $t_i\in\{0,1,\ldots,n-2\}$ and decompositions $\{S_i,S_i',R_i,R_i',C_i,C_i'\}$ and \linebreak $\{T_i,T_i',R_i,R_i',C_i,C_i'\}$ of $H_i$, with  $(S_1,S_1',T_1,T_1')=(S,S',T,T')$ and \linebreak $(S_2,S_2',T_2,T_2')=(\rho^m(S),\rho^m(S'),\rho^m(T),\rho^m(T'))$, that satisfy the conditions of Lemma \ref{4regular_odd};
\item both $\{E(S),E(S'),E(X\cup P_1)\}$ and $\{E(T),E(T'),E(X\cup P_1)\}$ are $\{a_\infty,b_\infty,a_0\}$-connectors by Lemma \ref{connectors}$(7)$ with $(\alpha,\beta,u)=(a_\infty,b_\infty,a_0)$; and 
\item both $\{E(\rho^m(S)),E(\rho^m(S')),E(X\cup P_2)\}$ and $\{E(\rho^m(T)),E(\rho^m(T')),E(X\cup P_2)\}$ are $\{a_\infty,b_\infty,a_1\}$-connectors by Lemma \ref{connectors}$(7)$ with $(\alpha,\beta,u)=(a_\infty,b_\infty,a_1)$.
\end{itemize}

\noindent{\bf Case 2: $n\geq 11$ is odd}

Let $n=2m+1$ and let $I$ be the $2$-factor of $K$ with 
$$I=(a_\infty, b_\infty,c)\cup (a_0,b_1,a_2,b_3,\ldots,a_{2m-3}b_{2m-2})\cup (b_0,a_1,b_2,a_3,\ldots,a_{2m-3},b_{2m-2}).$$
Define 
\begin{itemize}
\item $O_c=\{a_ic,b_ic\mid \text{ $i$ even}, i\in \Z_{2m}\}$;
\item $O_0=\{a_ib_i\mid \text{ $i$ even}, i\in \Z_{2m}\}$; 
\item $O_1=\{a_ia_{i+1},b_ib_{i+1}\mid \text{ $i$ even}, i\in \Z_{2m}\}$; 
\item $O_j=\{a_ia_{i+j},b_ib_{i+j},a_{i+j}b_{i},a_ib_{i+j}\mid \text{ $i$ even}, i\in \Z_{2m}\}$ for each $j=2,3,\ldots,m$;
\item $O_\infty = \{a_ia_{\infty},b_ib_{\infty},a_{\infty}b_{i},a_ib_{\infty}\mid \text{ $i$ even}, i\in \Z_{2m}\}$;
\item $O'_c=\{a_ic,b_ic\mid \text{ $i$ odd}, i\in \Z_{2m}\}$;
\item $O'_0=\{a_ib_i\mid \text{ $i$ odd}, i\in \Z_{2m}\}$; 
\item $O'_1=\{a_ia_{i+1},b_ib_{i+1}\mid \text{ $i$ odd}, i\in \Z_{2m}\}$; 
\item $O'_j=\{a_ia_{i+j},b_ib_{i+j},a_{i+j}b_{i},a_ib_{i+j}\mid \text{ $i$ odd}, i\in \Z_{2m}\}$ for each $j=2,3,\ldots,m$; and
\item $O'_\infty = \{a_ia_{\infty},b_ib_{\infty},a_{\infty}b_{i},a_ib_{\infty}\mid \text{ $i$ odd}, i\in \Z_{2m}\}$;
\end{itemize}
Observe that $\{O_c,O_0,O_1,\ldots,O_{m},O_\infty,O'_c,O'_0,O'_1,\ldots,O'_{m},O'_\infty\}$ partitions $E(K-I)$ (with $O_m=O'_m$ when $m$ is odd), and that each $O_i$ and $O'_i$ is the union of one or more edge orbits of $K-I$ under $\rho$. Furthermore, for any subgraph $H$ of $K-I$ we define  
$$\mathcal{O}(H)=(|E(H)\cap O_c|,|E(H)\cap O_0|,|E(H)\cap O_1|,\ldots,|E(H)\cap O_{m}|,|E(H)\cap O_\infty|)$$ 
and 
$$\mathcal{O}'(H)=(|E(H)\cap O'_c|,|E(H)\cap O'_0|,|E(H)\cap O'_1|,\ldots,|E(H)\cap O'_{m}|,|E(H)\cap O'_\infty|).$$ 
Observe that
$$\mathcal{O}(K-I)=(2m,m,2m,4m,4m,\ldots,4m,2m,4m)$$
and 
$$\mathcal{O}'(K-I)=(2m,m,2m,4m,4m,\ldots,4m,2m,4m)$$
when $m$ is even, and that 
$$\mathcal{O}(K-I)=(2m,m,2m,4m,4m,\ldots,4m,4m,4m)$$
and 
$$\mathcal{O}'(K-I)=(2m,m,2m,4m,4m,\ldots,4m,4m,4m)$$
when $m$ is odd.
It follows that the orbit of $^{A}H$ under $\rho$ decomposes $^{A}(K-I)$, whenever $H$ is a subgraph of $K-I$ such that 
\begin{itemize}
\item $\mathcal{O}(H)=(2,1,2,4,4,\ldots,4,2,4)$ and $\mathcal{O}'(H)=(2,1,2,4,4,\ldots,4,2,4)$ when $m$ is even; and
\item $\mathcal{O}(H)=(2,1,2,4,4,\ldots,4,4,4)$ and $\mathcal{O}'(H)=(2,1,2,4,4,\ldots,4,4,4)$ when $m$ is odd.
\end{itemize}

For each $i\in\{1,2\}$ let $Z_i=C_i\cup C_i' \cup S_i\cup S_i'=C_i\cup C_i' \cup T_i\cup T_i'$ and $Q_i=R_i\cup R_i' $ where
\begin{itemize}
\item $C_1=[b_{m},b_{m+1}]\cup [b_{m+3},a_{m+1}]\cup [a_{m+2},a_m]$; 
\item $C'_1=[b_{m+1},b_{m+3}]\cup [a_{m+1},a_{m+2}]\cup [a_m,b_m]$;
\item $S_1=[a_\infty,b_1,c,a_2,a_{2m-2}]\cup [b_\infty,b_{0},a_3,a_{2m-1}]$; 
\item $S'_1=[a_\infty,b_{0},a_2,a_{2m-1}]\cup [b_\infty,b_1,a_3,a_{2m-2}]$;
\item $T_1=[a_\infty,b_{0},a_2,a_{2m-2}]\cup [b_\infty,b_1,a_3,a_{2m-1}]$;
\item $T'_1=[a_\infty,b_1,c,a_2,a_{2m-1}]\cup [b_\infty,b_{0},a_3,a_{2m-2}]$; 
\item $R_1=[a_{2m-1},a_4,a_{2m-3},a_6,\ldots, a_{m+5},a_{m-2},a_{m+3}]$\\
$\cup$ $[a_{2m-2},a_5,a_{2m-4},a_7,\ldots,a_{m+4},a_{m-1},a_{m+2}]$;
\item $R'_1=[a_{2m-1},a_5,a_{2m-3},a_7,\ldots, a_{m+5},a_{m-1},a_{m+3}]$ \\
$\cup$ $[a_{2m-2},a_4,a_{2m-4},a_6,\ldots, a_{m+4},a_{m-2},a_{m+2}]$;
\end{itemize}
and
\begin{itemize}
\item $C_2=[a_{m},a_{m+1}]\cup [b_{m+1},b_{m+2}]\cup [b_{m},b_{m+3}]$; 
\item $C'_2=[a_{m+1},b_{m+1}]\cup [b_{m+2},b_{m}]\cup [b_{m+3},a_m]$;
\item $S_2=[a_\infty,a_1,c,b_2,b_{2m-2}]\cup [b_\infty,a_{0},a_3,b_{2m-1}]$; 
\item $S'_2=[a_\infty,a_{0},b_2,b_{2m-1}]\cup [b_\infty,a_1,a_3,b_{2m-2}]$;
\item $T_2=[a_\infty,a_{0},b_2,b_{2m-2}]\cup [b_\infty,a_1,a_3,b_{2m-1}]$;
\item $T'_2=[a_\infty,a_1,c,b_2,b_{2m-1}]\cup [b_\infty,a_{0},a_3,b_{2m-2}]$; 
\item $R_2=[b_{2m-1},b_4,b_{2m-3},b_6,\ldots, b_{m+5},b_{m-2},b_{m+3}]$\\
$\cup$ $[b_{2m-2},b_5,b_{2m-4},b_7,\ldots,b_{m+4},b_{m-1},b_{m+2}]$;
\item $R'_2=[b_{2m-1},b_5,b_{2m-3},b_7,\ldots, b_{m+5},b_{m-1},b_{m+3}]$ \\
$\cup$ $[b_{2m-2},b_4,b_{2m-4},b_6,\ldots, b_{m+4},b_{m-2},b_{m+2}]$;
\end{itemize} 
when $m$ is even, and 
\begin{itemize}
\item $C_1=[a_4,a_{2m-2}]\cup [b_{2m-2},b_{2m-1}]\cup [a_{2m-1},b_5]$; 
\item $C'_1=[a_{4},a_{2m-1}]\cup [b_{2m-1},b_{5}]\cup [a_{2m-2},b_{2m-2}]$;
\item $S_1=[a_\infty,a_1,c,b_2,a_{2m-3}]\cup [b_\infty,b_{0},a_3,a_{2m-4}]$; 
\item $S'_1=[a_\infty,b_{0},b_2,a_{2m-4}]\cup [b_\infty,a_1,a_3,a_{2m-3}]$;
\item $T_1=[a_\infty,b_{0},b_2,a_{2m-3}]\cup [b_\infty,a_1,a_3,a_{2m-4}]$;
\item $T'_1=[a_\infty,a_1,c,b_2,a_{2m-4}]\cup [b_\infty,b_{0},a_3,a_{2m-3}]$; 
\item $C_2=[b_4,b_{2m-1}]\cup [a_{2m-1},a_{2m-2}]\cup [b_{2m-2},b_5]$; 
\item $C'_2=[b_{4},b_{2m-2}]\cup [a_{2m-2},b_{5}]\cup [a_{2m-1},b_{2m-1}]$;
\item $S_2=[a_\infty,b_1,c,a_2,b_{2m-3}]\cup [b_\infty,a_{0},a_3,b_{2m-4}]$; 
\item $S'_2=[a_\infty,a_{0},a_2,b_{2m-4}]\cup [b_\infty,b_1,a_3,b_{2m-3}]$;
\item $T_2=[a_\infty,a_{0},a_2,b_{2m-3}]\cup [b_\infty,b_1,a_3,b_{2m-4}]$;
\item $T'_2=[a_\infty,b_1,c,a_2,b_{2m-4}]\cup [b_\infty,a_{0},a_3,b_{2m-3}]$; 
\end{itemize}
when $m$ is odd, with 
$R_1=[a_7,a_5]\cup [a_6,a_4]$, 
$R'_1=[a_7,a_4]\cup [a_6,a_5]$,
$R_2=[b_7,b_5]\cup [b_6,b_4]$ and
$R'_2=[b_7,b_4]\cup [b_6,b_5]$ when $m=5$; and
\begin{itemize}
\item $R_1=[a_{2m-3},a_7,a_{2m-7},a_{11},\ldots, a_{m+2},a_{m},a_{m+4},a_{m-4},a_{m+8},a_{m-8},\ldots,a_{2m-5},a_5]$\\
$\cup$ $[a_{2m-4},a_6,a_{2m-8},a_{10},\ldots, a_{m+1},a_{m-1},a_{m+3},a_{m-5},a_{m+7},a_{m-9},\ldots,a_{2m-6},a_4]$;
\item $R'_1=[a_{2m-3},a_6,a_{2m-7},a_{10},\ldots, a_{m+2},a_{m-1},a_{m+4},a_{m-5},a_{m+8},a_{m-9},\ldots,a_{2m-5},a_4]$\\
$\cup$ $[a_{2m-4},a_7,a_{2m-8},a_{11},\ldots, a_{m+1},a_{m},a_{m+3},a_{m-4},a_{m+7},a_{m-8},\ldots,a_{2m-6},a_5]$;
\item $R_2=[b_{2m-3},b_7,b_{2m-7},b_{11},\ldots, b_{m+2},b_{m},b_{m+4},b_{m-4},b_{m+8},b_{m-8},\ldots,b_{2m-5},b_5]$\\
$\cup$ $[b_{2m-4},b_6,b_{2m-8},b_{10},\ldots, b_{m+1},b_{m-1},b_{m+3},b_{m-5},b_{m+7},b_{m-9},\ldots,b_{2m-6},b_4]$;
\item $R'_2=[b_{2m-3},b_6,b_{2m-7},b_{10},\ldots, b_{m+2},b_{m-1},b_{m+4},b_{m-5},b_{m+8},b_{m-9},\ldots,b_{2m-5},b_4]$\\
$\cup$ $[b_{2m-4},b_7,b_{2m-8},b_{11},\ldots, b_{m+1},b_{m},b_{m+3},b_{m-4},b_{m+7},b_{m-8},\ldots,b_{2m-6},b_5]$;
\end{itemize}
when $m\equiv 1\pmod{4}$ and $m\geq 9$; and
\begin{itemize}
\item $R_1=[a_{2m-3},a_7,a_{2m-7},a_{11},\ldots, a_{m+4},a_{m},a_{m+2},a_{m-2},a_{m+6},a_{m-6},\ldots,a_{2m-5},a_5]$\\
$\cup$ $[a_{2m-4},a_6,a_{2m-8},a_{10},\ldots, a_{m+3},a_{m-1},a_{m+1},a_{m-3},a_{m+5},a_{m-7},\ldots,a_{2m-6},a_4]$;
\item $R'_1=[a_{2m-3},a_6,a_{2m-7},a_{10},\ldots, a_{m+4},a_{m-1},a_{m+2},a_{m-3},a_{m+6},a_{m-7},\ldots,a_{2m-5},a_4]$\\
$\cup$ $[a_{2m-4},a_7,a_{2m-8},a_{11},\ldots, a_{m+3},a_{m},a_{m+1},a_{m-2},a_{m+5},a_{m-6},\ldots,a_{2m-6},a_5]$;
\item $R_2=[b_{2m-3},b_7,b_{2m-7},b_{11},\ldots, b_{m+4},b_{m},b_{m+2},b_{m-2},b_{m+6},b_{m-6},\ldots,b_{2m-5},b_5]$\\
$\cup$ $[b_{2m-4},b_6,b_{2m-8},b_{10},\ldots, b_{m+3},b_{m-1},b_{m+1},b_{m-3},b_{m+5},b_{m-7},\ldots,b_{2m-6},b_4]$;
\item $R'_2=[b_{2m-3},b_6,b_{2m-7},b_{10},\ldots, b_{m+4},b_{m-1},b_{m+2},b_{m-3},b_{m+6},b_{m-7},\ldots,b_{2m-5},b_4]$\\
$\cup$ $[b_{2m-4},b_7,b_{2m-8},b_{11},\ldots, b_{m+3},b_{m},b_{m+1},b_{m-2},b_{m+5},b_{m-6},\ldots,b_{2m-6},b_5]$;
\end{itemize}
when $m\equiv 3\pmod{4}$.\\

\noindent Property $(1)$ follows by setting $Z=Z_1\cup Z_2$ and $Q=Q_1\cup Q_2$, and observing that $Z_i\cup Q_i$ is a $Q_i$-adjustable Hamilton fragment for each $i\in\{1,2\}$, by Lemma \ref{4regular_odd}, setting $H=Z_i\cup Q_i$ and 
\begin{itemize}
\item $t=m$ when $i=1$ and $m$ is even;
\item $t=m+1$ when $i=2$ and $m$ is even;
\item $t=2m-2$ when $i=1$ and $m$ is odd;
\item $t=2m-1$ when $i=2$ and $m$ is odd;
\end{itemize}
and noting in particular that $\{E(S_i),E(S_i')\}$ and $\{E(T_i),E(T_i')\}$ are $\{a_\infty,b_\infty\}$-connectors by Lemma \ref{connectors}$(6)$.\\ 

\noindent 
Property $(2)$ follows by noting that
\begin{itemize}
\item $\mathcal{O}(Z_1\cup Z_2)=(2,1,2,4,4,2,2,0,0,\ldots,0,4)$; 
\item $\mathcal{O}'(Z_1\cup Z_2)=(2,1,2,4,2,2,0,0,\ldots,0,4)$;
\item $\mathcal{O}(Q_1\cup Q_2)= (0,0,0,0,0,2,2,4,4,\ldots,4,2,0)$;
\item $\mathcal{O}'(Q_1\cup Q_2)= (0,0,0,0,2,2,4,4,\ldots,4,2,0)$;
\end{itemize}
and hence $\mathcal{O}(Z\cup Q)=\mathcal{O}'(Z\cup Q)=(2,1,2,4,4,\ldots,4,2,4)$ when $m$ is even; 
\begin{itemize}
\item $\mathcal{O}(Z_1\cup Z_2)=(2,1,2,2,2,4,4,4)$ and $\mathcal{O}'(Z_1\cup Z_2)=(2,1,0,2,4,4,4,4)$ and
\item $\mathcal{O}(Q_1\cup Q_2)=(0,0,0,2,2,0,0,0)$ and $\mathcal{O}'(Q_1\cup Q_2)=(0,0,2,2,0,0,0,0)$;
\end{itemize}
and hence $\mathcal{O}(Z\cup Q)=\mathcal{O}'(Z\cup Q)=(2,1,2,4,4,4,4,4)$ when $m=5$;
\begin{itemize}
\item $\mathcal{O}(Z_1\cup Z_2)=(2,1,2,2,2,0,0,4,4,4)$ and $\mathcal{O}'(Z_1\cup Z_2)=(2,1,0,2,0,0,4,4,4,4)$; and
\item $\mathcal{O}(Q_1\cup Q_2)=(0,0,0,2,2,4,4,0,0,0)$ and $\mathcal{O}'(Q_1\cup Q_2)=(0,0,2,2,4,4,0,0,0,0)$;
\end{itemize}
and hence $\mathcal{O}(Z\cup Q)=\mathcal{O}'(Z\cup Q)=(2,1,2,4,4,4,4,4,4,4)$ 
when $m=7$; and
\begin{itemize}
\item $\mathcal{O}(Z_1\cup Z_2)=(2,1,2,2,2,0,0,4,4,0,0,\ldots,0,4)$; 
\item $\mathcal{O}'(Z_1\cup Z_2)=(2,1,0,2,0,0,4,4,0,0,\ldots,0,4)$;
\item $\mathcal{O}(Q_1\cup Q_2)=(0,0,0,2,2,4,4,0,0,4,4,\ldots,4,0)$; 
\item $\mathcal{O}'(Q_1\cup Q_2)=(0,0,2,2,4,4,0,0,4,4,\ldots,4,0)$;
\end{itemize}
and hence $\mathcal{O}(Z\cup Q)=\mathcal{O}'(Z\cup Q)=(2,1,2,4,4,\ldots,4)$ when $m\geq 9$ is odd.\\ 

\noindent Property $(3)$ follows by noting that each of the twelve edges in $E(Z)$ belong to distinct edge orbits of $K-I$ under $\rho$.\\

\noindent Finally, we prove that property $(P2)$ is satisfied with $\D'=\{Z\cup Q_0\}$ (recall that $Q_0$ satifies $^A Q_0=\prescript{A}{}Q$ and that $Q=Q_1\cup Q_2$ with $Q_i= R_i\cup R_i'$ for $i=1,2$). 
Let $\{\hat{Q}_1,\hat{Q}_2\}$ be a decomposition of $Q_0$ satisfying $^A \hat{Q}_1=\prescript{A}{}Q_1$ and $^A \hat{Q}_2=\prescript{A}{}Q_2$, and for each $i\in\{1,2\}$ let $\{\hat{R}_i,\hat{R}'_i\}$ be a decomposition of $\hat{Q}_i$ satisfying $^A \hat{R}_i=\prescript{A}{}R_i$ and $^A \hat{R}'_i=\prescript{A}{}R'_i$. Then, when $m$ is even, property $(P2)$
follows from Lemma \ref{12regular_odd}, setting
\begin{itemize}
\item $H=I\cup (Z\cup Q_0)$; 
\item $t=m+1$;  
\item $P_1=[a_{m+1},a_{m},b_{m-1},a_{m-2},\ldots,a_{m+2},b_{m+1}]$;
\item $P_2=[a_{m+1},b_{m},a_{m-1},b_{m-2},\ldots,b_{m+2},b_{m+1}]$;
\item $H_1=Z_1\cup \hat{Q}_1$;
and
\item $H_2=(Z_2-\{a_{m}a_{m+1},b_{m+1}b_{m+2}\}+\{a_{m}b_{m+1},a_{m+1}b_{m+2}\})\cup \hat{Q}_2$;
\end{itemize} noting that
\begin{itemize}
\item from the proof of property $(1)$ above, it is easy to see that for $t_1=m$ and $t_2=m+1$ there are decompositions
$\{S_1,S_1',\hat{R}_1,\hat{R}_1',C_1,C_1'\}$ and $\{T_1,T_1',\hat{R}_1,\hat{R}_1',C_1,C_1'\}$ of $H_1$, and 
$\{S_2,S_2',\hat{R}_2,\hat{R}_2',C_2-\{a_{m}a_{m+1},b_{m+1}b_{m+2}\}+\{a_{m}b_{m+1},a_{m+1}b_{m+2}\},C_2'\}$ and $\{T_2,T_2',\hat{R}_2,\hat{R}_2',C_2-\{a_{m}a_{m+1},b_{m+1}b_{m+2}\}+\{a_{m}b_{m+1},a_{m+1}b_{m+2}\},C_2'\}$ of $H_2$ that satisfy the conditions of Lemma \ref{4regular_odd};  
\item both $\{E(S_1),E(S_1'),E(X\cup P_1)\}$ and $\{E(T_1),E(T_1'),E(X\cup P_1)\}$ are $\{a_\infty,b_\infty,b_1\}$-connectors by Lemma \ref{connectors}$(7)$ with $(\alpha,\beta,u)=(a_\infty,b_\infty,b_1)$; and 
\item both $\{E(S_2),E(S'_2),E(X\cup P_2)\}$ and $\{E(T_2),E(T'_2),E(X\cup P_2)\}$ are $\{a_\infty,b_\infty,a_1\}$-connectors by Lemma \ref{connectors}$(7)$ with $(\alpha,\beta,u)=(a_\infty,b_\infty,a_1)$.
\end{itemize}

\noindent When $m$ is odd, property $(P2)$
follows from Lemma \ref{12regular_odd}, setting
\begin{itemize}
\item $H=I\cup (Z\cup Q_0)$; 
\item $t=m+1$; 
\item $P_1=[a_{m+1},a_{m},b_{m-1},a_{m-2},\ldots,a_{m+2},b_{m+1}]$;
\item $P_2=[b_{m+1},b_{m},a_{m-1},b_{m-2},\ldots,b_{m+2},a_{m+1}]$;
\item $H_1=Z_1\cup \hat{Q}_1'-\{e_1\}+\{e_1'\}$; 
and
\item $H_2=Z_2\cup \hat{Q}_2'-\{e_2\}+\{e_2'\}$;
\end{itemize} where $e_1\in E(\hat{R}_1')$ and $e_2\in E(\hat{R}_2')$ with $\{e_1,e_2\}=\{a_ma_{m+1},b_mb_{m+1}\}$ and $\{e_1',e_2'\}=\{a_mb_{m+1},b_{m}a_{m+1}\}$, noting that
\begin{itemize}
\item it follows easily from the proof of property $(1)$ above, that for each $i\in\{1,2\}$ the decompositions $\{S_i,S_i',\hat{R}_i,\hat{R}_i'-\{e_i\}+\{e_i'\},C_1,C_1'\}$ and $\{T_i,T_i',\hat{R}_i,\hat{R}_i'-\{e_i\}+\{e_i'\},C_i,C_i'\}$ of $H_i$ satisfy the conditions of Lemma \ref{4regular_odd};  
\item both $\{E(S_1),E(S_1'),E(X\cup P_1)\}$ and $\{E(T_1),E(T_1'),E(X\cup P_1)\}$ are $\{a_\infty,b_\infty,a_1\}$-connectors by Lemma \ref{connectors}$(7)$ with $(\alpha,\beta,u)=(a_\infty,b_\infty,a_1)$; and 
\item both $\{E(S_2),E(S'_2),E(X\cup P_2)\}$ and $\{E(T_2),E(T'_2),E(X\cup P_2)\}$ are $\{a_\infty,b_\infty,b_1\}$-connectors by Lemma \ref{connectors}$(7)$ with $(\alpha,\beta,u)=(a_\infty,b_\infty,b_1)$.
\end{itemize}
\qed

\section{Decompositions into Hamilton Fragments for Small Degree Cases}\label{Section5}

\subsection{Decompositions of $K_{A_n\cup B_n}$ for $n\in\{2,3,4,5,7\}$}

The proof splits into cases according to the value of $n$.\\

\noindent{\bf The case n=2:}

Let  $H=K_{A_{2}\cup B_{2}}$, let $G$ be a $4$-regular graph and let $\{\F_0,\F_\infty\}$ be a directed $2$-factorisation of $G$ where $\F_\infty$ is a Hamilton cycle, say $\F_\infty=(v_1,v_2,\ldots,v_m)$. Our aim is to show that $L(G)$ decomposes into Hamilton cycles (and hence $H$ is a Hamilton fragment), and we do this by applying Lemma \ref{manyrepairs}
with $H_i=\sigma_{v_i}(H)$ and $V_i=\{a_\infty^{v_i},b_\infty^{v_i}\}$ for each $i=1,2,\ldots,m$.
To this end, observe that $\sigma_{v_1}(H),\sigma_{v_2}(H),\ldots,\sigma_{v_m}(H)$ are edge-disjoint subgraphs of $L(G)$,
and since $b_\infty^{v_i}=a_\infty^{v_{i+1}}$ for $i=1,2,\ldots,m-1$ it follows that $V_i\cap V_{i+1}\ne \emptyset$ for $i=1,2,\ldots,m-1$ as required. It remains to show there is a $2$-factorisation $\{J_1,J_2,J_3\}$ of $L(G)$ such that 
\begin{itemize}
\item $\S_\infty$ links $\{J_1,J_2,J_3\}$; and
\item $\sigma_{v_i}(H)$ induces an $\{a_\infty^{v_{i}},b_\infty^{v_i}\}$-connector in $\{J_1,J_2,J_3\}$, for $i=1,2,\ldots,m$.
\end{itemize} 

Let $\{X,Y,Y'\}$ and $\{X',Z,Z'\}$ be decompositions of $H$ defined by
\begin{itemize}
\item $E(X)=\{a_\infty b_\infty,a_0b_0\}$, $E(Y)=\{a_\infty a_0,a_0b_\infty\}$, $E(Y')=\{a_\infty b_0, b_0b_\infty\}$; and
\item $E(X')=\{a_\infty a_0, b_\infty b_0\}$, $E(Z)=\{a_\infty b_0,b_0a_0,a_0b_\infty\}$, $E(Z')=\{a_\infty b_\infty\}$;
\end{itemize}
let $U$ and $U'$ be disjoint subsets of $V(G)$ such that each of $U$ and $U'$ contains precisely one vertex from each connected component of $F_0$ and each connected component of $\F_0$ contains an edge oriented from $u$ to $v$ for some $u\in U$ and $v\in U'$, and let $\{J_1,J_2,J_3\}$ be the decomposition of $L(G)$ defined by
\begin{itemize}
\item $J_1=(\bigcup_{v\in V(G)\setminus( U\cup U')}\sigma_v(X))\cup (\bigcup_{v\in (U\cup U')}\sigma_v(X'))$;
\item $J_2=(\bigcup_{v\in V(G)\setminus( U\cup U')}\sigma_v(Y))\cup (\bigcup_{v\in (U)}\sigma_v(Z))\cup (\bigcup_{v\in (U')}\sigma_v(Z'))$;
\item $J_3=(\bigcup_{v\in V(G)\setminus( U\cup U')}\sigma_v(Y'))\cup (\bigcup_{v\in (U)}\sigma_v(Z'))\cup (\bigcup_{v\in (U')}\sigma_v(Z))$.
\end{itemize}   
It is a simple exercise to check that $\{J_1,J_2,J_3\}$ is a $2$-factorisation of $L(G)$ and that $\S_\infty$ links  $\{J_1,J_2,J_3\}$. Finally, it is easy to see that both $\{X,Y,Y'\}$ and $\{X',Z,Z'\}$ are $\{a_\infty,b_\infty\}$-connectors as required. The result follows.\\

\noindent{\bf The case n=3:}

Observe that $H=K_{A_{3}\cup B_{3}}$ is a Hamilton fragment by Lemma \ref{4regular_withHam}, setting $(r,s,t)=(5,0,1)$, $E(X_1)=\{b_\infty a_0,b_0a_1,b_1a_\infty\}$, 
$E(X_2)= \{a_\infty b_\infty,a_0b_0,a_1b_1\}$, 
$E(X_3)=\{a_\infty a_0, a_0b_1,b_1b_\infty\}$, 
$E(X_4)=\{a_\infty b_0, a_0a_1,a_1b_\infty\}$, and 
$E(X_5)=\{a_\infty a_1, b_1b_0,b_0b_\infty\}$, 
noting in particular that $\{E(X_1),E(X_2),E(X_3),E(X_4),E(X_5)\}$ is an $\{a_\infty,b_\infty\}$-connector by Lemma \ref{connectors} $(5)$, with $(\alpha,\beta,u,v,w,x)=(a_\infty,b_\infty,a_0,b_0,a_1,b_1)$.\\

\noindent{\bf The case n=4:}

Let $\{H_0,H_1\}$ be the decomposition of $K_{A_4\cup B_4}$ defined by
\begin{itemize}
\item $E(H_0)=\{a_\infty b_2,b_2b_0,a_0b_1,a_1 b_\infty\}\cup \{a_\infty b_\infty,a_0b_0,a_1b_2,b_2b_1\}$\\
$\cup \{a_\infty a_1, a_1b_0,a_0a_2,b_2b_\infty\}$;
\item $E(H_1)=\{a_\infty b_1, a_1a_0,a_0b_2,a_2b_\infty\}\cup \{a_\infty a_0,a_0 b_\infty, a_1b_1,a_2b_2\}$\\
$\cup \{a_\infty b_0, b_0a_2,a_2a_1,b_1b_\infty\}\cup \{a_\infty a_2,a_2b_1,b_1b_0,b_0b_\infty\}$. 
\end{itemize}
Then 
\begin{itemize}
\item $H_0$ is a Hamilton fragment by Lemma \ref{4regular_withHam}, setting $(r,s,t)=(3,0,1)$, \linebreak $E(X_1)=\{a_\infty b_2,b_2b_0,a_0b_1,a_1 b_\infty\}$, $E(X_2)= \{a_\infty b_\infty,a_0b_0,a_1b_2,b_2b_1\}$ and $E(X_3)=\{a_\infty a_1, a_1b_0,a_0a_2,b_2b_\infty\}$, noting in particular that $\{E(X_1),E(X_3),E(X_2)\}$ is an $\{a_\infty,b_\infty\}$-connector by Lemma \ref{connectors} $(2)$, with $(\alpha,\beta,u,v,w)=(a_\infty,b_\infty,b_2,b_0,a_1)$; and
\item $H_1$ is a Hamilton fragment by Lemma \ref{4regular_withHam}, setting $(r,s,t)=(4,1,2)$, \linebreak $E(X_1)=\{a_\infty b_1, a_1a_0,a_0b_2,a_2b_\infty\}$, $E(X_2)=\{a_\infty a_0,a_0 b_\infty, a_1b_1,a_2b_2\}$, $E(X_3)=\{a_\infty b_0, b_0a_2,a_2a_1,b_1b_\infty\}$  and $E(X_4)=\{a_\infty a_2,a_2b_1,b_1b_0,b_0b_\infty\}$, noting in particular that $\{E(X_4),E(X_3),E(X_1),E(X_2)\}$ is an $\{a_\infty,b_\infty\}$-connector by Lemma \ref{connectors} $(4)$, with $(\alpha,\beta,u,v,w,x)=(b_\infty,a_\infty,b_0,b_1,a_2,a_0).$
\end{itemize}

\noindent{\bf The case n=5:}

Let $\{H_0,H_1,H_2,H_3\}$ be the decomposition of $K_{A_5\cup B_5}$ defined by
\begin{itemize}
\item $E(H_0)=\{a_\infty b_3,b_3b_1,a_1a_2,b_2a_0,a_0b_\infty\} \cup \{a_\infty a_0,a_0b_3,b_3b_\infty,a_1b_1,a_2b_2\}$;
\item $E(H_1)=\{a_\infty a_2,a_2b_3,a_3b_0,a_0b_1,b_1b_\infty\} \cup \{a_\infty b_1,b_1a_2,a_2b_\infty,a_0b_0,a_3b_3\}$;
\item $E(H_2)=\{a_\infty a_1,a_1b_2,a_2a_0,a_0a_3,a_3b_\infty\}\cup\{a_\infty a_3,a_3a_2,a_2b_0,a_0a_1,a_1b_\infty\}$;
\item $E(H_3)=\{a_\infty b_0,b_0b_1,a_1b_3,b_3b_2,b_2b_\infty\}\cup \{a_\infty b_\infty,a_3b_1,b_1b_2,b_2b_0,b_0b_3\}$\\
$\cup \{a_\infty b_2,b_2a_3,a_3a_1,a_1b_0,b_0b_\infty\}$.
\end{itemize}
Then 
\begin{itemize}
\item $H_0$ is a Hamilton fragment by Lemma \ref{4regular_withHam}, setting $(r,s,t)=(2,1,2)$, 
 $E(X_1)= \{a_\infty b_3,b_3b_1,a_1a_2,b_2a_0,a_0b_\infty\} $, and $E(X_2)= \{a_\infty a_0,a_0b_3,b_3b_\infty,a_1b_1,a_2b_2\}$, noting that $\{E(X_1), E(X_2)\}$ is an $\{a_\infty, b_\infty\}$-connector by Lemma \ref{connectors} $(0)$ with $(\alpha,\beta,u,v,w)=(a_\infty,b_\infty,a_0,b_2,b_3)$;
 \item $H_1$ is a Hamilton fragment by Lemma \ref{4regular_withHam}, setting $(r,s,t)=(2,0,3)$, 
 $E(X_1)= \{a_\infty a_2,a_2b_3,a_3b_0,a_0b_1,b_1b_\infty\} $, and $E(X_2)= \{a_\infty b_1,b_1a_2,a_2b_\infty,a_0b_0,a_3b_3\}$, noting that $\{E(X_1), E(X_2)\}$ is an $\{a_\infty, b_\infty\}$-connector by Lemma \ref{connectors} $(0)$ with $(\alpha,\beta,u,v,w)=(a_\infty,b_\infty,b_1,a_0,a_2)$;
 \item $H_2$ is a Hamilton fragment by Lemma \ref{4regular_easy}, setting $E(S)=\{a_\infty a_1,a_1b_2,a_2a_0,a_0a_3,a_3b_\infty\}$, $E(S')=\{a_\infty a_3,a_3a_2,a_2b_0,a_0a_1,a_1b_\infty\}$ and $E(R)=E(R')=\emptyset$, noting that $\{E(S), E(S')\}$ is an $\{a_\infty, b_\infty\}$-connector by Lemma \ref{connectors} $(1)$, with $(\alpha,\beta,u,v,w)=(b_\infty,a_\infty,a_3,a_0,a_1)$; and
 \item $H_3$ is a Hamilton fragment by Lemma \ref{6regular_2}, setting $(s,t)=(1,3)$, 
$E(X_1)= \{a_\infty b_0,b_0b_1,a_1b_3,b_3b_2,b_2b_\infty\}$, 
$E(X_2)=\{a_\infty b_\infty,a_3b_1,b_1b_2,b_2b_0,b_0b_3\}$, \newline 
$E(X_1')=\{a_\infty b_0,,b_0b_3,a_3b_1,b_1b_2,b_2b_\infty\}$, 
$E(X_2')=\{a_\infty b_\infty,a_1b_3,b_3b_2,b_2b_0,b_0b_1\}$ and
$E(X_3)= \{a_\infty b_2,b_2a_3,a_3a_1,a_1b_0,b_0b_\infty\}$, noting that $\{E(X_1), E(X_2), E(X_3)\}$ and $\{E(X_1'), E(X_2'), E(X_3)\}$ are $\{a_\infty, b_\infty\}$-connectors by Lemma \ref{connectors} $(3)$, with $(\alpha,\beta,u,v,w,x)=(a_\infty,b_\infty,b_0,b_1,b_2,b_3)$ and $(\alpha,\beta,u,v,w,x)=(a_\infty,b_\infty,b_0,b_3,b_2,b_1)$, respectively. 
\end{itemize}

\noindent{\bf The case n=7:}

Let $\rho$ be the permutation on $K_{A_7\cup B_7}$ defined by 
$$\rho = (a_\infty) (b_\infty) (a_0\,\, a_2\,\,a_4) (a_1\,\, a_3\,\,a_5) (b_0\,\, b_2\,\,b_4) (b_1\,\, b_3\,\,b_5).$$
Let $\{H_0,H_1,H_2,H'_0,H'_1,H'_2,I\}$ be the decomposition of $K_{A_7\cup B_7}$ defined by
\begin{itemize}
\item $E(H_0)=\{a_\infty a_0,a_0a_2,a_2b_5,a_5a_4,b_4b_3,a_3b_1,b_1b_\infty\}\cup 
\{a_\infty b_1,b_1a_2,a_2b_3,a_3a_0,a_0b_\infty,a_4b_4, a_5b_5\}$;
\item $H_1=\rho(H_0)$;
\item $H_2=\rho(H_1)$;
\item $E(H_0')=\{a_\infty b_0,b_0b_2,b_2a_4,b_4b_5,b_5b_3,b_3a_1,a_1b_\infty\}\cup
\{a_\infty a_1,a_1b_2,b_2b_5,a_5a_3,a_3a_4,a_4b_0,b_0b_\infty\}$; and
\item $H_1'=\rho(H_0')$;
\item $H_2'=\rho(H_1')$;
\item $E(I)=\{a_\infty b_\infty,b_0a_1,a_1b_4,b_4a_5,a_5b_2,b_2a_3,a_3b_0\}$.
\end{itemize}
Observe that 
\begin{itemize}
\item $H_0$ is a Hamilton fragment by Lemma \ref{4regular}, setting $T=\{4,5\}$,
$E(S)=\{a_\infty a_0,a_0a_2, b_1b_\infty \}$,
$E(S')=\{a_\infty b_1,b_1a_2, a_0b_\infty\}$,
$E(R)=\{a_3b_1\}$,
$E(R')=\{a_3a_0\}$,
$E(C)=\{b_3b_4,a_4a_5,b_5a_2\}$, and
$E(C')=\{a_2b_3, b_4a_4,a_5b_5\}$,
noting in particular that $\{E(S),E(S')\}$ is an $\{a_\infty,b_\infty\}$-connector by Lemma \ref{connectors} $(1)$, with $(\alpha,\beta,u,v,w)=(a_\infty,b_\infty,a_0,a_2,b_1)$;
\item $H_0'$ is a Hamilton fragment by Lemma \ref{4regular_easy}, setting 
$E(S)=\{a_\infty b_0,b_0b_2,a_1b_\infty\}$,
$E(S')=\{a_\infty a_1,a_1b_2,b_0b_\infty\}$,
$E(R)=\{b_2a_4,b_4b_5,b_5b_3,b_3a_1\}$, and \newline
$E(R')=\{b_2b_5,a_5a_3,a_3a_4,a_4b_0\}$, noting in particular that $\{E(S),E(S')\}$ is an $\{a_\infty,b_\infty\}$-connector by Lemma \ref{connectors} $(1)$, with $(\alpha,\beta,u,v,w)=(a_\infty,b_\infty,b_0,b_2,a_1)$;
\item $H_0'\cup I$ is a Hamilton fragment by Lemma \ref{6regular_2}, setting $(s,t)=(4,5)$,\newline
$E(X_1)=\{a_\infty b_0, b_0b_2, b_2a_4, b_4b_5, b_5b_3, b_3a_1, a_1b_\infty\}$,\newline
$E(X_2)=\{a_\infty b_\infty, b_0a_1, a_1b_2, b_2b_5, a_5a_3, a_3a_4, a_4b_0\}$,\newline
$E(X_1')=\{b_\infty a_1, a_1b_2, b_2b_5, a_5a_3, a_3a_4, a_4b_0, b_0 a_\infty\}$,\newline
$E(X_2')=\{a_\infty b_\infty, b_0a_1 , b_0b_2, b_2a_4, b_4b_5, b_5b_3, b_3a_1 \}$, and\newline
$E(X_3)=\{a_\infty a_1,a_1b_4,b_4a_5,a_5b_2,b_2a_3,a_3b_0,b_0b_\infty\}$. \newline
To see that $\{E(X_1),E(X_2), E(X_3)\}$ and $\{E(X'_1),E(X'_2), E(X_3)\}$ are  $\{a_\infty, b_\infty\}$-connectors, suppose $\{F_1,F_2, F_3\}$ is a $2$-factorisation of 
some $6$-regular graph $G$ such that $E(X_1)\subseteq E(F_1)$, $E(X_2)\subseteq E(F_2)$ and $E(X_3)\subseteq E(F_3)$, or $E(X'_1)\subseteq E(F'_1)$, $E(X'_2)\subseteq E(F'_2)$ and $E(X_3)\subseteq E(F_3)$, respectively.
Note that $a_\y$ and $b_\y$ are currently in the same component of $F_2$ and in the same component of $F_3$. If $a_\y$ and $b_\y$ are also in the same component of $F_1$ then we are done. Otherwise, replace $F_1$ with $F_1 - \{a_\infty b_0, a_1 b_\infty\} + \{a_\infty a_1, b_0 b_\infty\}$ and replace $F_3$ with 
$F_3 - \{a_\infty a_1, b_0 b_\infty\} + \{a_\infty b_0, a_1 b_\infty\}$.

\end{itemize}
It follows that $H_1,H_2,H_1',H_2'$ are Hamilton fragments and thus $\{H_0,H_1,H_2,H'_0\cup I,H'_1,H'_2\}$ is the required decomposition of $K_{A_7\cup B_7}$.
\qed

\subsection{Decompositions of $K_{A_n\cup B_n\cup\{c\}}$ for $n\in\{1,2,\ldots,9\}$}

In this section we give a decomposition of $K_{A_n\cup B_n\cup\{c\}}$ into Hamilton fragments for each $n\in\{1,2,\ldots,9\}$.

\noindent{\bf The case n=1:}

Kotzig \cite{Kot} showed that a $3$-regular graph is Hamiltonian if and only if its line graph is Hamilton decomposable,
and it follows from this that $K_{A_1\cup B_1\cup\{c\}}$ itself is a Hamilton fragment.

For each of the cases $n=2,3,\ldots,9$, 
let $G$ be a $2n+1$-regular graph with vertex set $\{v_1,\ldots,v_m\}$, 
let $F$ be a $1-$factor of $G$ and let $\{\F_0,\F_1,\ldots,\F_{n-2},\F_\infty\}$ be a directed $2$-factorisation of $G-F$ where $\F_\infty$ is a Hamilton cycle, say $\F_\infty=(v_1,v_2,\ldots,v_m)$.

\noindent{\bf The case n=2:}

Let $H=K_{A_2\cup B_2\cup\{c\}}\cong K_5$. We show that $H$ is a Hamilton fragment. 
Our aim is to show that the subgraph $$L=\bigcup_{v\in V(G)}\sigma_v(H)$$
of $L(G)$ decomposes into Hamilton cycles, and we do this by applying Lemma \ref{manyrepairs} 
with $H_i=\sigma_{v_i}(H)$ and $V_i=\{a_\infty^{v_i},b_\infty^{v_i}\}$ for each $i=1,2,\ldots,m$.
To this end, observe that \linebreak $\sigma_{v_1}(H),\sigma_{v_2}(H),\ldots,\sigma_{v_m}(H)$ are edge-disjoint subgraphs of $L$,
and since $b_\infty^{v_i}=a_\infty^{v_{i+1}}$ for $i=1,2,\ldots,m-1$ it follows that $V_i\cap V_{i+1}\ne \emptyset$ for $i=1,2,\ldots,m-1$ as required. It remains to show there is a $2$-factorisation
$\{J_1,J_2,J_3,J_4\}$ of $L$ such that 
\begin{itemize}
\item $\{a_\infty^{v_1},a_\infty^{v_2},\ldots,a_\infty^{v_m}\}$ links $\{J_1,J_2,J_3,J_4\}$; and
\item $\sigma_{v_i}(H)$ induces an $\{a_\infty^{v_{i}},b_\infty^{v_i}\}$-connector in $\{J_1,J_2,J_3,J_4\}$, for $i=1,2,\ldots,m$.
\end{itemize} 
Let $\{U,U'\}$ be a partition of $V(G)$ such that both $U$ and $U'$ link $\{F,\F_0\}$ (such a partition exists by Lemma \ref{intersects}), 
and let $\{J_1,J_2,J_3,J_4\}$ be the decomposition of $L$ defined by 
\begin{itemize}
\item $J_1=\bigcup_{v\in U}\sigma_v([a_\y,a_0,c,b_\y])\cup \bigcup_{v\in U'}\sigma_v([a_\y,a_0,b_\y])$;
\item $J_2= \bigcup_{v\in U}\sigma_v([a_\y,c,b_0]\cup[b_\y,a_0])\cup \bigcup_{v\in U'}\sigma_v([a_\y,b_\y]\cup[a_0,b_0])$;
\item $J_3=\bigcup_{v\in U}\sigma_v([a_\y,b_0,b_\y])\cup \bigcup_{v\in U'}\sigma_v([a_\y,c,b_0,b_\y])$;
\item $J_4= \bigcup_{v\in U}\sigma_v([a_\y,b_\y]\cup[a_0,b_0])\cup \bigcup_{v\in U'}\sigma_v([a_\y,b_0]\cup[b_\y,c,a_0])$.
\end{itemize}
It is easily checked that $\{J_1,J_2,J_3,J_4\}$ is a $2$-factorisation of $L$ and that $\{a_\infty^{v_1},a_\infty^{v_2},\ldots,a_\infty^{v_m}\}$ links $\{J_1,J_2,J_3,J_4\}$.
It remains to show that for $i=1,2,\ldots,m$, $\sigma_{v_i}(H)$ induces an $\{a_\infty^{v_{i}},b_\infty^{v_i}\}$-connector in $\{J_1,J_2,J_3,J_4\}$. First note that $a_\y^{v_i}$ and $b_\y^{v_i}$ are in the same component of $J_1$ and $J_3$ for $i=1,2,\ldots,m$.
Also note that $a_\y^{v_i}$ and $b_\y^{v_i}$ are in the same component of $J_4$ for each $v_i\in U$, and that
$a_\y^{v_i}$ and $b_\y^{v_i}$ are in the same component of $J_2$ for each $v_i\in U'$.

If $a_\y^{v_i}$ and $b_\y^{v_i}$ are in the distinct components of $J_4$ for some $v_i\in U'$, then we move the two edges 
$a_\y^{v_i}b_0^{v_i}$ and $b_\y^{v_i}c^{v_i}$ from $J_4$ to $J_3$ and move the two edges 
$a_\y^{v_i}c^{v_i}$ and $b_\y^{v_i}b_0^{v_i}$ from $J_3$ to $J_4$. 
If $a_\y^{v_i}$ and $b_\y^{v_i}$ are in the distinct components of $J_2$ for some $v_i\in U$, then we move the two edges 
$a_\y^{v_i}c^{v_i}$ and $b_\y^{v_i}a_0^{v_i}$ from $J_2$ to $J_1$ and move the two edges 
$a_\y^{v_i}a_0^{v_i}$ and $b_\y^{v_i}c^{v_i}$ from $J_1$ to $J_2$. 
It is easily checked that this has the desired effect, and thus 
$\sigma_{v_i}(H)$ induces an $\{a_\infty^{v_{i}},b_\infty^{v_i}\}$-connector in $\{J_1,J_2,J_3,J_4\}$.

\noindent{\bf The case n=3:}

Let $H=K_{A_3\cup B_3\cup\{c\}}\cong K_7$. We show that $H$ is a Hamilton fragment. 
Our aim is to show that the subgraph $$L=\bigcup_{v\in V(G)}\sigma_v(H)$$
of $L(G)$ decomposes into Hamilton cycles, and we do this by applying Lemma \ref{manyrepairs} 
with $H_i=\sigma_{v_i}(H)$ and $V_i=\{a_\infty^{v_i},b_\infty^{v_i}\}$ for each $i=1,2,\ldots,m$.
To this end, observe that \linebreak $\sigma_{v_1}(H),\sigma_{v_2}(H),\ldots,\sigma_{v_m}(H)$ are edge-disjoint subgraphs of $L$,
and since $b_\infty^{v_i}=a_\infty^{v_{i+1}}$ for $i=1,2,\ldots,m-1$ it follows that $V_i\cap V_{i+1}\ne \emptyset$ for $i=1,2,\ldots,m-1$ as required. It remains to show there is a $2$-factorisation
$\{J_1,J_2,J_3,J_4,J_5,J_6\}$ of $L$ such that 
\begin{itemize}
\item $\{a_\infty^{v_1},a_\infty^{v_2},\ldots,a_\infty^{v_m}\}$ links $\{J_1,J_2,J_3,J_4,J_5,J_6\}$; and
\item $\sigma_{v_i}(H)$ induces an $\{a_\infty^{v_{i}},b_\infty^{v_i}\}$-connector in $\{J_1,J_2,J_3,J_4,J_5,J_6\}$, for $i=1,2,\ldots,m$.
\end{itemize} 
Let $\{U,U'\}$ be a partition of $V(G)$ such that both $U$ and $U'$ link $\{F,\F_0\}$,
let $\{V,V'\}$ be a partition of $V(G)$ such that both $V$ and $V'$ link $\{F,\F_1\}$ 
(such partitions exist by Lemma \ref{intersects}),
and let $\{J_1,J_2,J_3,J_4,J_5,J_6\}$ be the decomposition of $L$ defined by 
\begin{itemize}
\item 
$
\begin{array}[t]{lll}
J_1&=&\bigcup_{v\in U\cap V}\sigma_v([a_\y,b_\y]\cup[a_1,b_0,c,b_1])\ \cup \\
&&\bigcup_{v\in U'\cap V}\sigma_v([a_\y,b_\y]\cup[a_1,c,b_0,b_1])\ \cup \\
&&\bigcup_{v\in U\cap V'}\sigma_v([a_\y,b_0,b_1]\cup[b_\y,a_1])\ \cup \\
&&\bigcup_{v\in U'\cap V'}\sigma_v([a_\y,b_0,a_1]\cup[b_\y,b_1]);
\end{array}
$

\item 
$
\begin{array}[t]{lll}
J_2&=&\bigcup_{v\in U\cap V}\sigma_v([a_\y,a_0,b_\y]\cup[a_1,b_1])\ \cup \\
&&\bigcup_{v\in U'\cap V}\sigma_v([a_\y,a_0,b_\y]\cup[a_1,b_1])\ \cup \\
&&\bigcup_{v\in U\cap V'}\sigma_v([a_\y,a_1]\cup[b_\y,a_0,c,b_1])\ \cup \\
&&\bigcup_{v\in U'\cap V'}\sigma_v([a_\y,b_1]\cup[b_\y,a_0,c,a_1]);
\end{array}
$

\item 
$
\begin{array}[t]{lll}
J_3&=&\bigcup_{v\in U\cap V}\sigma_v([a_\y,b_1,b_\y]\cup[a_0,b_0])\ \cup \\
&&\bigcup_{v\in U'\cap V}\sigma_v([a_\y,c,b_1,a_0]\cup[b_\y,b_0])\ \cup \\
&&\bigcup_{v\in U\cap V'}\sigma_v([a_\y,b_1,b_\y]\cup[a_0,b_0])\ \cup \\
&&\bigcup_{v\in U'\cap V'}\sigma_v([a_\y,a_0]\cup[b_\y,c,b_1,b_0]);
\end{array}
$

\item 
$
\begin{array}[t]{lll}
J_4&=&\bigcup_{v\in U\cap V}\sigma_v([a_\y,a_1]\cup[b_\y,c,a_0,b_1])\ \cup \\
&&\bigcup_{v\in U'\cap V}\sigma_v([a_\y,b_1]\cup[b_\y,c,a_0,a_1])\ \cup \\
&&\bigcup_{v\in U\cap V'}\sigma_v([a_\y,b_\y]\cup[a_1,a_0,b_1])\ \cup \\
&&\bigcup_{v\in U'\cap V'}\sigma_v([a_\y,b_\y]\cup[a_1,a_0,b_1]);
\end{array}
$

\item 
$
\begin{array}[t]{lll}
J_5&=&\bigcup_{v\in U\cap V}\sigma_v([a_\y,c,a_1,a_0]\cup[b_\y,b_0])\ \cup \\
&&\bigcup_{v\in U'\cap V}\sigma_v([a_\y,a_1,b_\y]\cup[a_0,b_0])\ \cup \\
&&\bigcup_{v\in U\cap V'}\sigma_v([a_\y,a_0]\cup[b_\y,c,a_1,b_0])\ \cup \\
&&\bigcup_{v\in U'\cap V'}\sigma_v([a_\y,a_1,b_\y]\cup[a_0,b_0]);
\end{array}
$

\item 
$
\begin{array}[t]{lll}
J_6&=&\bigcup_{v\in U\cap V}\sigma_v([a_\y,b_0,b_1]\cup[b_\y,a_1])\ \cup \\
&&\bigcup_{v\in U'\cap V}\sigma_v([a_\y,b_0,a_1]\cup[b_\y,b_1])\ \cup \\
&&\bigcup_{v\in U\cap V'}\sigma_v([a_\y,c,b_0,b_\y]\cup[a_1,b_1])\ \cup \\
&&\bigcup_{v\in U'\cap V'}\sigma_v([a_\y,c,b_0,b_\y]\cup[a_1,b_1]);
\end{array}
$
\end{itemize}
It is easily checked that $\{J_1,J_2,J_3,J_4,J_5,J_6\}$ is a $2$-factorisation of $L$ and that $\{a_\infty^{v_1},a_\infty^{v_2},\ldots,a_\infty^{v_m}\}$ links $\{J_1,J_2,J_3,J_4,J_5,J_6\}$.
It remains to show that for $i=1,2,\ldots,m$, $\sigma_{v_i}(H)$ induces an $\{a_\infty^{v_{i}},b_\infty^{v_i}\}$-connector in $\{J_1,J_2,J_3,J_4,J_5,J_6\}$. 
There are four cases to consider: $v_i\in U\cap V$, $v_i\in U'\cap V$, $v_i\in U\cap V'$ and $v_i\in U'\cap V'$.

First suppose $v_i\in U\cap V$. In this case $a_\y^{v_i}$ and $b_\y^{v_i}$ are in the same component of $J_1$, $J_2$ and $J_3$,
and we have eight cases for the remaining $2$-factors. Namely, for each of $J_4$, $J_5$ and $J_6$, $a_\y^{v_i}$ and $b_\y^{v_i}$
are either in the same component or they are not. We number these eight cases as in the following table. 

\vspace{0.3cm}

\begin{center}
\begin{tabular}{c|c|c|c|}
&$J_4$&$J_5$&$J_6$\\
\hline
1.1& same&same&same\\
\hline
1.2& distinct&same&same\\
\hline
1.3& same&distinct&same\\
\hline
1.4& same&same&distinct\\
\hline
1.5& same&distinct&distinct\\
\hline
1.6& distinct&same&distinct\\
\hline
1.7& distinct&distinct&same\\
\hline
1.8& distinct&distinct&distinct\\
\hline
\end{tabular}
\end{center}

Depending on which of cases 1.1--1.8 that we are in, 
we can reallocate the edges of $\sigma_{v_i}(H)$ to the factors $J_1,J_2,\ldots,J_6$ 
as indicated in the following table to obtain a new 
$2$-factorisation of $L$ with the desired properties. If $J_x$ ($x\in\{1,2,\ldots,6\}$) is not listed for a particular case, then the edges of 
$\sigma_{v_i}(H)$ that are in $J_x$ are unchanged. 

\begin{center}
\begin{tabular}{|c|l|}
\hline
1.1& \\
\hline
1.2& $J_4:\sigma_{v_i}([a_\y,c,b_\y]\cup[a_1,a_0,b_1])$\\
& $J_5:\sigma_{v_i}([a_\y,a_1,c,a_0]\cup[b_\y,b_0])$\\
\hline
1.3a& $J_1:\sigma_{v_i}([a_\y,c,b_\y]\cup[a_1,b_0,b_1])$\\

& $J_2:\sigma_{v_i}([a_\y,b_\y]\cup[a_1,a_0,b_1])$\\

& $J_4:\sigma_{v_i}([a_\y,a_0,c,b_1]\cup[b_\y,a_1])$\\

& $J_5:\sigma_{v_i}([a_\y,a_1,c,b_0]\cup[b_\y,a_0])$\\

& $J_6:\sigma_{v_i}([a_\y,b_0,b_\y]\cup[a_1,b_1])$\\

&if $a_\y^{v_i}$ and $b_\y^{v_i}$ are in the same component of $J_4\setminus\sigma_{v_i}(H)$\\

&and $a_\y^{v_i}$ and $a_1^{v_i}$ are in the same component of $J_6\setminus\sigma_{v_i}(H)$.\\

1.3b& $J_4:\sigma_{v_i}([a_\y,c,a_0,b_1]\cup[b_\y,a_1])$\\

& $J_5:\sigma_{v_i}([a_\y,b_0]\cup[b_\y,c,a_1,a_0])$\\

& $J_6:\sigma_{v_i}([a_\y,a_1]\cup[b_\y,b_0,b_1])$\\

&if $a_\y^{v_i}$ and $b_\y^{v_i}$ are in the same component of $J_4\setminus\sigma_{v_i}(H)$\\

&and $a_\y^{v_i}$ and $b_\y^{v_i}$ are in the same component of $J_6\setminus\sigma_{v_i}(H)$.\\

1.3c& $J_1:\sigma_{v_i}([a_\y,c,b_\y]\cup[a_1,b_0,b_1])$\\

& $J_2:\sigma_{v_i}([a_\y,b_\y]\cup[a_1,a_0,b_1])$\\

& $J_4:\sigma_{v_i}([a_\y,a_0,b_\y]\cup[a_1,c,b_1])$\\

& $J_5:\sigma_{v_i}([a_\y,a_1,b_\y]\cup[a_0,c,b_0])$\\

& $J_6:\sigma_{v_i}([a_\y,b_0,b_\y]\cup[a_1,b_1])$\\

&if $a_\y^{v_i}$ and $b_1^{v_i}$ are in the same component of $J_4\setminus\sigma_{v_i}(H)$\\

&and $a_\y^{v_i}$ and $a_1^{v_i}$ are in the same component of $J_6\setminus\sigma_{v_i}(H)$.\\

1.3d& $J_4:\sigma_{v_i}([a_\y,c,b_\y]\cup[a_1,a_0,b_1])$\\

& $J_5:\sigma_{v_i}([a_\y,b_0]\cup[b_\y,a_1,c,a_0])$\\

& $J_6:\sigma_{v_i}([a_\y,a_1]\cup[b_\y,b_0,b_1])$\\

&if $a_\y^{v_i}$ and $b_1^{v_i}$ are in the same component of $J_4\setminus\sigma_{v_i}(H)$\\

&and $a_\y^{v_i}$ and $b_\y^{v_i}$ are in the same component of $J_6\setminus\sigma_{v_i}(H)$.\\

\hline
\end{tabular}
\end{center}

\begin{center}
\begin{tabular}{|c|l|}
\hline
1.4& $J_1:\sigma_{v_i}([a_\y,c,b_\y]\cup[a_1,b_0,b_1])$\\
& $J_2:\sigma_{v_i}([a_\y,b_\y]\cup[a_1,a_0,b_1])$\\

& $J_4:\sigma_{v_i}([a_\y,a_1]\cup[b_\y,a_0,c,b_1])$\\

& $J_5:\sigma_{v_i}([a_\y,a_0]\cup[b_\y,a_1,c,b_0])$\\

& $J_6:\sigma_{v_i}([a_\y,b_0,b_\y]\cup[a_1,b_1])$\\

\hline

1.5& $J_1:\sigma_{v_i}([a_\y,b_0,b_\y]\cup[a_1,c,b_1])$\\

& $J_5:\sigma_{v_i}([a_\y,c,b_0]\cup[b_\y,a_1,a_0])$\\

& $J_6:\sigma_{v_i}([a_\y,b_\y]\cup[a_1,b_0,b_1])$\\
\hline
1.6& $J_1:\sigma_{v_i}([a_\y,b_0,b_\y]\cup[a_1,c,b_1])$\\

& $J_4:\sigma_{v_i}([a_\y,c,a_0,b_1]\cup[b_\y,a_1])$\\

& $J_5:\sigma_{v_i}([a_\y,a_1,a_0]\cup[b_\y,c,b_0])$\\

& $J_6:\sigma_{v_i}([a_\y,b_\y]\cup[a_1,b_0,b_1])$\\
\hline
1.7a& $J_4:\sigma_{v_i}([a_\y,c,b_\y]\cup[a_1,a_0,b_1])$\\

& $J_5:\sigma_{v_i}([a_\y,b_0]\cup[b_\y,a_1,c,a_0])$\\

& $J_6:\sigma_{v_i}([a_\y,a_1]\cup[b_\y,b_0,b_1])$\\

&if $a_\y^{v_i}$ and $b_\y^{v_i}$ are in the same component of $J_6\setminus\sigma_{v_i}(H)$.\\

1.7b& $J_1:\sigma_{v_i}([a_\y,b_0,b_\y]\cup[a_1,c,b_1])$\\

& $J_4:\sigma_{v_i}([a_\y,c,b_\y]\cup[a_1,a_0,b_1])$\\

& $J_5:\sigma_{v_i}([a_\y,a_1,b_\y]\cup[a_0,c,b_0])$\\

& $J_6:\sigma_{v_i}([a_\y,b_\y]\cup[a_1,b_0,b_1])$\\

&if $a_\y^{v_i}$ and $a_1^{v_i}$ are in the same component of $J_6\setminus\sigma_{v_i}(H)$.\\
\hline
1.8& $J_4:\sigma_{v_i}([a_\y,c,b_\y]\cup[a_1,a_0,b_1])$\\

& $J_5:\sigma_{v_i}([a_\y,b_0]\cup[b_\y,a_1,c,a_0])$\\

& $J_6:\sigma_{v_i}([a_\y,a_1]\cup[b_\y,b_0,b_1])$\\
\hline
\end{tabular}
\end{center}

Now suppose $v_i\in U'\cap V$. In this case $a_\y^{v_i}$ and $b_\y^{v_i}$ are in the same component of $J_1$, $J_2$ and $J_5$,
and we have eight cases for the remaining $2$-factors. Namely, for each of $J_3$, $J_4$ and $J_6$, $a_\y^{v_i}$ and $b_\y^{v_i}$
are either in the same component or they are not. We number these eight cases as in the following table. 

\vspace{0.3cm}

\begin{center}
\begin{tabular}{c|c|c|c|}
&$J_3$&$J_4$&$J_6$\\
\hline
2.1& same&same&same\\
\hline
2.2& distinct&same&same\\
\hline
2.3& same&distinct&same\\
\hline
2.4& same&same&distinct\\
\hline
2.5& distinct&distinct&same\\
\hline
2.6& distinct&same&distinct\\
\hline
2.7& same&distinct&distinct\\
\hline
2.8& distinct&distinct&distinct\\
\hline
\end{tabular}
\end{center}

Depending on which of cases 2.1--2.8 that we are in, 
we can reallocate the edges of $\sigma_{v_i}(H)$ to the factors $J_1,J_2,\ldots,J_6$ 
as indicated in the following table to obtain a new 
$2$-factorisation of $L$ with the desired properties. If $J_x$ ($x\in\{1,2,\ldots,6\}$) is not listed for a particular case, then the edges of 
$\sigma_{v_i}(H)$ that are in $J_x$ are unchanged.

\begin{center}
\begin{tabular}{|c|l|}
\hline
2.1& \\
\hline
2.2a& $J_1:\sigma_{v_i}([a_\y,b_0,c,b_\y]\cup[a_1,b_1])$\\

& $J_2:\sigma_{v_i}([a_\y,b_\y]\cup[a_1,a_0,b_1])$\\

& $J_3:\sigma_{v_i}([a_\y,c,b_1,b_0]\cup[b_\y,a_0])$\\

& $J_4:\sigma_{v_i}([a_\y,a_0,c,a_1]\cup[b_\y,b_1])$\\

& $J_6:\sigma_{v_i}([a_\y,b_1]\cup[b_\y,b_0,a_1])$\\

&if $a_\y^{v_i}$ and $b_\y^{v_i}$ are in the same component of $J_4\setminus\sigma_{v_i}(H)$\\

&and $a_\y^{v_i}$ and $b_\y^{v_i}$  are in the same component of $J_6\setminus\sigma_{v_i}(H)$.\\

2.2b& $J_1:\sigma_{v_i}([a_\y,c,b_0,b_\y]\cup[a_1,b_1])$\\

& $J_2:\sigma_{v_i}([a_\y,b_\y]\cup[a_1,a_0,b_1])$\\

& $J_3:\sigma_{v_i}([a_\y,b_1,b_0]\cup[b_\y,c,a_0])$\\

& $J_4:\sigma_{v_i}([a_\y,a_0,b_\y]\cup[a_1,c,b_1])$\\

&if $a_\y^{v_i}$ and $a_1^{v_i}$ are in the same component of $J_4\setminus\sigma_{v_i}(H)$\\

&and $a_\y^{v_i}$ and $b_\y^{v_i}$ are in the same component of $J_6\setminus\sigma_{v_i}(H)$.\\

2.2c& $J_1:\sigma_{v_i}([a_\y,c,b_\y]\cup[a_1,b_0,b_1])$\\

& $J_2:\sigma_{v_i}([a_\y,b_\y]\cup[a_1,a_0,b_1])$\\

& $J_3:\sigma_{v_i}([a_\y,b_1,c,b_0]\cup[b_\y,a_0])$\\

& $J_4:\sigma_{v_i}([a_\y,a_0,c,a_1]\cup[b_\y,b_1])$\\

& $J_6:\sigma_{v_i}([a_\y,b_0,b_\y]\cup[a_1,b_1])$\\

&if $a_\y^{v_i}$ and $b_\y^{v_i}$ are in the same component of $J_4\setminus\sigma_{v_i}(H)$\\

&and $a_\y^{v_i}$ and $b_1^{v_i}$ are in the same component of $J_6\setminus\sigma_{v_i}(H)$.\\

2.2d& $J_1:\sigma_{v_i}([a_\y,c,b_\y]\cup[a_1,b_0,b_1])$\\

& $J_2:\sigma_{v_i}([a_\y,b_\y]\cup[a_1,a_0,b_1])$\\

& $J_3:\sigma_{v_i}([a_\y,b_1,b_\y]\cup[a_0,c,b_0])$\\

& $J_4:\sigma_{v_i}([a_\y,a_0,b_\y]\cup[a_1,c,b_1])$\\

& $J_6:\sigma_{v_i}([a_\y,b_0,b_\y]\cup[a_1,b_1])$\\

&if $a_\y^{v_i}$ and $a_1^{v_i}$ are in the same component of $J_4\setminus\sigma_{v_i}(H)$\\

&and $a_\y^{v_i}$ and $b_1^{v_i}$ are in the same component of $J_6\setminus\sigma_{v_i}(H)$.\\
\hline
\end{tabular}
\end{center}

\begin{center}
\begin{tabular}{|c|l|}

\hline
2.3& $J_3:\sigma_{v_i}([a_\y,b_1,c,a_0]\cup[b_\y,b_0])$\\

& $J_4:\sigma_{v_i}([a_\y,c,b_\y]\cup[a_1,a_0,b_1])$\\

\hline
2.4& $J_1:\sigma_{v_i}([a_\y,c,b_\y]\cup[a_1,b_0,b_1])$\\

& $J_2:\sigma_{v_i}([a_\y,b_\y]\cup[a_1,a_0,b_1])$\\

& $J_3:\sigma_{v_i}([a_\y,a_0]\cup[b_\y,b_1,c,b_0])$\\

& $J_4:\sigma_{v_i}([a_\y,b_1]\cup[b_\y,a_0,c,a_1])$\\

& $J_6:\sigma_{v_i}([a_\y,b_0,b_\y]\cup[a_1,b_1])$\\

\hline

2.5& $J_1:\sigma_{v_i}([a_\y,c,b_0,b_\y]\cup[a_1,b_1])$\\

& $J_2:\sigma_{v_i}([a_\y,b_\y]\cup[a_1,a_0,b_1])$\\

& $J_3:\sigma_{v_i}([a_\y,b_1,b_0]\cup[b_\y,c,a_0])$\\

& $J_4:\sigma_{v_i}([a_\y,a_0,b_\y]\cup[a_1,c,b_1])$\\

\hline
2.6& $J_1:\sigma_{v_i}([a_\y,b_0,b_\y]\cup[a_1,c,b_1])$\\

& $J_3:\sigma_{v_i}([a_\y,c,b_0]\cup[b_\y,b_1,a_0])$\\

& $J_6:\sigma_{v_i}([a_\y,b_\y]\cup[a_1,b_0,b_1])$\\

\hline
2.7& $J_1:\sigma_{v_i}([a_\y,b_0,c,b_\y]\cup[a_1,b_1])$\\

& $J_2:\sigma_{v_i}([a_\y,b_\y]\cup[a_1,a_0,b_1])$\\

& $J_3:\sigma_{v_i}([a_\y,c,a_0]\cup[b_\y,b_1,b_0])$\\

& $J_4:\sigma_{v_i}([a_\y,a_0,b_\y]\cup[a_1,c,b_1])$\\

& $J_6:\sigma_{v_i}([a_\y,b_1]\cup[b_\y,b_0,a_1])$\\

\hline
2.8& $J_1:\sigma_{v_i}([a_\y,b_0,c,b_\y]\cup[a_1,b_1])$\\

& $J_2:\sigma_{v_i}([a_\y,b_\y]\cup[a_1,a_0,b_1])$\\

& $J_3:\sigma_{v_i}([a_\y,c,b_1,b_0]\cup[b_\y,a_0])$\\

& $J_4:\sigma_{v_i}([a_\y,a_0,c,a_1]\cup[b_\y,b_1])$\\

& $J_6:\sigma_{v_i}([a_\y,b_1]\cup[b_\y,b_0,a_1])$\\
\hline
\end{tabular}
\end{center}

Now suppose $v_i\in U\cap V'$. In this case $a_\y^{v_i}$ and $b_\y^{v_i}$ are in the same component of $J_3$, $J_4$ and $J_6$,
and we have eight cases for the remaining $2$-factors. Namely, for each of $J_1$, $J_2$ and $J_5$, $a_\y^{v_i}$ and $b_\y^{v_i}$
are either in the same component or they are not. We number these eight cases as in the following table. 

\vspace{0.3cm}

\begin{center}
\begin{tabular}{c|c|c|c|}
&$J_1$&$J_2$&$J_5$\\
\hline
3.1& same&same&same\\
\hline
3.2& same&distinct&same\\
\hline
3.3& distinct&same&same\\
\hline
3.4& same&same&distinct\\
\hline
3.5& distinct&distinct&same\\
\hline
3.6& same&distinct&distinct\\
\hline
3.7& distinct&same&distinct\\
\hline
3.8& distinct&distinct&distinct\\
\hline
\end{tabular}
\end{center}

Depending on which of cases 3.1--3.8 that we are in, 
we can reallocate the edges of $\sigma_{v_i}(H)$ to the factors $J_1,J_2,\ldots,J_6$ 
as indicated in the following table to obtain a new 
$2$-factorisation of $L$ with the desired properties. If $J_x$ ($x\in\{1,2,\ldots,6\}$) is not listed for a particular case, then the edges of 
$\sigma_{v_i}(H)$ that are in $J_x$ are unchanged.

\begin{center}
\begin{tabular}{|c|l|}
\hline
3.1& \\
\hline
3.2& $J_2:\sigma_{v_i}([a_\y,c,b_\y]\cup[a_1,a_0,b_1])$\\

& $J_4:\sigma_{v_i}([a_\y,a_0,b_\y]\cup[a_1,b_1])$\\

& $J_5:\sigma_{v_i}([a_\y,a_1,c,a_0]\cup[b_\y,b_0])$\\

& $J_6:\sigma_{v_i}([a_\y,b_\y]\cup[a_1,b_0,c,b_1])$\\

\hline
3.3& $J_1:\sigma_{v_i}([a_\y,b_0,b_\y]\cup[a_1,b_1])$\\

& $J_5:\sigma_{v_i}([a_\y,a_0]\cup[b_\y,a_1,c,b_0])$\\

& $J_6:\sigma_{v_i}([a_\y,c,b_\y]\cup[a_1,b_0,b_1])$\\

\hline
3.4a& $J_1:\sigma_{v_i}([a_\y,a_1]\cup[b_\y,b_0,b_1])$\\

& $J_2:\sigma_{v_i}([a_\y,a_0,c,b_1]\cup[b_\y,a_1])$\\

& $J_5:\sigma_{v_i}([a_\y,c,a_1,b_0]\cup[b_\y,a_0])$\\

& $J_6:\sigma_{v_i}([a_\y,b_0,c,b_\y]\cup[a_1,b_1])$\\

&if $a_\y^{v_i}$ and $b_\y^{v_i}$ are in the same component of $J_1\setminus\sigma_{v_i}(H)$\\

&and $a_\y^{v_i}$ and $b_\y^{v_i}$  are in the same component of $J_2\setminus\sigma_{v_i}(H)$.\\

3.4b& $J_2:\sigma_{v_i}([a_\y,a_0,b_\y]\cup[a_1,c,b_1])$\\

& $J_5:\sigma_{v_i}([a_\y,a_1,b_0]\cup[b_\y,c,a_0])$\\

&if $a_\y^{v_i}$ and $b_\y^{v_i}$ are in the same component of $J_1\setminus\sigma_{v_i}(H)$\\

&and $a_\y^{v_i}$ and $b_1^{v_i}$  are in the same component of $J_2\setminus\sigma_{v_i}(H)$.\\

3.4c& $J_1:\sigma_{v_i}([a_\y,b_\y]\cup[a_1,b_0,b_1])$\\

& $J_2:\sigma_{v_i}([a_\y,a_1]\cup[b_\y,c,a_0,b_1])$\\

& $J_4:\sigma_{v_i}([a_\y,a_0,b_\y]\cup[a_1,b_1])$\\

& $J_5:\sigma_{v_i}([a_\y,c,b_0]\cup[b_\y,a_1,a_0])$\\

& $J_6:\sigma_{v_i}([a_\y,b_0,b_\y]\cup[a_1,c,b_1])$\\

&if $a_\y^{v_i}$ and $a_1^{v_i}$ are in the same component of $J_1\setminus\sigma_{v_i}(H)$\\

&and $a_\y^{v_i}$ and $b_\y^{v_i}$ are in the same component of $J_2\setminus\sigma_{v_i}(H)$.\\

3.4d& $J_1:\sigma_{v_i}([a_\y,b_\y]\cup[a_1,b_0,b_1])$\\

& $J_2:\sigma_{v_i}([a_\y,c,b_\y]\cup[a_1,a_0,b_1])$\\

& $J_4:\sigma_{v_i}([a_\y,a_0,b_\y]\cup[a_1,b_1])$\\

& $J_5:\sigma_{v_i}([a_\y,a_1,b_\y]\cup[a_0,c,b_0])$\\

& $J_6:\sigma_{v_i}([a_\y,b_0,b_\y]\cup[a_1,c,b_1])$\\

&if $a_\y^{v_i}$ and $a_1^{v_i}$ are in the same component of $J_1\setminus\sigma_{v_i}(H)$\\

&and $a_\y^{v_i}$ and $b_1^{v_i}$ are in the same component of $J_2\setminus\sigma_{v_i}(H)$.\\
\hline
\end{tabular}
\end{center}

\begin{center}
\begin{tabular}{|c|l|}

\hline

3.5& $J_1:\sigma_{v_i}([a_\y,b_\y]\cup[a_1,b_0,b_1])$\\

& $J_2:\sigma_{v_i}([a_\y,c,a_0,b_1]\cup[b_\y,a_1])$\\

& $J_4:\sigma_{v_i}([a_\y,a_0,b_\y]\cup[a_1,b_1])$\\

& $J_5:\sigma_{v_i}([a_\y,a_1,a_0]\cup[b_\y,c,b_0])$\\

& $J_6:\sigma_{v_i}([a_\y,b_0,b_\y]\cup[a_1,c,b_1])$\\

\hline
3.6& $J_2:\sigma_{v_i}([a_\y,a_0,b_\y]\cup[a_1,c,b_1])$\\

& $J_5:\sigma_{v_i}([a_\y,a_1,b_0]\cup[b_\y,c,a_0])$\\

\hline
3.7& $J_1:\sigma_{v_i}([a_\y,b_\y]\cup[a_1,b_0,b_1])$\\

& $J_2:\sigma_{v_i}([a_\y,a_1]\cup[b_\y,c,a_0,b_1])$\\

& $J_4:\sigma_{v_i}([a_\y,a_0,b_\y]\cup[a_1,b_1])$\\

& $J_5:\sigma_{v_i}([a_\y,c,b_0]\cup[b_\y,a_1,a_0])$\\

& $J_6:\sigma_{v_i}([a_\y,b_0,b_\y]\cup[a_1,c,b_1])$\\

\hline
3.8& $J_1:\sigma_{v_i}([a_\y,b_\y]\cup[a_1,b_0,b_1])$\\

& $J_2:\sigma_{v_i}([a_\y,c,b_\y]\cup[a_1,a_0,b_1])$\\

& $J_4:\sigma_{v_i}([a_\y,a_0,b_\y]\cup[a_1,b_1])$\\

& $J_5:\sigma_{v_i}([a_\y,a_1,b_\y]\cup[a_0,c,b_0])$\\

& $J_6:\sigma_{v_i}([a_\y,b_0,b_\y]\cup[a_1,c,b_1])$\\
\hline
\end{tabular}
\end{center}

Finally suppose $v_i\in U'\cap V'$. In this case $a_\y^{v_i}$ and $b_\y^{v_i}$ are in the same component of $J_4$, $J_5$ and $J_6$,
and we have eight cases for the remaining $2$-factors. Namely, for each of $J_1$, $J_2$ and $J_3$, $a_\y^{v_i}$ and $b_\y^{v_i}$
are either in the same component or they are not. We number these eight cases as in the following table. 

\vspace{0.3cm}

\begin{center}
\begin{tabular}{c|c|c|c|}
&$J_1$&$J_2$&$J_3$\\
\hline
4.1& same&same&same\\
\hline
4.2& distinct&same&same\\
\hline
4.3& same&same&distinct\\
\hline
4.4& same&distinct&same\\
\hline
4.5& distinct&same&distinct\\
\hline
4.6& distinct&distinct&same\\
\hline
4.7& same&distinct&distinct\\
\hline
4.8& distinct&distinct&distinct\\
\hline
\end{tabular}
\end{center}

Depending on which of cases 4.1--4.8 that we are in, 
we can reallocate the edges of $\sigma_{v_i}(H)$ to the factors $J_1,J_2,\ldots,J_6$ 
as indicated in the following table to obtain a new 
$2$-factorisation of $L$ with the desired properties. If $J_x$ ($x\in\{1,2,\ldots,6\}$) is not listed for a particular case, then the edges of 
$\sigma_{v_i}(H)$ that are in $J_x$ are unchanged.

\begin{center}
\begin{tabular}{|c|l|}
\hline
4.1& \\
\hline
4.2& $J_1:\sigma_{v_i}([a_\y,b_0,b_\y]\cup[a_1,b_1])$\\

& $J_3:\sigma_{v_i}([a_\y,a_0]\cup[b_\y,b_1c,b_0])$\\

& $J_6:\sigma_{v_i}([a_\y,c,b_\y]\cup[a_1,b_0,b_1])$\\

\hline
4.3a& $J_1:\sigma_{v_i}([a_\y,b_1]\cup[b_\y,b_0,a_1])$\\

& $J_2:\sigma_{v_i}([a_\y,c,a_0,a_1]\cup[b_\y,b_1])$\\

& $J_3:\sigma_{v_i}([a_\y,b_0]\cup[b_\y,c,b_1,a_0])$\\

& $J_4:\sigma_{v_i}([a_\y,a_0,b_\y]\cup[a_1,b_1])$\\

& $J_6:\sigma_{v_i}([a_\y,b_\y]\cup[a_1,c,b_0,b_1])$\\

&if $a_\y^{v_i}$ and $b_\y^{v_i}$ are in the same component of $J_1\setminus\sigma_{v_i}(H)$\\

&and $a_\y^{v_i}$ and $b_\y^{v_i}$  are in the same component of $J_2\setminus\sigma_{v_i}(H)$.\\

4.3b& $J_1:\sigma_{v_i}([a_\y,b_\y]\cup[a_1,b_0,b_1])$\\

& $J_2:\sigma_{v_i}([a_\y,b_1]\cup[b_\y,c,a_0,a_1])$\\

& $J_3:\sigma_{v_i}([a_\y,c,b_0]\cup[b_\y,b_1,a_0])$\\

& $J_4:\sigma_{v_i}([a_\y,a_0,b_\y]\cup[a_1,b_1])$\\

& $J_6:\sigma_{v_i}([a_\y,b_0,b_\y]\cup[a_1,c,b_1])$\\

&if $a_\y^{v_i}$ and $b_1^{v_i}$ are in the same component of $J_1\setminus\sigma_{v_i}(H)$\\

&and $a_\y^{v_i}$ and $b_\y^{v_i}$  are in the same component of $J_2\setminus\sigma_{v_i}(H)$.\\

4.3c& $J_1:\sigma_{v_i}([a_\y,b_1]\cup[b_\y,b_0,a_1])$\\

& $J_2:\sigma_{v_i}([a_\y,c,b_\y]\cup[a_1,a_0,b_1])$\\

& $J_3:\sigma_{v_i}([a_\y,b_0]\cup[b_\y,b_1,c,a_0])$\\

& $J_4:\sigma_{v_i}([a_\y,a_0,b_\y]\cup[a_1,b_1])$\\

& $J_6:\sigma_{v_i}([a_\y,b_\y]\cup[a_1,c,b_0,b_1])$\\

&if $a_\y^{v_i}$ and $b_\y^{v_i}$ are in the same component of $J_1\setminus\sigma_{v_i}(H)$\\

&and $a_\y^{v_i}$ and $a_1^{v_i}$ are in the same component of $J_2\setminus\sigma_{v_i}(H)$.\\

4.3d& $J_1:\sigma_{v_i}([a_\y,b_\y]\cup[a_1,b_0,b_1])$\\

& $J_2:\sigma_{v_i}([a_\y,c,b_\y]\cup[a_1,a_0,b_1])$\\

& $J_3:\sigma_{v_i}([a_\y,b_1,b_\y]\cup[a_0,c,b_0])$\\

& $J_4:\sigma_{v_i}([a_\y,a_0,b_\y]\cup[a_1,b_1])$\\

& $J_6:\sigma_{v_i}([a_\y,b_0,b_\y]\cup[a_1,c,b_1])$\\

&if $a_\y^{v_i}$ and $b_1^{v_i}$ are in the same component of $J_1\setminus\sigma_{v_i}(H)$\\

&and $a_\y^{v_i}$ and $a_1^{v_i}$ are in the same component of $J_2\setminus\sigma_{v_i}(H)$.\\
\hline
\end{tabular}
\end{center}

\begin{center}
\begin{tabular}{|c|l|}

\hline

4.4& $J_2:\sigma_{v_i}([a_\y,c,b_\y]\cup[a_1,a_0,b_1])$\\

& $J_3:\sigma_{v_i}([a_\y,b_1,c,a_0]\cup[b_\y,b_0])$\\

& $J_4:\sigma_{v_i}([a_\y,a_0,b_\y]\cup[a_1,b_1])$\\

& $J_6:\sigma_{v_i}([a_\y,b_\y]\cup[a_1,c,b_0,b_1])$\\

\hline
4.5& $J_1:\sigma_{v_i}([a_\y,b_\y]\cup[a_1,b_0,b_1])$\\

& $J_2:\sigma_{v_i}([a_\y,b_1]\cup[b_\y,c,a_0,a_1])$\\

& $J_3:\sigma_{v_i}([a_\y,c,b_0]\cup[b_\y,b_1,a_0])$\\

& $J_4:\sigma_{v_i}([a_\y,a_0,b_\y]\cup[a_1,b_1])$\\

& $J_6:\sigma_{v_i}([a_\y,b_0,b_\y]\cup[a_1,c,b_1])$\\

\hline
4.6& $J_1:\sigma_{v_i}([a_\y,b_\y]\cup[a_1,b_0,b_1])$\\

& $J_2:\sigma_{v_i}([a_\y,c,a_0,a_1]\cup[b_\y,b_1])$\\

& $J_3:\sigma_{v_i}([a_\y,b_1,a_0]\cup[b_\y,c,b_0])$\\

& $J_4:\sigma_{v_i}([a_\y,a_0,b_\y]\cup[a_1,b_1])$\\

& $J_6:\sigma_{v_i}([a_\y,b_0,b_\y]\cup[a_1,c,b_1])$\\

\hline
4.7& $J_2:\sigma_{v_i}([a_\y,a_0,b_\y]\cup[a_1,c,b_1])$\\

& $J_3:\sigma_{v_i}([a_\y,b_1,b_0]\cup[b_\y,c,a_0])$\\

\hline
4.8& $J_1:\sigma_{v_i}([a_\y,b_1]\cup[b_\y,b_0,a_1])$\\

& $J_2:\sigma_{v_i}([a_\y,c,a_0,a_1]\cup[b_\y,b_1])$\\

& $J_3:\sigma_{v_i}([a_\y,b_0]\cup[b_\y,c,b_1,a_0])$\\

& $J_4:\sigma_{v_i}([a_\y,a_0,b_\y]\cup[a_1,b_1])$\\

& $J_6:\sigma_{v_i}([a_\y,b_\y]\cup[a_1,c,b_0,b_1])$\\
\hline
\end{tabular}
\end{center}

\vspace{5cm}

\noindent{\bf The case n=4:}
Let $H$ be the union of the paths 
$$[a_{\infty},a_1,b_0,b_{\infty}], [a_2,b_2], [b_{\infty},a_1,c,b_2], [a_2,b_0,a_{\infty}],$$
and let $H'=K_{A_4\cup B_4\cup\{c\}}-H$ so that $\{H,H'\}$ is a decomposition of $K_{A_4\cup B_4\cup\{c\}}\cong K_9$.
We show that each of $H$ and $H'$ is a Hamilton fragment. 

For $H$, our aim is to show that the subgraph $$L=\bigcup_{v\in V(G)}\sigma_v(H)$$
of $L(G)$ decomposes into Hamilton cycles, and we do this by applying Lemma \ref{manyrepairs} 
with $H_i=\sigma_{v_i}(H)$ and $V_i=\{a_\infty^{v_i},b_\infty^{v_i}\}$ for each $i=1,2,\ldots,m$.
To this end, observe that \linebreak $\sigma_{v_1}(H),\sigma_{v_2}(H),\ldots,\sigma_{v_m}(H)$ are edge-disjoint subgraphs of $L$,
and since $b_\infty^{v_i}=a_\infty^{v_{i+1}}$ for $i=1,2,\ldots,m-1$ it follows that $V_i\cap V_{i+1}\ne \emptyset$ for $i=1,2,\ldots,m-1$ as required. It remains to show there is a $2$-factorisation
$\{J_1,J_2\}$ of $L$ such that 
\begin{itemize}
\item $\{a_\infty^{v_1},a_\infty^{v_2},\ldots,a_\infty^{v_m}\}$ links $\{J_1,J_2\}$; and
\item $\sigma_{v_i}(H)$ induces an $\{a_\infty^{v_{i}},b_\infty^{v_i}\}$-connector in $\{J_1,J_2\}$, for $i=1,2,\ldots,m$.
\end{itemize} 
Let $\{U,U'\}$ be a partition of $V(G)$ such that both $U$ and $U'$ link $\{F,\F_0\}$
(such a partition exists by Lemma \ref{intersects}),
and let $\{J_1,J_2\}$ be the decomposition of $L$ defined by 
$$J_1=\bigcup_{v\in U}\sigma_v([a_{\infty},a_1,b_0,b_{\infty}] \cup [a_2,b_2])\ \cup \ 
\bigcup_{v\in U'}\sigma_v([b_{\infty},a_1,c,b_2] \cup [a_2,b_0,a_{\infty}])$$
and
$$J_2=\bigcup_{v\in U}\sigma_v([b_{\infty},a_1,c,b_2] \cup [a_2,b_0,a_{\infty}])\ \cup \ 
\bigcup_{v\in U'}\sigma_v(a_{\infty},a_1,b_0,b_{\infty}] \cup [a_2,b_2]).$$

It is easily checked that $\{J_1,J_2\}$ is a $2$-factorisation of $L$ and that 
$\{a_\infty^{v_1},a_\infty^{v_2},\ldots,a_\infty^{v_m}\}$ links $\{J_1,J_2\}$. By Lemma \ref{connectors} (0), $\sigma_{v_i}(H)$ induces an $\{a_\y^{v_i}, b_\y^{v_i}\}$-connector in $\{J_1, J_2\}$ with $(\alpha, \beta, u,v,w) = (a_\y, b_\y, a_1, c, b_0)$ for $v_i\in U$ and $(\alpha, \beta, u,v,w) = (b_\y, a_\y, a_1, c, b_0)$ for $v_i\in U'$. 

For $H'$, our aim is to show that the subgraph $$L=\bigcup_{v\in V(G)}\sigma_v(H')$$
of $L(G)$ decomposes into Hamilton cycles, and we do this by applying Lemma \ref{manyrepairs} 
with $H_i=\sigma_{v_i}(H')$ and $V_i=\{a_\infty^{v_i},b_\infty^{v_i}\}$ for each $i=1,2,\ldots,m$.
To this end, observe that \linebreak $\sigma_{v_1}(H'),\sigma_{v_2}(H'),\ldots,\sigma_{v_m}(H')$ are edge-disjoint subgraphs of $L$,
and since $b_\infty^{v_i}=a_\infty^{v_{i+1}}$ for $i=1,2,\ldots,m-1$ it follows that $V_i\cap V_{i+1}\ne \emptyset$ for $i=1,2,\ldots,m-1$ as required. It remains to show there is a $2$-factorisation
$\{J_1,J_2,J_3,J_4,J_5,J_6\}$ of $L$ such that 
\begin{itemize}
\item $\{a_\infty^{v_1},a_\infty^{v_2},\ldots,a_\infty^{v_m}\}$ links $\{J_1,J_2,J_3,J_4,J_5,J_6\}$; and
\item $\sigma_{v_i}(H')$ induces an $\{a_\infty^{v_{i}},b_\infty^{v_i}\}$-connector in $\{J_1,J_2,J_3,J_4,J_5,J_6\}$, for $i=1,2,\ldots,m$.
\end{itemize} 
Let $\{U,U'\}$ be a partition of $V(G)$ such that both $U$ and $U'$ link $\{F,\F_1\}$,
let $\{V,V'\}$ be a partition of $V(G)$ such that both $V$ and $V'$ link $\{F,\F_2\}$ 
(such partitions exist by Lemma \ref{intersects}),
and let $\{J_1,J_2,J_3,J_4,J_5,J_6\}$ be the decomposition of $L$ defined by 

\begin{itemize}
\item 
$
\begin{array}[t]{lll}
J_1&=&\bigcup_{v\in U\cap V}\sigma_v([a_1,a_0,a_2,b_1] \cup [a_{\infty},b_{\infty}])\ \cup \\
&&\bigcup_{v\in U'\cap V}\sigma_v([a_{\infty},c,a_2,b_1] \cup [a_1,a_0,b_{\infty}])\ \cup \\
&&\bigcup_{v\in U\cap V'}\sigma_v([a_1,a_0,a_2,b_1] \cup [a_{\infty},b_{\infty}])\ \cup \\
&&\bigcup_{v\in U'\cap V'}\sigma_v([a_{\infty},c,a_2,b_1] \cup [a_1,a_0,b_{\infty}]);
\end{array}
$
\item 
$
\begin{array}[t]{lll}
J_2&=&\bigcup_{v\in U\cap V}\sigma_v([a_1,b_2,b_0,c,a_{\infty}] \cup [b_{\infty},b_1])\ \cup \\
&&\bigcup_{v\in U'\cap V}\sigma_v([a_1,b_2,b_0,b_1] \cup [a_{\infty},b_{\infty}])\ \cup \\
&&\bigcup_{v\in U\cap V'}\sigma_v([a_1,b_2,b_0,c,a_{\infty}] \cup [b_{\infty},b_1])\ \cup \\
&&\bigcup_{v\in U'\cap V'}\sigma_v([a_1,b_2,b_0,b_1] \cup [a_{\infty},b_{\infty}]);
\end{array}
$
\item 
$
\begin{array}[t]{lll}
J_3&=&\bigcup_{v\in U\cap V}\sigma_v([a_{\infty},a_2,c,a_0,b_{\infty}] \cup [a_1,b_1])\ \cup \\
&&\bigcup_{v\in U'\cap V}\sigma_v([a_1,a_2,b_{\infty}] \cup [a_{\infty},a_0,b_1])\ \cup \\
&&\bigcup_{v\in U\cap V'}\sigma_v([a_{\infty},a_2,c,a_0,b_{\infty}] \cup [a_1,b_1])\ \cup \\
&&\bigcup_{v\in U'\cap V'}\sigma_v([a_1,a_2,b_{\infty}] \cup [a_{\infty},a_0,b_1]);
\end{array}
$
\item 
$
\begin{array}[t]{lll}
J_4&=&\bigcup_{v\in U\cap V}\sigma_v([a_1,a_2,b_{\infty}] \cup [a_{\infty},a_0,b_1])\ \cup \\
&&\bigcup_{v\in U'\cap V}\sigma_v([a_{\infty},a_2,a_0,c,b_{\infty}] \cup [a_1,b_1])\ \cup \\
&&\bigcup_{v\in U\cap V'}\sigma_v([a_1,a_2,b_{\infty}] \cup [a_{\infty},a_0,b_1])\ \cup \\
&&\bigcup_{v\in U'\cap V'}\sigma_v([a_{\infty},a_2,a_0,c,b_{\infty}] \cup [a_1,b_1]);
\end{array}
$
\item 
$
\begin{array}[t]{lll}
J_5&=&\bigcup_{v\in U\cap V}\sigma_v([a_{\infty},b_1,b_2,b_{\infty}] \cup [a_0,b_0])\ \cup \\
&&\bigcup_{v\in U'\cap V}\sigma_v([a_{\infty},b_1,b_2,b_{\infty}] \cup [a_0,b_0])\ \cup \\
&&\bigcup_{v\in U\cap V'}\sigma_v([b_{\infty},c,b_1,b_0] \cup [a_0,b_2,a_{\infty}])\ \cup \\
&&\bigcup_{v\in U'\cap V'}\sigma_v([b_{\infty},b_1,c,b_0] \cup [a_0,b_2,a_{\infty}]);
\end{array}
$
\item 
$
\begin{array}[t]{lll}
J_6&=&\bigcup_{v\in U\cap V}\sigma_v([b_{\infty},c,b_1,b_0] \cup [a_0,b_2,a_{\infty}])\ \cup \\
&&\bigcup_{v\in U'\cap V}\sigma_v([b_{\infty},b_1,c,b_0] \cup [a_0,b_2,a_{\infty}])\ \cup \\
&&\bigcup_{v\in U\cap V'}\sigma_v([a_{\infty},b_1,b_2,b_{\infty}] \cup [a_0,b_0])\ \cup \\
&&\bigcup_{v\in U'\cap V'}\sigma_v([a_{\infty},b_1,b_2,b_{\infty}] \cup [a_0,b_0]);
\end{array}
$

\end{itemize}
It is easily checked that $\{J_1,J_2,J_3,J_4,J_5,J_6\}$ is a $2$-factorisation of $L$ and that $\{a_\infty^{v_1},a_\infty^{v_2},\ldots,a_\infty^{v_m}\}$ links $\{J_1,J_2,J_3,J_4,J_5,J_6\}$.
It remains to show that for $i=1,2,\ldots,m$, $\sigma_{v_i}(H')$ induces an $\{a_\infty^{v_{i}},b_\infty^{v_i}\}$-connector in $\{J_1,J_2,J_3,J_4,J_5,J_6\}$. 
There are four cases to consider: $v_i\in U\cap V$, $v_i\in U'\cap V$, $v_i\in U\cap V'$ and $v_i\in U'\cap V'$.

First suppose $v_i\in U\cap V$. In this case $a_\y^{v_i}$ and $b_\y^{v_i}$ are in the same component of $J_1$, $J_3$ and $J_5$,
and we have eight cases for the remaining $2$-factors. Namely, for each of $J_2$, $J_4$ and $J_6$, $a_\y^{v_i}$ and $b_\y^{v_i}$
are either in the same component or they are not. We number these eight cases as in the following table. 

\vspace{0.3cm}

\begin{center}
\begin{tabular}{c|c|c|c|}
&$J_2$&$J_4$&$J_6$\\
\hline
1.1& same&same&same\\
\hline
1.2& distinct&same&same\\
\hline
1.3& same&distinct&same\\
\hline
1.4& same&same&distinct\\
\hline
1.5& same&distinct&distinct\\
\hline
1.6& distinct&same&distinct\\
\hline
1.7& distinct&distinct&same\\
\hline
1.8& distinct&distinct&distinct\\
\hline
\end{tabular}
\end{center}

Depending on which of cases 1.1--1.8 that we are in, 
we can reallocate the edges of $\sigma_{v_i}(H)$ to the factors $J_1,J_2,\ldots,J_6$ 
as indicated in the following table to obtain a new 
$2$-factorisation of $L$ with the desired properties. If $J_x$ ($x\in\{1,2,\ldots,6\}$) is not listed for a particular case, then the edges of 
$\sigma_{v_i}(H)$ that are in $J_x$ are unchanged. 

\begin{center}
\begin{tabular}{|c|l|}
\hline
1.1& \\
\hline
1.2& $J_2:\sigma_{v_i}([a_1,b_2,b_0,b_1] \cup [a_{\infty},c,b_{\infty}])$\\
& $J_6:\sigma_{v_i}([b_{\infty},b_1,c,b_0] \cup [a_0,b_2,a_{\infty}])$\\
\hline
1.3& $J_3:\sigma_{v_i}([a_{\infty},a_0,c,a_2,b_{\infty}] \cup [a_1,b_1])$\\
& $J_4:\sigma_{v_i}([a_1,a_2,a_{\infty}] \cup [b_{\infty},a_0,b_1])$\\
\hline
1.4& $J_2:\sigma_{v_i}([b_{\infty},c,b_0,b_1] \cup [a_1,b_2,a_{\infty}])$\\
& $J_5:\sigma_{v_i}([a_0,b_2,b_0] \cup [a_{\infty},b_1,b_{\infty}])$\\
& $J_6:\sigma_{v_i}([a_{\infty},c,b_1,b_2,b_{\infty}] \cup [a_0,b_0])$\\
\hline
1.5& $J_2:\sigma_{v_i}([b_{\infty},c,b_0,b_1] \cup [a_1,b_2,a_{\infty}])$\\
& $J_3:\sigma_{v_i}([a_{\infty},a_0,c,a_2,b_{\infty}] \cup [a_1,b_1])$\\
& $J_4:\sigma_{v_i}([a_1,a_2,a_{\infty}] \cup [b_{\infty},a_0,b_1])$\\
& $J_5:\sigma_{v_i}([a_0,b_2,b_0] \cup [a_{\infty},b_1,b_{\infty}])$\\
& $J_6:\sigma_{v_i}([a_{\infty},c,b_1,b_2,b_{\infty}] \cup [a_0,b_0])$\\
\hline
1.6& $J_2:\sigma_{v_i}([a_1,b_2,b_0,b_1] \cup [a_{\infty},c,b_{\infty}])$\\
& $J_5:\sigma_{v_i}([a_{\infty},b_2,b_1,b_{\infty}] \cup [a_0,b_0])$\\
& $J_6:\sigma_{v_i}([a_{\infty},b_1,c,b_0] \cup [a_0,b_2,b_{\infty}])$\\
\hline
1.7& $J_2:\sigma_{v_i}([a_1,b_2,b_0,b_1] \cup [a_{\infty},c,b_{\infty}])$\\
& $J_3:\sigma_{v_i}([a_{\infty},a_0,c,a_2,b_{\infty}] \cup [a_1,b_1])$\\
& $J_4:\sigma_{v_i}([a_1,a_2,a_{\infty}] \cup [b_{\infty},a_0,b_1])$\\
& $J_6:\sigma_{v_i}([b_{\infty},b_1,c,b_0] \cup [a_0,b_2,a_{\infty}])$\\
\hline
1.8& $J_2:\sigma_{v_i}([a_1,b_2,b_0,b_1] \cup [a_{\infty},c,b_{\infty}])$\\
& $J_3:\sigma_{v_i}([a_{\infty},a_0,c,a_2,b_{\infty}] \cup [a_1,b_1])$\\
& $J_4:\sigma_{v_i}([a_1,a_2,a_{\infty}] \cup [b_{\infty},a_0,b_1])$\\
& $J_5:\sigma_{v_i}([a_{\infty},b_2,b_1,b_{\infty}] \cup [a_0,b_0])$\\
& $J_6:\sigma_{v_i}([a_{\infty},b_1,c,b_0] \cup [a_0,b_2,b_{\infty}])$\\
\hline
\end{tabular}
\end{center}

Now suppose $v_i\in U'\cap V$. In this case $a_\y^{v_i}$ and $b_\y^{v_i}$ are in the same component of $J_2$, $J_4$ and $J_5$,
and we have eight cases for the remaining $2$-factors. Namely, for each of $J_1$, $J_3$ and $J_6$, $a_\y^{v_i}$ and $b_\y^{v_i}$
are either in the same component or they are not. We number these eight cases as in the following table. 

\vspace{0.3cm}

\begin{center}
\begin{tabular}{c|c|c|c|}
&$J_1$&$J_3$&$J_6$\\
\hline
2.1& same&same&same\\
\hline
2.2& distinct&same&same\\
\hline
2.3& same&distinct&same\\
\hline
2.4& same&same&distinct\\
\hline
2.5& distinct&distinct&same\\
\hline
2.6& distinct&same&distinct\\
\hline
2.7& same&distinct&distinct\\
\hline
2.8& distinct&distinct&distinct\\
\hline
\end{tabular}
\end{center}

Depending on which of cases 2.1--2.8 that we are in, 
we can reallocate the edges of $\sigma_{v_i}(H)$ to the factors $J_1,J_2,\ldots,J_6$ 
as indicated in the following table to obtain a new 
$2$-factorisation of $L$ with the desired properties. If $J_x$ ($x\in\{1,2,\ldots,6\}$) is not listed for a particular case, then the edges of 
$\sigma_{v_i}(H)$ that are in $J_x$ are unchanged.

\begin{center}
\begin{tabular}{|c|l|}
\hline
2.1& \\
\hline
2.2& $J_1:\sigma_{v_i}([a_1,a_0,a_2,b_1] \cup [a_{\infty},c,b_{\infty}])$\\
& $J_4:\sigma_{v_i}([a_{\infty},a_2,c,a_0,b_{\infty}] \cup [a_1,b_1])$\\
\hline
2.3a& $J_1:\sigma_{v_i}([a_1,a_0,c,a_{\infty}] \cup [b_{\infty},a_2,b_1])$\\

& $J_3:\sigma_{v_i}([a_1,a_2,a_{\infty}] \cup [b_{\infty},a_0,b_1])$\\

& $J_4:\sigma_{v_i}([a_{\infty},a_0,a_2,c,b_{\infty}] \cup [a_1,b_1])$\\

&if $a_\y^{v_i}$ and $b_\y^{v_i}$ are in the same component of $J_1\setminus\sigma_{v_i}(H)$.\\
2.3b& $J_1:\sigma_{v_i}([a_{\infty},c,a_2,b_{\infty}] \cup [a_1,a_0,b_1])$\\

& $J_3:\sigma_{v_i}([a_{\infty},a_2,a_0,b_{\infty}] \cup [a_1,b_1])$\\

& $J_4:\sigma_{v_i}([a_{\infty},a_0,c,b_{\infty}] \cup [a_1,a_2,b_1])$\\

&if $a_\y^{v_i}$ and $a_1^{v_i}$ are in the same component of $J_1\setminus\sigma_{v_i}(H)$.\\
\hline
2.4 & $J_5:\sigma_{v_i}([a_{\infty},b_2,b_1,b_{\infty}] \cup [a_0,b_0])$\\
& $J_6:\sigma_{v_i}([a_{\infty},b_1,c,b_0] \cup [a_0,b_2,b_{\infty}])$\\
\hline
2.5 & $J_1:\sigma_{v_i}([a_{\infty},c,a_2,b_{\infty}] \cup [a_1,a_0,b_1])$\\
& $J_3:\sigma_{v_i}([a_{\infty},a_2,a_0,b_{\infty}] \cup [a_1,b_1])$\\
& $J_4:\sigma_{v_i}([a_{\infty},a_0,c,b_{\infty}] \cup [a_1,a_2,b_1])$\\
\hline
2.6 & $J_1:\sigma_{v_i}([a_1,a_0,a_2,b_1] \cup [a_{\infty},c,b_{\infty}])$\\
& $J_4:\sigma_{v_i}([a_{\infty},a_2,c,a_0,b_{\infty}] \cup [a_1,b_1])$\\
& $J_5:\sigma_{v_i}([a_{\infty},b_2,b_1,b_{\infty}] \cup [a_0,b_0])$\\
& $J_6:\sigma_{v_i}([a_{\infty},b_1,c,b_0] \cup [a_0,b_2,b_{\infty}])$\\
\hline
2.7a& $J_1:\sigma_{v_i}([a_1,a_0,c,a_{\infty}] \cup [b_{\infty},a_2,b_1])$\\

& $J_3:\sigma_{v_i}([a_1,a_2,a_{\infty}] \cup [b_{\infty},a_0,b_1])$\\

& $J_4:\sigma_{v_i}([a_{\infty},a_0,a_2,c,b_{\infty}] \cup [a_1,b_1])$\\

& $J_5:\sigma_{v_i}([a_{\infty},b_2,b_1,b_{\infty}] \cup [a_0,b_0])$\\

& $J_6:\sigma_{v_i}([a_{\infty},b_1,c,b_0] \cup [a_0,b_2,b_{\infty}])$\\

&if $a_\y^{v_i}$ and $b_\y^{v_i}$ are in the same component of $J_1\setminus\sigma_{v_i}(H)$.\\
2.7b& $J_1:\sigma_{v_i}([a_{\infty},c,a_2,b_{\infty}] \cup [a_1,a_0,b_1])$\\

& $J_3:\sigma_{v_i}([a_{\infty},a_2,a_0,b_{\infty}] \cup [a_1,b_1])$\\

& $J_4:\sigma_{v_i}([a_{\infty},a_0,c,b_{\infty}] \cup [a_1,a_2,b_1])$\\

& $J_5:\sigma_{v_i}([a_{\infty},b_2,b_1,b_{\infty}] \cup [a_0,b_0])$\\

& $J_6:\sigma_{v_i}([a_{\infty},b_1,c,b_0] \cup [a_0,b_2,b_{\infty}])$\\

&if $a_\y^{v_i}$ and $a_1^{v_i}$ are in the same component of $J_1\setminus\sigma_{v_i}(H)$.\\
\hline
2.8& $J_1:\sigma_{v_i}([a_1,a_0,c,a_{\infty}] \cup [b_{\infty},a_2,b_1])$\\

& $J_3:\sigma_{v_i}([a_1,a_2,a_{\infty}] \cup [b_{\infty},a_0,b_1])$\\

& $J_4:\sigma_{v_i}([a_{\infty},a_0,a_2,c,b_{\infty}] \cup [a_1,b_1])$\\

& $J_5:\sigma_{v_i}([a_{\infty},b_2,b_1,b_{\infty}] \cup [a_0,b_0])$\\

& $J_6:\sigma_{v_i}([a_{\infty},b_1,c,b_0] \cup [a_0,b_2,b_{\infty}])$\\
\hline
\end{tabular}
\end{center}

Now suppose $v_i\in U\cap V'$. In this case $a_\y^{v_i}$ and $b_\y^{v_i}$ are in the same component of $J_1$, $J_3$ and $J_6$,
and we have eight cases for the remaining $2$-factors. Namely, for each of $J_2$, $J_4$ and $J_5$, $a_\y^{v_i}$ and $b_\y^{v_i}$
are either in the same component or they are not. We number these eight cases as in the following table. 

\vspace{0.3cm}

\begin{center}
\begin{tabular}{c|c|c|c|}
&$J_2$&$J_4$&$J_5$\\
\hline
3.1& same&same&same\\
\hline
3.2& distinct&same&same\\
\hline
3.3& same&distinct&same\\
\hline
3.4& same&same&distinct\\
\hline
3.5& distinct&distinct&same\\
\hline
3.6& distinct&same&distinct\\
\hline
3.7& same&distinct&distinct\\
\hline
3.8& distinct&distinct&distinct\\
\hline
\end{tabular}
\end{center}

Depending on which of cases 3.1--3.8 that we are in, 
we can reallocate the edges of $\sigma_{v_i}(H)$ to the factors $J_1,J_2,\ldots,J_6$ 
as indicated in the following table to obtain a new 
$2$-factorisation of $L$ with the desired properties. If $J_x$ ($x\in\{1,2,\ldots,6\}$) is not listed for a particular case, then the edges of 
$\sigma_{v_i}(H)$ that are in $J_x$ are unchanged.

\begin{center}
\begin{tabular}{|c|l|}
\hline
3.1& \\
\hline
3.2& $J_2:\sigma_{v_i}([a_1,b_2,b_0,b_1] \cup [a_{\infty},c,b_{\infty}])$\\
& $J_5:\sigma_{v_i}([b_{\infty},b_1,c,b_0] \cup [a_0,b_2,a_{\infty}])$\\
\hline
3.3& $J_3:\sigma_{v_i}([a_{\infty},a_0,c,a_2,b_{\infty}] \cup [a_1,b_1])$\\
& $J_4:\sigma_{v_i}([a_1,a_2,a_{\infty}] \cup [b_{\infty},a_0,b_1])$\\
\hline
3.4& $J_2:\sigma_{v_i}([b_{\infty},c,b_0,b_1] \cup [a_1,b_2,a_{\infty}])$\\
& $J_5:\sigma_{v_i}([a_{\infty},c,b_1,b_{\infty}] \cup [a_0,b_2,b_0])$\\
\hline
3.5& $J_2:\sigma_{v_i}([a_1,b_2,b_0,b_1] \cup [a_{\infty},c,b_{\infty}])$\\
& $J_3:\sigma_{v_i}([a_{\infty},a_0,c,a_2,b_{\infty}] \cup [a_1,b_1])$\\
& $J_4:\sigma_{v_i}([a_1,a_2,a_{\infty}] \cup [b_{\infty},a_0,b_1])$\\
& $J_5:\sigma_{v_i}([b_{\infty},b_1,c,b_0] \cup [a_0,b_2,a_{\infty}])$\\
\hline
3.6& $J_2:\sigma_{v_i}([a_1,b_2,b_0,b_1] \cup [a_{\infty},c,b_{\infty}])$\\
& $J_5:\sigma_{v_i}([a_{\infty},b_1,c,b_0] \cup [a_0,b_2,b_{\infty}])$\\
& $J_6:\sigma_{v_i}([a_{\infty},b_2,b_1,b_{\infty}] \cup [a_0,b_0])$\\
\hline
3.7& $J_2:\sigma_{v_i}([b_{\infty},c,b_0,b_1] \cup [a_1,b_2,a_{\infty}])$\\
& $J_3:\sigma_{v_i}([a_{\infty},a_0,c,a_2,b_{\infty}] \cup [a_1,b_1])$\\
& $J_4:\sigma_{v_i}([a_1,a_2,a_{\infty}] \cup [b_{\infty},a_0,b_1])$\\
& $J_5:\sigma_{v_i}([a_{\infty},c,b_1,b_{\infty}] \cup [a_0,b_2,b_0])$\\
\hline
3.8& $J_2:\sigma_{v_i}(a_1,b_2,b_0,b_1] \cup [a_{\infty},c,b_{\infty}])$\\
& $J_3:\sigma_{v_i}([a_{\infty},a_0,c,a_2,b_{\infty}] \cup [a_1,b_1])$\\
& $J_4:\sigma_{v_i}([a_1,a_2,a_{\infty}] \cup [b_{\infty},a_0,b_1])$\\
& $J_5:\sigma_{v_i}([a_{\infty},b_1,c,b_0] \cup [a_0,b_2,b_{\infty}])$\\
& $J_6:\sigma_{v_i}([a_{\infty},b_2,b_1,b_{\infty}] \cup [a_0,b_0])$\\
\hline
\end{tabular}
\end{center}

Now suppose $v_i\in U'\cap V'$. In this case $a_\y^{v_i}$ and $b_\y^{v_i}$ are in the same component of $J_2$, $J_4$ and $J_6$,
and we have eight cases for the remaining $2$-factors. Namely, for each of $J_1$, $J_3$ and $J_5$, $a_\y^{v_i}$ and $b_\y^{v_i}$
are either in the same component or they are not. We number these eight cases as in the following table. 

\vspace{0.3cm}

\begin{center}
\begin{tabular}{c|c|c|c|}
&$J_1$&$J_3$&$J_5$\\
\hline
4.1& same&same&same\\
\hline
4.2& distinct&same&same\\
\hline
4.3& same&distinct&same\\
\hline
4.4& same&same&distinct\\
\hline
4.5& distinct&distinct&same\\
\hline
4.6& distinct&same&distinct\\
\hline
4.7& same&distinct&distinct\\
\hline
4.8& distinct&distinct&distinct\\
\hline
\end{tabular}
\end{center}

Depending on which of cases 4.1--4.8 that we are in, 
we can reallocate the edges of $\sigma_{v_i}(H)$ to the factors $J_1,J_2,\ldots,J_6$ 
as indicated in the following table to obtain a new 
$2$-factorisation of $L$ with the desired properties. If $J_x$ ($x\in\{1,2,\ldots,6\}$) is not listed for a particular case, then the edges of 
$\sigma_{v_i}(H)$ that are in $J_x$ are unchanged.

\begin{center}
\begin{tabular}{|c|l|}
\hline
4.1& \\
\hline
4.2& $J_1:\sigma_{v_i}(a_1,a_0,a_2,b_1] \cup [a_{\infty},c,b_{\infty}])$\\
& $J_4:\sigma_{v_i}([a_{\infty},a_2,c,a_0,b_{\infty}] \cup [a_1,b_1])$\\
\hline
4.3a& $J_1:\sigma_{v_i}([a_1,a_0,c,a_{\infty}] \cup [b_{\infty},a_2,b_1])$\\
& $J_3:\sigma_{v_i}([a_1,a_2,a_{\infty}] \cup [b_{\infty},a_0,b_1])$\\
& $J_4:\sigma_{v_i}([a_{\infty},a_0,a_2,c,b_{\infty}] \cup [a_1,b_1])$\\
& if $a_\y^{v_i}$ and $b_\y^{v_i}$ are in the same component of $J_1\setminus\sigma_{v_i}(H)$. \\
4.3b& $J_1:\sigma_{v_i}([a_{\infty},c,a_2,b_{\infty}] \cup [a_1,a_0,b_1])$\\
& $J_3:\sigma_{v_i}([a_{\infty},a_2,a_0,b_{\infty}] \cup [a_1,b_1])$\\
& $J_4:\sigma_{v_i}([a_{\infty},a_0,c,b_{\infty}] \cup [a_1,a_2,b_1])$\\
& if $a_\y^{v_i}$ and $a_1^{v_i}$ are in the same component of $J_1\setminus\sigma_{v_i}(H)$.\\
\hline
4.4& $J_5:\sigma_{v_i}([a_{\infty},b_1,c,b_0] \cup [a_0,b_2,b_{\infty}])$\\
& $J_6:\sigma_{v_i}([a_{\infty},b_2,b_1,b_{\infty}] \cup [a_0,b_0])$\\
\hline
4.5& $J_1:\sigma_{v_i}([a_{\infty},c,a_2,b_{\infty}] \cup [a_1,a_0,b_1])$\\
& $J_3:\sigma_{v_i}([a_{\infty},a_2,a_0,b_{\infty}] \cup [a_1,b_1])$\\
& $J_4:\sigma_{v_i}([a_{\infty},a_0,c,b_{\infty}] \cup [a_1,a_2,b_1])$\\
\hline
4.6& $J_1:\sigma_{v_i}([a_1,a_0,a_2,b_1] \cup [a_{\infty},c,b_{\infty}])$\\
& $J_4:\sigma_{v_i}([a_{\infty},a_2,c,a_0,b_{\infty}] \cup [a_1,b_1])$\\
& $J_5:\sigma_{v_i}([a_{\infty},b_1,c,b_0] \cup [a_0,b_2,b_{\infty}])$\\
& $J_6:\sigma_{v_i}([a_{\infty},b_2,b_1,b_{\infty}] \cup [a_0,b_0])$\\
\hline
4.7a& $J_1:\sigma_{v_i}([a_1,a_0,c,a_{\infty}] \cup [b_{\infty},a_2,b_1])$\\
& $J_3:\sigma_{v_i}([a_1,a_2,a_{\infty}] \cup [b_{\infty},a_0,b_1])$\\
& $J_4:\sigma_{v_i}([a_{\infty},a_0,a_2,c,b_{\infty}] \cup [a_1,b_1])$\\
& $J_5:\sigma_{v_i}([a_{\infty},b_1,c,b_0] \cup [a_0,b_2,b_{\infty}])$\\
& $J_6:\sigma_{v_i}([a_{\infty},b_2,b_1,b_{\infty}] \cup [a_0,b_0])$\\
& if $a_\y^{v_i}$ and $b_\y^{v_i}$ are in the same component of $J_1\setminus\sigma_{v_i}(H)$. \\
4.7b& $J_1:\sigma_{v_i}([a_{\infty},c,a_2,b_{\infty}] \cup [a_1,a_0,b_1])$\\
& $J_3:\sigma_{v_i}([a_{\infty},a_2,a_0,b_{\infty}] \cup [a_1,b_1])$\\
& $J_4:\sigma_{v_i}([a_{\infty},a_0,c,b_{\infty}] \cup [a_1,a_2,b_1])$\\
& $J_5:\sigma_{v_i}([a_{\infty},b_1,c,b_0] \cup [a_0,b_2,b_{\infty}])$\\
& $J_6:\sigma_{v_i}([a_{\infty},b_2,b_1,b_{\infty}] \cup [a_0,b_0])$\\
& if $a_\y^{v_i}$ and $a_1^{v_i}$ are in the same component of $J_1\setminus\sigma_{v_i}(H)$.\\
\hline
4.8& $J_1:\sigma_{v_i}([a_1,a_0,c,a_{\infty}] \cup [b_{\infty},a_2,b_1])$\\
& $J_3:\sigma_{v_i}([a_1,a_2,a_{\infty}] \cup [b_{\infty},a_0,b_1])$\\
& $J_4:\sigma_{v_i}([a_{\infty},a_0,a_2,c,b_{\infty}] \cup [a_1,b_1])$\\
& $J_5:\sigma_{v_i}([a_{\infty},b_1,c,b_0] \cup [a_0,b_2,b_{\infty}])$\\
& $J_6:\sigma_{v_i}([a_{\infty},b_2,b_1,b_{\infty}] \cup [a_0,b_0])$\\
\hline
\end{tabular}
\end{center}

\vspace{5cm}

The cases $n=5,6,7,8,9$ are similar and use the following lemma.

\begin{lemma}\label{4RegularHamiltonFragments}
If $H$ is a subgraph of $K_{A_n\cup B_n\cup\{c\}}$, $t\in \{0,1,\ldots,n-2\}$, 
and there exist four decompositions 
$\D_1=\{X_1,Y_1\}$, $\D'_1=\{X'_1,Y'_1\}$, $\D_2=\{X_2,Y_2\}$ and $\D'_2=\{X'_2,Y'_2\}$  of $H$ such that 
\begin{itemize}
\item[(a)] $^{A}X_1$ is a $2$-regular graph, $^{A}Y_1$ is a $2$-regular graph,
and
$\{V(^{A}X_1),V(^{A}Y_1)\}=\{N_n,N_n\cup\{c\}\}$;
\item[(b)] for each vertex $v$ in $K_{A_n\cup B_n}$, 
$\deg_{X_1}(v)=\deg_{X'_1}(v)=\deg_{X_2}(v)=\deg_{X'_2}(v)$ and $\deg_{Y_1}(v)=\deg_{Y'_1}(v)=\deg_{Y_2}(v)=\deg_{Y'_2}(v)$;
\item[(c)] $\deg_{X_1}(c)= \deg_{X_2}(c) = \deg_{Y'_1}(c)= \deg_{Y'_2}(c)$ and 
$\deg_{Y_1}(c)= \deg_{Y_2}(c) = \deg_{X'_1}(c)= \deg_{X'_2}(c)$;
\item[(d)] each of $X_1$ and $X_2$ is the vertex disjoint union of an $a_t,b_t$-path and an $a_\y,b_\y$-path; and
\begin{itemize}
\item[(d$_1$)] $Y_1$ is  
the vertex disjoint union of an $a_t,a_\y$-path and a $b_t,b_\y$-path,
and $Y_2$ is  
the vertex disjoint union of either an $a_t,b_\y$-path and a $b_t,a_\y$-path,
or an $a_t,b_t$-path and an $a_\y,b_\y$-path; or
\item[(d$_2$)] $Y_1$ is  
the vertex disjoint union of an $a_t,b_\y$-path and a $b_t,a_\y$-path,
and $Y_2$ is  
the vertex disjoint union of either an $a_t,a_\y$-path and a $b_t,b_\y$-path,
or an $a_t,b_t$-path and an $a_\y,b_\y$-path,
\end{itemize}
\item[(e)] each of $Y'_1$ and $Y'_2$ is the vertex disjoint union of an $a_t,b_t$-path and an $a_\y,b_\y$-path; and
\begin{itemize}
\item[(e$_1$)] $X'_1$ is  
the vertex disjoint union of an $a_t,a_\y$-path and a $b_t,b_\y$-path,
and $X'_2$ is  
the vertex disjoint union of either an $a_t,b_\y$-path and a $b_t,a_\y$-path,
or an $a_t,b_t$-path and an $a_\y,b_\y$-path; or
\item[(e$_2$)] $X'_1$ is  
the vertex disjoint union of an $a_t,b_\y$-path and a $b_t,a_\y$-path,
and $X'_2$ is  
the vertex disjoint union of either an $a_t,a_\y$-path and a $b_t,b_\y$-path,
or an $a_t,b_t$-path and an $a_\y,b_\y$-path,
\end{itemize}
\end{itemize}
then $H$ is a Hamilton fragment. 
\end{lemma}

\proof
Our aim is to show that the subgraph $$L=\bigcup_{v\in V(G)}\sigma_v(H)$$
of $L(G)$ decomposes into Hamilton cycles, and we do this by applying Lemma \ref{manyrepairs} 
with $H_i=\sigma_{v_i}(H)$ and $V_i=\{a_\infty^{v_i},b_\infty^{v_i}\}$ for each $i=1,2,\ldots,m$.
To this end, observe that \linebreak $\sigma_{v_1}(H),\sigma_{v_2}(H),\ldots,\sigma_{v_m}(H)$ are edge-disjoint subgraphs of $L$,
and since $b_\infty^{v_i}=a_\infty^{v_{i+1}}$ for $i=1,2,\ldots,m-1$ it follows that $V_i\cap V_{i+1}\ne \emptyset$ for $i=1,2,\ldots,m-1$ as required. It remains to show there is a $2$-factorisation
$\{J_1,J_2\}$ of $L$ such that 
\begin{itemize}
\item $\{a_\infty^{v_1},a_\infty^{v_2},\ldots,a_\infty^{v_m}\}$ links $\{J_1,J_2\}$; and
\item $\sigma_{v_i}(H)$ induces an $\{a_\infty^{v_{i}},b_\infty^{v_i}\}$-connector in $\{J_1,J_2\}$, for $i=1,2,\ldots,m$.
\end{itemize} 
Let $\{U,U'\}$ be a partition of $V(G)$ such that both $U$ and $U'$ link $\{F,\F_0\}$
(such a partition exists by Lemma \ref{intersects}),
and let $\{J_1,J_2\}$ be the decomposition of $L$ defined by 
$$J_1=\bigcup_{v\in U}\sigma_v(X_1)\ \cup \ 
\bigcup_{v\in U'}\sigma_v(X'_1)$$
and
$$J_2=\bigcup_{v\in U}\sigma_v(Y_1)\ \cup \ 
\bigcup_{v\in U'}\sigma_v(Y'_1).$$

It is easily checked that $\{J_1,J_2\}$ is a $2$-factorisation of $L$ and that 
$\{a_\infty^{v_1},a_\infty^{v_2},\ldots,a_\infty^{v_m}\}$ links $\{J_1,J_2\}$.
It remains to show that for $i=1,2,\ldots,m$, $\sigma_{v_i}(H)$ induces an $\{a_\infty^{v_{i}},b_\infty^{v_i}\}$-connector in 
$\{J_1,J_2\}$. 
There are two cases to consider: $v_i\in U$ and $v_i\in U'$.

First suppose $v_i\in U$. In this case $a_\y^{v_i}$ and $b_\y^{v_i}$ are in the same component of $J_1$.
If $a_\y^{v_i}$ and $b_\y^{v_i}$ are also in the same component of $J_2$, then we are done. 
If $a_\y^{v_i}$ and $b_\y^{v_i}$ are not in the same component of $J_2$, then reallocate the edges of 
$\sigma_{v_i}(H)$ to the factors $J_1$ and $J_2$ as follows to obtain a $2$-factorisation of $L$ with the desired properties. 
$$J_1:\sigma_{v_i}(X_2)\qquad J_2:\sigma_{v_i}(Y_2)$$
 
Now suppose $v_i\in U'$. In this case $a_\y^{v_i}$ and $b_\y^{v_i}$ are in the same component of $J_2$.
If $a_\y^{v_i}$ and $b_\y^{v_i}$ are also in the same component of $J_1$, then we are done. 
If $a_\y^{v_i}$ and $b_\y^{v_i}$ are not in the same component of $J_1$, then reallocate the edges of 
$\sigma_{v_i}(H)$ to the factors $J_1$ and $J_2$ as follows to obtain a $2$-factorisation of $L$ with the desired properties. 
$$J_1:\sigma_{v_i}(X'_2)\qquad J_2:\sigma_{v_i}(Y'_2)$$

\qed

\noindent{\bf The case n=5:}

For $i=1,2$, let $X_1(i)$, $Y_1(i)$, $X'_1(i)$, $Y'_1(i)$, $X_2(i)$, $Y_2(i)$, $X'_2(i)$, $Y'_2(i)$ be the subgraphs of $K_{A_n\cup B_n\cup\{c\}}$ given by the union of
the paths listed in the following tables, and let $H^i$ be the subgraph of $K_{A_n\cup B_n\cup\{c\}}$ with edge set
$E(X_1(i))\cup E(Y_1(i))$. Applying Lemma \ref{4RegularHamiltonFragments} with $H=H^i$, $X_1=X_1(i)$, $Y_1=Y_1(i)$, $X'_1=X'_1(i)$, $Y'_1=Y_1(i)$,$X_2=X_2(i)$, $Y_2=Y_2(i)$, $X'_2=X_2(i)$, $Y'_2=Y_2(i)$  shows that each 
$H^i$ is a Hamilton fragment. The value of $t$ can be deduced from the ends of the given paths.

\vspace{0.3cm}

$\begin{array}{|c|c|}
\hline
X_1(1)& [a_{\infty},b_1,c,a_2,a_0,b_{\infty}] \cup [a_3,b_3]\\
\hline
Y_1(1)&[b_{\infty},b_1,a_2,b_3] \cup [a_3,a_0,a_{\infty}]\\
\hline
X'_1(1)& Y_1(1) \\
\hline
Y'_1(1) & X_1(1)\\
\hline
X_2(1)& [a_{\infty},a_0,a_2,c,b_1,b_{\infty}] \cup [a_3,b_3]\\
\hline
Y_2(1)& [a_{\infty},b_1,a_2,b_3] \cup [a_3,a_0,b_{\infty}]\\
\hline
X'_2(1) & X_2(1) \\
\hline
Y'_2(1) & Y_2(1) \\
\hline
\end{array}
$
$\begin{array}{|c|c|}
\hline
X_1(2)& [a_{\infty},a_1,b_0,c,b_3,b_{\infty}] \cup [a_2,b_2]\\
\hline
Y_1(2)&[a_{\infty},b_3,b_0,b_2] \cup [a_2,a_1,b_{\infty}]\\
\hline
X'_1(2)& Y_1(2) \\
\hline
Y'_1(2) & X_1(2)\\
\hline
X_2(2)& [a_{\infty},b_3,c,b_0,a_1,b_{\infty}] \cup [a_2,b_2]\\
\hline
Y_2(2)& [b_{\infty},b_3,b_0,b_2] \cup [a_2,a_1,a_{\infty}]\\
\hline
X'_2(2) & X_2(2) \\
\hline
Y'_2(2) & Y_2(2) \\
\hline
\end{array}
$

Note that $H^1$ is edge disjoint from $H^2$ and let $H'=K_{A_5\cup B_5\cup \{c\}}-{H^1\cup H^2}$ so that $\{H', H^1, H^2\}$ is a decomposition of $K_{A_5\cup B_5\cup \{c\}}\cong K_{11}$. We now show that $H'$ is a Hamilton fragment. 

Our aim is to show that the subgraph $$L=\bigcup_{v\in V(G)}\sigma_v(H')$$
of $L(G)$ decomposes into Hamilton cycles, and we do this by applying Lemma \ref{manyrepairs} 
with $H_i=\sigma_{v_i}(H')$ and $V_i=\{a_\infty^{v_i},b_\infty^{v_i}\}$ for each $i=1,2,\ldots,m$.
To this end, observe that \linebreak $\sigma_{v_1}(H'),\sigma_{v_2}(H'),\ldots,\sigma_{v_m}(H')$ are edge-disjoint subgraphs of $L$,
and since $b_\infty^{v_i}=a_\infty^{v_{i+1}}$ for $i=1,2,\ldots,m-1$ it follows that $V_i\cap V_{i+1}\ne \emptyset$ for $i=1,2,\ldots,m-1$ as required. It remains to show there is a $2$-factorisation
$\{J_1,J_2,J_3,J_4,J_5,J_6\}$ of $L$ such that 
\begin{itemize}
\item $\{a_\infty^{v_1},a_\infty^{v_2},\ldots,a_\infty^{v_m}\}$ links $\{J_1,J_2,J_3,J_4,J_5,J_6\}$; and
\item $\sigma_{v_i}(H')$ induces an $\{a_\infty^{v_{i}},b_\infty^{v_i}\}$-connector in $\{J_1,J_2,J_3,J_4,J_5,J_6\}$, for $i=1,2,\ldots,m$.
\end{itemize} 
Let $\{U,U'\}$ be a partition of $V(G)$ such that both $U$ and $U'$ link $\{F,\F_1\}$,
let $\{V,V'\}$ be a partition of $V(G)$ such that both $V$ and $V'$ link $\{F,\F_2\}$ 
(such partitions exist by Lemma \ref{intersects}),
and let $\{J_1,J_2,J_3,J_4,J_5,J_6\}$ be the decomposition of $L$ defined by 

\begin{itemize}
\item 
$
\begin{array}[t]{lll}
J_1&=&\bigcup_{v\in U\cap V}\sigma_v([[a_1,b_2,a_0,b_3,b_1] \cup [a_{\infty},b_{\infty}])\ \cup \\
&&\bigcup_{v\in U'\cap V}\sigma_v([b_{\infty},c,a_0,b_1]\cup[a_1,b_3,b_2,a_{\infty}])\ \cup \\
&&\bigcup_{v\in U\cap V'}\sigma_v([a_1,b_3,a_0,b_2,b_1]\cup [a_{\infty},b_{\infty}])\ \cup \\
&&\bigcup_{v\in U'\cap V'}\sigma_v([a_1,b_3,a_0,c,a_{\infty}]\cup[b_{\infty},b_2,b_1]);
\end{array}
$
\item 
$
\begin{array}[t]{lll}
J_2&=&\bigcup_{v\in U\cap V}\sigma_v([b_{\infty},c,a_0,b_1] \cup [a_1,b_3,b_2,a_{\infty}])\ \cup \\
&&\bigcup_{v\in U'\cap V}\sigma_v([a_1,b_2,a_0,b_3,b_1]\cup[a_{\infty},b_{\infty}])\ \cup \\
&&\bigcup_{v\in U\cap V'}\sigma_v([a_1,a_0,c,b_{\infty}]\cup[a_{\infty},b_2,b_3,b_1])\ \cup \\
&&\bigcup_{v\in U'\cap V'}\sigma_v([a_1,a_0,b_2,b_3,b_1]\cup[a_{\infty},b_{\infty}]);
\end{array}
$
\item 
$
\begin{array}[t]{lll}
J_3&=&\bigcup_{v\in U\cap V}\sigma_v([a_{\infty},a_2,b_0,a_3,b_{\infty}] \cup [a_1,b_1])\ \cup \\
&&\bigcup_{v\in U'\cap V}\sigma_v([a_1,a_3,c,a_{\infty}]\cup[b_{\infty},a_2,b_0,b_1])\ \cup \\
&&\bigcup_{v\in U\cap V'}\sigma_v([a_{\infty},a_2,b_0,a_3,b_{\infty}]\cup[a_1,b_1])\ \cup \\
&&\bigcup_{v\in U'\cap V'}\sigma_v([a_1,c,a_3,a_{\infty}]\cup[b_{\infty},a_2,b_0,b_1]);
\end{array}
$
\item 
$
\begin{array}[t]{lll}
J_4&=&\bigcup_{v\in U\cap V}\sigma_v([a_1,c,a_3,a_2,b_{\infty}] \cup [a_{\infty},b_0,b_1])\ \cup \\
&&\bigcup_{v\in U'\cap V}\sigma_v([a_{\infty},a_2,a_3,b_0,b_{\infty}]\cup[a_1,b_1])\ \cup \\
&&\bigcup_{v\in U\cap V'}\sigma_v([a_1,c,a_3,a_2,b_{\infty}]\cup[a_{\infty},b_0,b_1])\ \cup \\
&&\bigcup_{v\in U'\cap V'}\sigma_v([a_{\infty},a_2,a_3,b_0,b_{\infty}]\cup[a_1,b_1]);
\end{array}
$
\item 
$
\begin{array}[t]{lll}
J_5&=&\bigcup_{v\in U\cap V}\sigma_v([a_{\infty},a_3,b_1,b_2,b_{\infty}]\cup [a_0,b_0])\ \cup \\
&&\bigcup_{v\in U'\cap V}\sigma_v([a_{\infty},a_3,b_1,b_2,b_{\infty}]\cup[a_0,b_0])\ \cup \\
&&\bigcup_{v\in U\cap V'}\sigma_v([a_0,b_1,a_3,b_2,c,a_{\infty}]\cup[b_{\infty},b_0])\ \cup \\
&&\bigcup_{v\in U'\cap V'}\sigma_v([a_0,b_1,a_3,b_2,c,b_{\infty}]\cup[a_{\infty},b_0]);
\end{array}
$
\item 
$
\begin{array}[t]{lll}
J_6&=&\bigcup_{v\in U\cap V}\sigma_v([a_0,a_1,a_3,b_2,c,a_{\infty}] \cup [b_{\infty},b_0])\ \cup \\
&&\bigcup_{v\in U'\cap V}\sigma_v([a_0,a_1,c,b_2,a_3,b_{\infty}]\cup[a_{\infty},b_0])\ \cup \\
&&\bigcup_{v\in U\cap V'}\sigma_v([a_{\infty},a_3,a_1,b_2,b_{\infty}]\cup[a_0,b_0])\ \cup \\
&&\bigcup_{v\in U'\cap V'}\sigma_v([a_{\infty},b_2,a_1,a_3,b_{\infty}]\cup[a_0,b_0]);
\end{array}
$

\end{itemize}
It is easily checked that $\{J_1,J_2,J_3,J_4,J_5,J_6\}$ is a $2$-factorisation of $L$ and that $\{a_\infty^{v_1},a_\infty^{v_2},\ldots,a_\infty^{v_m}\}$ links $\{J_1,J_2,J_3,J_4,J_5,J_6\}$.
It remains to show that for $i=1,2,\ldots,m$, $\sigma_{v_i}(H')$ induces an $\{a_\infty^{v_{i}},b_\infty^{v_i}\}$-connector in $\{J_1,J_2,J_3,J_4,J_5,J_6\}$. 
There are four cases to consider: $v_i\in U\cap V$, $v_i\in U'\cap V$, $v_i\in U\cap V'$ and $v_i\in U'\cap V'$.

First suppose $v_i\in U\cap V$. In this case $a_\y^{v_i}$ and $b_\y^{v_i}$ are in the same component of $J_1$, $J_3$ and $J_5$,
and we have eight cases for the remaining $2$-factors. Namely, for each of $J_2$, $J_4$ and $J_6$, $a_\y^{v_i}$ and $b_\y^{v_i}$
are either in the same component or they are not. We number these eight cases as in the following table. 

\vspace{0.3cm}

\begin{center}
\begin{tabular}{c|c|c|c|}
&$J_2$&$J_4$&$J_6$\\
\hline
1.1& same&same&same\\
\hline
1.2& distinct&same&same\\
\hline
1.3& same&distinct&same\\
\hline
1.4& same&same&distinct\\
\hline
1.5& same&distinct&distinct\\
\hline
1.6& distinct&same&distinct\\
\hline
1.7& distinct&distinct&same\\
\hline
1.8& distinct&distinct&distinct\\
\hline
\end{tabular}
\end{center}

Depending on which of cases 1.1--1.8 that we are in, 
we can reallocate the edges of $\sigma_{v_i}(H)$ to the factors $J_1,J_2,\ldots,J_6$ 
as indicated in the following table to obtain a new 
$2$-factorisation of $L$ with the desired properties. If $J_x$ ($x\in\{1,2,\ldots,6\}$) is not listed for a particular case, then the edges of 
$\sigma_{v_i}(H)$ that are in $J_x$ are unchanged. 

\begin{center}
\begin{tabular}{|c|l|}
\hline
1.1& \\
\hline
1.2& $J_1:\sigma_{v_i}([a_1,b_2,b_3,a_0,b_1] \cup [a_{\infty},b_{\infty}])$\\
& $J_2:\sigma_{v_i}([a_{\infty},b_2,a_0,c,b_{\infty}] \cup [a_1,b_3,b_1])$\\
\hline
1.3& $J_3:\sigma_{v_i}([a_{\infty},b_0,a_2,a_3,b_{\infty}] \cup [a_1,b_1])$\\
& $J_4:\sigma_{v_i}([a_1,c,a_3,b_0,b_1] \cup [a_{\infty},a_2,b_{\infty}])$\\
\hline
1.4a& $J_3:\sigma_{v_i}([a_1,a_3,b_0,b_1] \cup [a_{\infty},a_2,b_{\infty}])$\\
& $J_4:\sigma_{v_i}([a_{\infty},c,a_3,a_2,b_0,b_{\infty}] \cup [a_1,b_1])$\\
& $J_6:\sigma_{v_i}([a_0,a_1,c,b_2,a_3,b_{\infty}] \cup [a_{\infty},b_0])$\\
& if $a_\y^{v_i}$ and $a_1^{v_i}$ are in the same component of $J_4\setminus\sigma_{v_i}(H')$.\\
1.4b& $J_3:\sigma_{v_i}([a_{\infty},a_2,a_3,b_0,b_{\infty}] \cup [a_1,b_1])$\\
& $J_4:\sigma_{v_i}([a_1,a_3,c,a_{\infty}] \cup [b_{\infty},a_2,b_0,b_1])$\\
& $J_6:\sigma_{v_i}([a_0,a_1,c,b_2,a_3,b_{\infty}] \cup [a_{\infty},b_0])$\\
& if $a_\y^{v_i}$ and $b_\y^{v_i}$ are in the same component of $J_4\setminus\sigma_{v_i}(H')$.\\
\hline
1.5& $J_3:\sigma_{v_i}([a_1,a_3,b_0,b_1] \cup [a_{\infty},a_2,b_{\infty}])$\\
& $J_4:\sigma_{v_i}([a_{\infty},c,a_3,a_2,b_0,b_{\infty}] \cup [a_1,b_1])$\\
& $J_6:\sigma_{v_i}([a_0,a_1,c,b_2,a_3,b_{\infty}] \cup [a_{\infty},b_0])$\\
\hline
1.6a& $J_1:\sigma_{v_i}([a_1,a_0,b_3,b_2,b_1] \cup [a_{\infty},b_{\infty}])$\\
& $J_2:\sigma_{v_i}([a_{\infty},b_2,a_0,c,b_{\infty}] \cup [a_1,b_3,b_1])$\\
& $J_3:\sigma_{v_i}([a_{\infty},a_2,a_3,b_0,b_{\infty}] \cup [a_1,b_1])$\\
& $J_4:\sigma_{v_i}([a_1,c,a_3,b_1] \cup [a_{\infty},b_0,a_2,b_{\infty}])$\\
& $J_5:\sigma_{v_i}([a_{\infty},a_3,b_2,b_{\infty}] \cup [a_0,b_1,b_0])$\\
& $J_6:\sigma_{v_i}([a_{\infty},c,b_2,a_1,a_3,b_{\infty}] \cup [a_0,b_0])$\\
& if $a_\y^{v_i}$ and $a_1^{v_i}$ are in the same component of $J_4\setminus\sigma_{v_i}(H')$.\\
1.6b& $J_1:\sigma_{v_i}([a_1,b_2,b_3,a_0,b_1] \cup [a_{\infty},b_{\infty}])$\\
& $J_2:\sigma_{v_i}([a_{\infty},b_2,a_0,c,b_{\infty}] \cup [a_1,b_3,b_1])$\\
& $J_3:\sigma_{v_i}([a_{\infty},a_2,a_3,b_0,b_{\infty}] \cup [a_1,b_1])$\\
& $J_4:\sigma_{v_i}([a_1,a_3,c,a_{\infty}] \cup [b_{\infty},a_2,b_0,b_1])$\\
& $J_6:\sigma_{v_i}([a_0,a_1,c,b_2,a_3,b_{\infty}] \cup [a_{\infty},b_0])$\\
& if $a_\y^{v_i}$ and $b_\y^{v_i}$ are in the same component of $J_4\setminus\sigma_{v_i}(H')$.\\
\hline
1.7& $J_1:\sigma_{v_i}([a_1,b_2,b_3,a_0,b_1] \cup [a_{\infty},b_{\infty}])$\\
& $J_2:\sigma_{v_i}([a_{\infty},b_2,a_0,c,b_{\infty}] \cup [a_1,b_3,b_1])$\\
& $J_3:\sigma_{v_i}([a_{\infty},b_0,a_2,a_3,b_{\infty}] \cup [a_1,b_1])$\\
& $J_4:\sigma_{v_i}([a_1,c,a_3,b_0,b_1] \cup [a_{\infty},a_2,b_{\infty}])$\\
& $J_6:\sigma_{v_i}([a_0,a_1,a_3,b_2,c,a_{\infty}] \cup [b_{\infty},b_0])$\\
\hline
1.8& $J_1:\sigma_{v_i}([a_1,b_2,b_3,a_0,b_1] \cup [a_{\infty},b_{\infty}])$\\
& $J_2:\sigma_{v_i}([a_{\infty},b_2,a_0,c,b_{\infty}] \cup [a_1,b_3,b_1])$\\
& $J_3:\sigma_{v_i}([a_{\infty},a_2,a_3,b_0,b_{\infty}] \cup [a_1,b_1])$\\
& $J_4:\sigma_{v_i}([a_1,a_3,c,a_{\infty}] \cup [b_{\infty},a_2,b_0,b_1])$\\
& $J_6:\sigma_{v_i}([a_0,a_1,c,b_2,a_3,b_{\infty}] \cup [a_{\infty},b_0])$\\
\hline

\end{tabular}
\end{center}

Next suppose $v_i\in U'\cap V$. In this case $a_\y^{v_i}$ and $b_\y^{v_i}$ are in the same component of $J_2$, $J_4$ and $J_5$,
and we have eight cases for the remaining $2$-factors. Namely, for each of $J_1$, $J_3$ and $J_6$, $a_\y^{v_i}$ and $b_\y^{v_i}$
are either in the same component or they are not. We number these eight cases as in the following table. 

\vspace{0.3cm}

\begin{center}
\begin{tabular}{c|c|c|c|}
&$J_1$&$J_3$&$J_6$\\
\hline
2.1& same&same&same\\
\hline
2.2& distinct&same&same\\
\hline
2.3& same&distinct&same\\
\hline
2.4& same&same&distinct\\
\hline
2.5& same&distinct&distinct\\
\hline
2.6& distinct&same&distinct\\
\hline
2.7& distinct&distinct&same\\
\hline
2.8& distinct&distinct&distinct\\
\hline
\end{tabular}
\end{center}

Depending on which of cases 2.1--2.8 that we are in, 
we can reallocate the edges of $\sigma_{v_i}(H)$ to the factors $J_1,J_2,\ldots,J_6$ 
as indicated in the following table to obtain a new 
$2$-factorisation of $L$ with the desired properties. If $J_x$ ($x\in\{1,2,\ldots,6\}$) is not listed for a particular case, then the edges of 
$\sigma_{v_i}(H')$ that are in $J_x$ are unchanged. 

\begin{center}
\begin{tabular}{|c|l|}
\hline
2.1& \\
\hline
2.2& $J_1:\sigma_{v_i}([a_{\infty},c,a_0,b_1] \cup [a_1,b_3,b_2,b_{\infty}])$\\
& $J_3:\sigma_{v_i}([a_1,c,a_3,a_{\infty}] \cup [b_{\infty},a_2,b_0,b_1])$\\
& $J_5:\sigma_{v_i}([a_{\infty},b_2,b_1,a_3,b_{\infty}] \cup [a_0,b_0])$\\
& $J_6:\sigma_{v_i}([a_0,a_1,a_3,b_2,c,b_{\infty}] \cup [a_{\infty},b_0])$\\
\hline
2.3& $J_3:\sigma_{v_i}([a_{\infty},c,a_3,a_2,b_0,b_{\infty}] \cup [a_1,b_1])$\\
& $J_4:\sigma_{v_i}([a_1,a_3,b_0,b_1] \cup [a_{\infty},a_2,b_{\infty}])$\\
\hline
2.4a& $J_1:\sigma_{v_i}([a_1,b_3,a_0,c,a_{\infty}] \cup [b_{\infty},b_2,b_1])$\\
& $J_2:\sigma_{v_i}([a_1,a_0,b_2,b_3,b_1] \cup [a_{\infty},b_{\infty}])$\\
& $J_3:\sigma_{v_i}([a_1,c,a_3,b_1] \cup [a_{\infty},b_0,a_2,b_{\infty}])$\\
& $J_5:\sigma_{v_i}([a_{\infty},b_2,a_3,b_{\infty}] \cup [a_0,b_1,b_0])$\\
& $J_6:\sigma_{v_i}([a_{\infty},a_3,a_1,b_2,c,b_{\infty}] \cup [a_0,b_0])$\\
& if $a_\y^{v_i}$ and $a_1^{v_i}$ are in the same component of $J_3\setminus\sigma_{v_i}(H')$.\\
2.4b & $J_3:\sigma_{v_i}([a_1,c,a_3,b_{\infty}] \cup [a_{\infty},a_2,b_0,b_1])$\\
& $J_4:\sigma_{v_i}([a_{\infty},b_0,a_3,a_2,b_{\infty}] \cup [a_1,b_1])$\\
& $J_6:\sigma_{v_i}([a_0,a_1,a_3,b_2,c,a_{\infty}] \cup [b_{\infty},b_0])$\\
& if $a_\y^{v_i}$ and $b_\y^{v_i}$ are in the same component of $J_3\setminus\sigma_{v_i}(H')$.\\
\hline
2.5 & $J_3:\sigma_{v_i}([a_1,c,a_3,b_{\infty}] \cup [a_{\infty},a_2,b_0,b_1])$\\
& $J_4:\sigma_{v_i}([a_{\infty},b_0,a_3,a_2,b_{\infty}] \cup [a_1,b_1])$\\
& $J_6:\sigma_{v_i}([a_0,a_1,a_3,b_2,c,a_{\infty}] \cup [b_{\infty},b_0])$\\
\hline
2.6a& $J_1:\sigma_{v_i}([a_{\infty},b_2,a_0,c,b_{\infty}] \cup [a_1,b_3,b_1])$\\
& $J_2:\sigma_{v_i}([a_1,a_0,b_3,b_2,b_1] \cup [a_{\infty},b_{\infty}])$\\
& $J_3:\sigma_{v_i}([a_1,c,a_3,b_1] \cup [a_{\infty},b_0,a_2,b_{\infty}])$\\
& $J_5:\sigma_{v_i}([a_{\infty},a_3,b_2,b_{\infty}] \cup [a_0,b_1,b_0])$\\
& $J_6:\sigma_{v_i}([a_{\infty},c,b_2,a_1,a_3,b_{\infty}] \cup [a_0,b_0])$\\
& if $a_\y^{v_i}$ and $b_1^{v_i}$ are in the same component of $J_3\setminus\sigma_{v_i}(H')$.\\
2.6b &  $J_1:\sigma_{v_i}([a_{\infty},c,a_0,b_1] \cup [a_1,b_3,b_2,b_{\infty}])$\\
& $J_3:\sigma_{v_i}([a_1,a_3,c,b_{\infty}] \cup [a_{\infty},a_2,b_0,b_1])$\\
& $J_4:\sigma_{v_i}([a_{\infty},b_0,a_3,a_2,b_{\infty}] \cup [a_1,b_1])$\\
& $J_5:\sigma_{v_i}([a_{\infty},b_2,b_1,a_3,b_{\infty}] \cup [a_0,b_0])$\\
& $J_6:\sigma_{v_i}([a_0,a_1,c,b_2,a_3,a_{\infty}] \cup [b_{\infty},b_0])$\\
& if $a_\y^{v_i}$ and $b_\y^{v_i}$ are in the same component of $J_3\setminus\sigma_{v_i}(H')$.\\
\hline
\end{tabular}
\end{center}

\begin{center}
\begin{tabular}{|c|l|}
\hline
2.7& $J_1:\sigma_{v_i}([a_{\infty},b_2,a_0,c,b_{\infty}] \cup [a_1,b_3,b_1])$\\
& $J_2:\sigma_{v_i}([a_1,b_2,b_3,a_0,b_1] \cup [a_{\infty},b_{\infty}])$\\
& $J_3:\sigma_{v_i}([a_{\infty},c,a_3,a_2,b_0,b_{\infty}] \cup [a_1,b_1])$\\
& $J_4:\sigma_{v_i}([a_1,a_3,b_0,b_1] \cup [a_{\infty},a_2,b_{\infty}])$\\
\hline
2.8& $J_1:\sigma_{v_i}([a_{\infty},b_2,a_0,c,b_{\infty}] \cup [a_1,b_3,b_1])$\\
& $J_2:\sigma_{v_i}([a_1,a_0,b_3,b_2,b_1] \cup [a_{\infty},b_{\infty}])$\\
& $J_3:\sigma_{v_i}([a_1,c,a_3,b_1] \cup [a_{\infty},b_0,a_2,b_{\infty}])$\\
& $J_5:\sigma_{v_i}([a_{\infty},a_3,b_2,b_{\infty}] \cup [a_0,b_1,b_0])$\\
& $J_6:\sigma_{v_i}([a_{\infty},c,b_2,a_1,a_3,b_{\infty}] \cup [a_0,b_0])$\\
\hline
\end{tabular}
\end{center}

Next suppose $v_i\in U\cap V'$. In this case $a_\y^{v_i}$ and $b_\y^{v_i}$ are in the same component of $J_1$, $J_3$ and $J_6$,
and we have eight cases for the remaining $2$-factors. Namely, for each of $J_2$, $J_4$ and $J_5$, $a_\y^{v_i}$ and $b_\y^{v_i}$
are either in the same component or they are not. We number these eight cases as in the following table. 

\vspace{0.3cm}

\begin{center}
\begin{tabular}{c|c|c|c|}
&$J_2$&$J_4$&$J_5$\\
\hline
3.1& same&same&same\\
\hline
3.2& distinct&same&same\\
\hline
3.3& same&distinct&same\\
\hline
3.4& same&same&distinct\\
\hline
3.5& same&distinct&distinct\\
\hline
3.6& distinct&same&distinct\\
\hline
3.7& distinct&distinct&same\\
\hline
3.8& distinct&distinct&distinct\\
\hline
\end{tabular}
\end{center}

Depending on which of cases 3.1--3.8 that we are in, 
we can reallocate the edges of $\sigma_{v_i}(H)$ to the factors $J_1,J_2,\ldots,J_6$ 
as indicated in the following table to obtain a new 
$2$-factorisation of $L$ with the desired properties. If $J_x$ ($x\in\{1,2,\ldots,6\}$) is not listed for a particular case, then the edges of 
$\sigma_{v_i}(H')$ that are in $J_x$ are unchanged. 

\begin{center}
\begin{tabular}{|c|l|}
\hline
3.1& \\
\hline
3.2& $J_1:\sigma_{v_i}([a_1,a_0,b_3,b_2,b_1] \cup [a_{\infty},b_{\infty}])$\\
& $J_2:\sigma_{v_i}([a_{\infty},b_2,a_0,c,b_{\infty}] \cup [a_1,b_3,b_1])$\\
\hline
3.3& $J_3:\sigma_{v_i}([a_{\infty},b_0,a_2,a_3,b_{\infty}] \cup [a_1,b_1])$\\
& $J_4:\sigma_{v_i}([a_1,c,a_3,b_0,b_1] \cup [a_{\infty},a_2,b_{\infty}])$\\
\hline
3.4a& $J_1:\sigma_{v_i}([a_1,a_0,b_3,b_2,b_1] \cup [a_{\infty},b_{\infty}])$\\
& $J_2:\sigma_{v_i}([a_{\infty},b_2,a_0,c,b_{\infty}] \cup [a_1,b_3,b_1])$\\
& $J_3:\sigma_{v_i}([a_{\infty},a_2,a_3,b_0,b_{\infty}] \cup [a_1,b_1])$\\
& $J_4:\sigma_{v_i}([a_1,c,a_3,b_1] \cup [a_{\infty},b_0,a_2,b_{\infty}])$\\
& $J_5:\sigma_{v_i}([a_{\infty},c,b_2,a_3,b_{\infty}] \cup [a_0,b_1,b_0])$\\
& if $a_\y^{v_i}$ and $a_1^{v_i}$ are in the same component of $J_2\setminus\sigma_{v_i}(H')$ and \\
& $a_\y^{v_i}$ and $a_1^{v_i}$ are in the same component of $J_4\setminus\sigma_{v_i}(H')$.\\
3.4b& $J_1:\sigma_{v_i}([a_1,a_0,b_3,b_2,b_1] \cup [a_{\infty},b_{\infty}])$\\
& $J_2:\sigma_{v_i}([a_{\infty},c,a_0,b_2,b_{\infty}] \cup [a_1,b_3,b_1])$\\
& $J_3:\sigma_{v_i}([a_{\infty},a_2,a_3,b_0,b_{\infty}] \cup [a_1,b_1])$\\
& $J_4:\sigma_{v_i}([a_1,c,a_3,a_{\infty}] \cup [b_{\infty},a_2,b_0,b_1])$\\
& $J_5:\sigma_{v_i}([a_0,b_1,a_3,b_2,c,b_{\infty}] \cup [a_{\infty},b_0])$\\
& $J_6:\sigma_{v_i}([a_{\infty},b_2,a_1,a_3,b_{\infty}] \cup [a_0,b_0])$\\
& if $a_\y^{v_i}$ and $a_1^{v_i}$ are in the same component of $J_2\setminus\sigma_{v_i}(H')$ and \\
& $a_\y^{v_i}$ and $b_\y^{v_i}$ are in the same component of $J_4\setminus\sigma_{v_i}(H')$.\\
3.4c& $J_2:\sigma_{v_i}([a_1,a_0,c,a_{\infty}] \cup [b_{\infty},b_2,b_3,b_1])$\\
& $J_3:\sigma_{v_i}([a_{\infty},a_3,a_2,b_0,b_{\infty}] \cup [a_1,b_1])$\\
& $J_4:\sigma_{v_i}([a_1,c,a_3,b_0,b_1] \cup [a_{\infty},a_2,b_{\infty}])$\\
& $J_5:\sigma_{v_i}([a_0,b_1,a_3,b_2,c,b_{\infty}] \cup [a_{\infty},b_0])$\\
& $J_6:\sigma_{v_i}([a_{\infty},b_2,a_1,a_3,b_{\infty}] \cup [a_0,b_0])$\\
& if $a_\y^{v_i}$ and $b_\y^{v_i}$ are in the same component of $J_2\setminus\sigma_{v_i}(H')$ and \\
& $a_\y^{v_i}$ and $a_1^{v_i}$ are in the same component of $J_4\setminus\sigma_{v_i}(H')$.\\
3.4d& $J_2:\sigma_{v_i}([a_1,a_0,c,a_{\infty}] \cup [b_{\infty},b_2,b_3,b_1])$\\
& $J_3:\sigma_{v_i}([a_{\infty},a_2,a_3,b_0,b_{\infty}] \cup [a_1,b_1])$\\
& $J_4:\sigma_{v_i}([a_1,c,a_3,a_{\infty}] \cup [b_{\infty},a_2,b_0,b_1])$\\
& $J_5:\sigma_{v_i}([a_0,b_1,a_3,b_2,c,b_{\infty}] \cup [a_{\infty},b_0])$\\
& $J_6:\sigma_{v_i}([a_{\infty},b_2,a_1,a_3,b_{\infty}] \cup [a_0,b_0])$\\
& if $a_\y^{v_i}$ and $b_\y^{v_i}$ are in the same component of $J_2\setminus\sigma_{v_i}(H')$ and \\
& $a_\y^{v_i}$ and $b_\y^{v_i}$ are in the same component of $J_4\setminus\sigma_{v_i}(H')$.\\
\hline
\end{tabular}
\end{center}

\begin{center}
\begin{tabular}{|c|l|}
\hline
3.5a & $J_1:\sigma_{v_i}([a_1,a_0,b_3,b_2,b_1] \cup [a_{\infty},b_{\infty}])$\\
& $J_2:\sigma_{v_i}([a_{\infty},b_2,a_0,c,b_{\infty}] \cup [a_1,b_3,b_1])$\\
& $J_3:\sigma_{v_i}([a_{\infty},a_2,a_3,b_0,b_{\infty}] \cup [a_1,b_1])$\\
& $J_4:\sigma_{v_i}([a_1,c,a_3,b_1] \cup [a_{\infty},b_0,a_2,b_{\infty}])$\\
& $J_5:\sigma_{v_i}([a_{\infty},c,b_2,a_3,b_{\infty}] \cup [a_0,b_1,b_0])$\\
& if $a_\y^{v_i}$ and $a_1^{v_i}$ are in the same component of $J_2\setminus\sigma_{v_i}(H')$.\\
3.5b& $J_2:\sigma_{v_i}([a_1,a_0,c,a_{\infty}] \cup [b_{\infty},b_2,b_3,b_1])$\\
& $J_3:\sigma_{v_i}([a_{\infty},a_3,a_2,b_0,b_{\infty}] \cup [a_1,b_1])$\\
& $J_4:\sigma_{v_i}([a_1,c,a_3,b_0,b_1] \cup [a_{\infty},a_2,b_{\infty}])$\\
& $J_5:\sigma_{v_i}([a_0,b_1,a_3,b_2,c,b_{\infty}] \cup [a_{\infty},b_0])$\\
& $J_6:\sigma_{v_i}([a_{\infty},b_2,a_1,a_3,b_{\infty}] \cup [a_0,b_0])$\\
& if $a_\y^{v_i}$ and $b_\y^{v_i}$ are in the same component of $J_2\setminus\sigma_{v_i}(H')$.\\
\hline
3.6a & $J_1:\sigma_{v_i}([a_1,a_0,b_3,b_2,b_1] \cup [a_{\infty},b_{\infty}])$\\
& $J_2:\sigma_{v_i}([a_{\infty},b_2,a_0,c,b_{\infty}] \cup [a_1,b_3,b_1])$\\
& $J_3:\sigma_{v_i}([a_{\infty},a_2,a_3,b_0,b_{\infty}] \cup [a_1,b_1])$\\
& $J_4:\sigma_{v_i}([a_1,c,a_3,b_1] \cup [a_{\infty},b_0,a_2,b_{\infty}])$\\
& $J_5:\sigma_{v_i}([a_{\infty},c,b_2,a_3,b_{\infty}] \cup [a_0,b_1,b_0])$\\
& if $a_\y^{v_i}$ and $a_1^{v_i}$ are in the same component of $J_4\setminus\sigma_{v_i}(H')$.\\
3.6b& $J_2:\sigma_{v_i}([a_1,a_0,c,a_{\infty}] \cup [b_{\infty},b_2,b_3,b_1])$\\
& $J_3:\sigma_{v_i}([a_{\infty},a_2,a_3,b_0,b_{\infty}] \cup [a_1,b_1])$\\
& $J_4:\sigma_{v_i}([a_1,c,a_3,a_{\infty}] \cup [b_{\infty},a_2,b_0,b_1])$\\
& $J_5:\sigma_{v_i}([a_0,b_1,a_3,b_2,c,b_{\infty}] \cup [a_{\infty},b_0])$\\
& $J_6:\sigma_{v_i}([a_{\infty},b_2,a_1,a_3,b_{\infty}] \cup [a_0,b_0])$\\
& if $a_\y^{v_i}$ and $b_\y^{v_i}$ are in the same component of $J_4\setminus\sigma_{v_i}(H')$.\\
\hline
3.7& $J_1:\sigma_{v_i}([a_1,b_3,b_2,a_0,b_1] \cup [a_{\infty},b_{\infty}])$\\
& $J_2:\sigma_{v_i}([a_{\infty},b_2,c,b_{\infty}] \cup [a_1,a_0,b_3,b_1])$\\
& $J_3:\sigma_{v_i}([a_{\infty},a_2,a_3,b_0,b_{\infty}] \cup [a_1,b_1])$\\
& $J_4:\sigma_{v_i}([a_1,c,a_3,b_1] \cup [a_{\infty},b_0,a_2,b_{\infty}])$\\
& $J_5:\sigma_{v_i}([b_{\infty},a_3,b_2,b_1,b_0] \cup [a_0,c,a_{\infty}])$\\
\hline
3.8& $J_1:\sigma_{v_i}([a_1,a_0,b_3,b_2,b_1] \cup [a_{\infty},b_{\infty}])$\\
& $J_2:\sigma_{v_i}([a_{\infty},b_2,a_0,c,b_{\infty}] \cup [a_1,b_3,b_1])$\\
& $J_3:\sigma_{v_i}([a_{\infty},a_2,a_3,b_0,b_{\infty}] \cup [a_1,b_1])$\\
& $J_4:\sigma_{v_i}([a_1,c,a_3,b_1] \cup [a_{\infty},b_0,a_2,b_{\infty}])$\\
& $J_5:\sigma_{v_i}([a_{\infty},c,b_2,a_3,b_{\infty}] \cup [a_0,b_1,b_0])$\\
\hline
\end{tabular}
\end{center}

Finally suppose $v_i\in U'\cap V'$. In this case $a_\y^{v_i}$ and $b_\y^{v_i}$ are in the same component of $J_2$, $J_4$ and $J_6$,
and we have eight cases for the remaining $2$-factors. Namely, for each of $J_1$, $J_3$ and $J_5$, $a_\y^{v_i}$ and $b_\y^{v_i}$
are either in the same component or they are not. We number these eight cases as in the following table. 

\vspace{0.3cm}

\begin{center}
\begin{tabular}{c|c|c|c|}
&$J_1$&$J_3$&$J_5$\\
\hline
4.1& same&same&same\\
\hline
4.2& distinct&same&same\\
\hline
4.3& same&distinct&same\\
\hline
4.4& same&same&distinct\\
\hline
4.5& same&distinct&distinct\\
\hline
4.6& distinct&same&distinct\\
\hline
4.7& distinct&distinct&same\\
\hline
4.8& distinct&distinct&distinct\\
\hline
\end{tabular}
\end{center}

Depending on which of cases 4.1--4.8 that we are in, 
we can reallocate the edges of $\sigma_{v_i}(H)$ to the factors $J_1,J_2,\ldots,J_6$ 
as indicated in the following table to obtain a new 
$2$-factorisation of $L$ with the desired properties. If $J_x$ ($x\in\{1,2,\ldots,6\}$) is not listed for a particular case, then the edges of 
$\sigma_{v_i}(H')$ that are in $J_x$ are unchanged. 

\begin{center}
\begin{tabular}{|c|l|}
\hline
4.1& \\
\hline
4.2& $J_1:\sigma_{v_i}([a_{\infty},c,a_0,b_2,b_{\infty}] \cup [a_1,b_3,b_1])$\\
& $J_2:\sigma_{v_i}([a_1,a_0,b_3,b_2,b_1] \cup [a_{\infty},b_{\infty}])$\\
\hline
4.3& $J_3:\sigma_{v_i}([a_1,c,a_3,b_0,b_1] \cup [a_{\infty},a_2,b_{\infty}])$\\
& $J_4:\sigma_{v_i}([a_{\infty},a_3,a_2,b_0,b_{\infty}] \cup [a_1,b_1])$\\
\hline
4.4a& $J_3:\sigma_{v_i}([a_1,c,a_3,b_1] \cup [a_{\infty},b_0,a_2,b_{\infty}])$\\
& $J_5:\sigma_{v_i}([a_{\infty},a_3,b_2,c,b_{\infty}] \cup [a_0,b_1,b_0])$\\
& if $a_\y^{v_i}$ and $b_1^{v_i}$ are in the same component of $J_3\setminus\sigma_{v_i}(H')$.\\
4.4b& $J_1:\sigma_{v_i}([a_{\infty},b_2,a_0,c,b_{\infty}] \cup [a_1,b_3,b_1])$\\
& $J_2:\sigma_{v_i}([a_1,a_0,b_3,b_2,b_1] \cup [a_{\infty},b_{\infty}])$\\
& $J_3:\sigma_{v_i}([a_1,c,a_3,b_{\infty}] \cup [a_{\infty},a_2,b_0,b_1])$\\
& $J_4:\sigma_{v_i}([a_{\infty},b_0,a_3,a_2,b_{\infty}] \cup [a_1,b_1])$\\
& $J_5:\sigma_{v_i}([a_0,b_1,a_3,b_2,c,a_{\infty}] \cup [b_{\infty},b_0])$\\
& $J_6:\sigma_{v_i}([a_{\infty},a_3,a_1,b_2,b_{\infty}] \cup [a_0,b_0])$\\
& if $a_\y^{v_i}$ and $b_1^{v_i}$ are in the same component of $J_1\setminus\sigma_{v_i}(H')$ and \\
& $a_\y^{v_i}$ and $b_\y^{v_i}$ are in the same component of $J_3\setminus\sigma_{v_i}(H')$.\\
4.4c & $J_1:\sigma_{v_i}([a_1,b_3,a_0,c,b_{\infty}] \cup [a_{\infty},b_2,b_1])$\\
& $J_3:\sigma_{v_i}([a_1,c,a_3,b_{\infty}] \cup [a_{\infty},a_2,b_0,b_1])$\\
& $J_4:\sigma_{v_i}([a_{\infty},b_0,a_3,a_2,b_{\infty}] \cup [a_1,b_1])$\\
& $J_5:\sigma_{v_i}([a_0,b_1,a_3,b_2,c,a_{\infty}] \cup [b_{\infty},b_0])$\\
& $J_6:\sigma_{v_i}([a_{\infty},a_3,a_1,b_2,b_{\infty}] \cup [a_0,b_0])$\\
& if $a_\y^{v_i}$ and $b_\y^{v_i}$ are in the same component of $J_1\setminus\sigma_{v_i}(H')$ and \\
& $a_\y^{v_i}$ and $b_\y^{v_i}$ are in the same component of $J_3\setminus\sigma_{v_i}(H')$.\\
\hline
4.5& $J_3:\sigma_{v_i}([a_1,c,a_3,b_1] \cup [a_{\infty},b_0,a_2,b_{\infty}])$\\
& $J_5:\sigma_{v_i}([a_{\infty},a_3,b_2,c,b_{\infty}] \cup [a_0,b_1,b_0])$\\
\hline
\end{tabular}
\end{center}

\begin{center}
\begin{tabular}{|c|l|}
\hline
4.6a& $J_1:\sigma_{v_i}([a_{\infty},b_2,a_0,c,b_{\infty}] \cup [a_1,b_3,b_1])$\\
& $J_2:\sigma_{v_i}([a_1,a_0,b_3,b_2,b_1] \cup [a_{\infty},b_{\infty}])$\\
& $J_3:\sigma_{v_i}([a_1,c,a_3,b_1] \cup [a_{\infty},b_0,a_2,b_{\infty}])$\\
& $J_5:\sigma_{v_i}([a_{\infty},c,b_2,a_3,b_{\infty}] \cup [a_0,b_1,b_0])$\\
& $J_6:\sigma_{v_i}([a_{\infty},a_3,a_1,b_2,b_{\infty}] \cup [a_0,b_0])$\\
& if $a_\y^{v_i}$ and $b_1^{v_i}$ are in the same component of $J_3\setminus\sigma_{v_i}(H')$ \\
4.6b & $J_1:\sigma_{v_i}([a_1,b_3,a_0,c,b_{\infty}] \cup [a_{\infty},b_2,b_1])$\\
& $J_3:\sigma_{v_i}([a_1,c,a_3,b_{\infty}] \cup [a_{\infty},a_2,b_0,b_1])$\\
& $J_4:\sigma_{v_i}([a_{\infty},b_0,a_3,a_2,b_{\infty}] \cup [a_1,b_1])$\\
& $J_5:\sigma_{v_i}([a_0,b_1,a_3,b_2,c,a_{\infty}] \cup [b_{\infty},b_0])$\\
& $J_6:\sigma_{v_i}([a_{\infty},a_3,a_1,b_2,b_{\infty}] \cup [a_0,b_0])$\\
& if $a_\y^{v_i}$ and $b_\y^{v_i}$ are in the same component of $J_3\setminus\sigma_{v_i}(H')$.\\
\hline
4.7& $J_1:\sigma_{v_i}([a_{\infty},c,a_0,b_2,b_{\infty}] \cup [a_1,b_3,b_1])$\\
& $J_2:\sigma_{v_i}([a_1,a_0,b_3,b_2,b_1] \cup [a_{\infty},b_{\infty}])$\\
& $J_3:\sigma_{v_i}([a_1,c,a_3,b_0,b_1] \cup [a_{\infty},a_2,b_{\infty}])$\\
& $J_4:\sigma_{v_i}([a_{\infty},a_3,a_2,b_0,b_{\infty}] \cup [a_1,b_1])$\\
\hline
4.8& $J_1:\sigma_{v_i}([a_{\infty},b_2,a_0,c,b_{\infty}] \cup [a_1,b_3,b_1])$\\
& $J_2:\sigma_{v_i}([a_1,a_0,b_3,b_2,b_1] \cup [a_{\infty},b_{\infty}])$\\
& $J_3:\sigma_{v_i}([a_1,c,a_3,b_1] \cup [a_{\infty},b_0,a_2,b_{\infty}])$\\
& $J_5:\sigma_{v_i}([a_{\infty},c,b_2,a_3,b_{\infty}] \cup [a_0,b_1,b_0])$\\
& $J_6:\sigma_{v_i}([a_{\infty},a_3,a_1,b_2,b_{\infty}] \cup [a_0,b_0])$\\
\hline
\end{tabular}
\end{center}

\noindent{\bf The case n=6:}
For $i=1,2,3$, let $X_1(i)$, $Y_1(i)$, $X'_1(i)$, $Y'_1(i)$, $X_2(i)$, $Y_2(i)$, $X'_2(i)$, $Y'_2(i)$ be the subgraphs of $K_{A_n\cup B_n\cup\{c\}}$ given by the union of
the paths listed in the following tables, and let $H^i$ be the subgraph of $K_{A_n\cup B_n\cup\{c\}}$ with edge set
$E(X_1(i))\cup E(Y_1(i))$. Applying Lemma \ref{4RegularHamiltonFragments} with $H=H^i$, $X_1=X_1(i)$, $Y_1=Y_1(i)$, $X'_1=X'_1(i)$, $Y'_1=Y_1(i)$,$X_2=X_2(i)$, $Y_2=Y_2(i)$, $X'_2=X_2(i)$, $Y'_2=Y_2(i)$  shows that each 
$H^i$ is a Hamilton fragment. The value of $t$ can be deduced from the ends of the given paths.

\vspace{0.3cm}

$\begin{array}{|c|c|}
\hline
X_1(1)& [a_{\infty},a_1,a_2,b_0,c,b_3,b_{\infty}] \cup [a_4,b_4]\\
\hline
Y_1(1)&[a_{\infty},a_2,b_3,a_1,b_4] \cup [a_4,b_0,b_{\infty}]\\
\hline
X'_1(1)& Y_1(1) \\
\hline
Y'_1(1) & X_1(1)\\
\hline
X_2(1)& [a_{\infty},a_1,a_2,b_3,c,b_0,b_{\infty}] \cup [a_4,b_4]\\
\hline
Y_2(1)& [b_{\infty},b_3,a_1,b_4] \cup [a_4,b_0,a_2,a_{\infty}]\\
\hline
X'_2(1) & Y_2(1) \\
\hline
Y'_2(1) & X_2(1) \\
\hline
\end{array}
$
$\begin{array}{|c|c|}
\hline
X_1(2)& [a_{\infty},a_4,b_2,b_1,c,a_0,b_{\infty}] \cup [a_3,b_3]\\
\hline
Y_1(2)&[a_3,a_4,a_0,b_1,b_{\infty}] \cup [a_{\infty},b_2,b_3]\\
\hline
X'_1(2)& Y_1(2) \\
\hline
Y'_1(2) & X_1(2)\\
\hline
X_2(2)& [a_{\infty},b_2,a_4,a_0,c,b_1,b_{\infty}] \cup [a_3,b_3]\\
\hline
Y_2(2)& [b_{\infty},a_0,b_1,b_2,b_3] \cup [a_3,a_4,a_{\infty}]\\
\hline
X'_2(2) & Y_2(2) \\
\hline
Y'_2(2) & X_2(2) \\
\hline
\end{array}
$

$\begin{array}{|c|c|}
\hline
X_1(3)& [a_{\infty},b_4,b_0,a_1,a_3,b_{\infty}] \cup [a_2,b_2]\\
\hline
Y_1(3)&[a_2,b_4,c,a_1,b_{\infty}] \cup [a_{\infty},b_0,a_3,b_2]\\
\hline
X'_1(3)& Y_1(3) \\
\hline
Y'_1(3) & X_1(3)\\
\hline
X_2(3)& [a_{\infty},b_4,b_0,a_3,a_1,b_{\infty}] \cup [a_2,b_2]\\
\hline
Y_2(3)& [a_2,b_4,c,a_1,b_0,a_{\infty}] \cup [b_{\infty},a_3,b_2]\\
\hline
X'_2(3) & Y_2(3) \\
\hline
Y'_2(3) & X_2(3) \\
\hline
\end{array}
$

Note that $H^1$, $H^2$ and $H^3$ are pairwise edge disjoint and let $H'=K_{A_5\cup B_5\cup \{c\}}-{H^1\cup H^2\cup H^3}$ so that $\{H', H^1, H^2, H^3\}$ is a decomposition of $K_{A_6\cup B_6\cup \{c\}}\cong K_{13}$. We now show that $H'$ is a Hamilton fragment. 

Our aim is to show that the subgraph $$L=\bigcup_{v\in V(G)}\sigma_v(H')$$
of $L(G)$ decomposes into Hamilton cycles, and we do this by applying Lemma \ref{manyrepairs} 
with $H_i=\sigma_{v_i}(H')$ and $V_i=\{a_\infty^{v_i},b_\infty^{v_i}\}$ for each $i=1,2,\ldots,m$.
To this end, observe that \linebreak $\sigma_{v_1}(H'),\sigma_{v_2}(H'),\ldots,\sigma_{v_m}(H')$ are edge-disjoint subgraphs of $L$,
and since $b_\infty^{v_i}=a_\infty^{v_{i+1}}$ for $i=1,2,\ldots,m-1$ it follows that $V_i\cap V_{i+1}\ne \emptyset$ for $i=1,2,\ldots,m-1$ as required. It remains to show there is a $2$-factorisation
$\{J_1,J_2,J_3,J_4,J_5,J_6\}$ of $L$ such that 
\begin{itemize}
\item $\{a_\infty^{v_1},a_\infty^{v_2},\ldots,a_\infty^{v_m}\}$ links $\{J_1,J_2,J_3,J_4,J_5,J_6\}$; and
\item $\sigma_{v_i}(H')$ induces an $\{a_\infty^{v_{i}},b_\infty^{v_i}\}$-connector in $\{J_1,J_2,J_3,J_4,J_5,J_6\}$, for $i=1,2,\ldots,m$.
\end{itemize} 
Let $\{U,U'\}$ be a partition of $V(G)$ such that both $U$ and $U'$ link $\{F,\F_1\}$,
let $\{V,V'\}$ be a partition of $V(G)$ such that both $V$ and $V'$ link $\{F,\F_2\}$ 
(such partitions exist by Lemma \ref{intersects}),
and let $\{J_1,J_2,J_3,J_4,J_5,J_6\}$ be the decomposition of $L$ defined by 

\begin{itemize}
\item 
$
\begin{array}[t]{lll}
J_1&=&\bigcup_{v\in U\cap V}\sigma_v([a_1,a_4,c,b_2,b_0,b_3,b_1] \cup [a_{\infty},b_{\infty}])\ \cup \\
&&\bigcup_{v\in U'\cap V}\sigma_v([a_1,a_4,b_3,b_0,b_2,b_{\infty}] \cup [a_{\infty},b_1])\ \cup \\
&&\bigcup_{v\in U\cap V'}\sigma_v([a_1,a_4,c,b_2,b_0,b_3,b_1] \cup [a_{\infty},b_{\infty}])\ \cup \\
&&\bigcup_{v\in U'\cap V'}\sigma_v([b_{\infty},b_2,b_0,b_1] \cup [a_1,a_4,b_3,a_{\infty}]);
\end{array}
$
\item 
$
\begin{array}[t]{lll}
J_2&=&\bigcup_{v\in U\cap V}\sigma_v([a_1,a_0,b_4,b_2,b_{\infty}] \cup [a_{\infty},c,a_3,b_1])\ \cup \\
&&\bigcup_{v\in U'\cap V}\sigma_v([a_1,b_2,a_0,a_3,b_4,b_1] \cup [a_{\infty},b_{\infty}])\ \cup \\
&&\bigcup_{v\in U\cap V'}\sigma_v([a_1,a_0,b_4,b_2,b_{\infty}] \cup [a_{\infty},c,a_3,b_1])\ \cup \\
&&\bigcup_{v\in U'\cap V'}\sigma_v([a_1,b_2,a_0,b_4,a_3,b_1] \cup [a_{\infty},b_{\infty}]);
\end{array}
$
\item 
$
\begin{array}[t]{lll}
J_3&=&\bigcup_{v\in U\cap V}\sigma_v([a_{\infty},a_3,a_0,a_2,a_4,b_{\infty}] \cup [a_1,b_1])\ \cup \\
&&\bigcup_{v\in U'\cap V}\sigma_v([a_1,a_0,a_2,c,a_4,b_{\infty}] \cup [a_{\infty},a_3,b_1])\ \cup \\
&&\bigcup_{v\in U\cap V'}\sigma_v([a_{\infty},a_0,a_3,a_2,a_4,b_{\infty}] \cup [a_1,b_1])\ \cup \\
&&\bigcup_{v\in U'\cap V'}\sigma_v([a_1,a_0,a_3,c,a_{\infty}] \cup [b_{\infty},a_2,a_4,b_1]);
\end{array}
$
\item 
$
\begin{array}[t]{lll}
J_4&=&\bigcup_{v\in U\cap V}\sigma_v([a_1,b_2,a_0,b_3,b_4,b_{\infty}] \cup [a_{\infty},b_1])\ \cup \\
&&\bigcup_{v\in U'\cap V}\sigma_v([a_{\infty},a_0,b_3,b_4,b_2,c,b_{\infty}] \cup [a_1,b_1])\ \cup \\
&&\bigcup_{v\in U\cap V'}\sigma_v([a_1,b_2,a_0,b_3,b_4,b_{\infty}] \cup [a_{\infty},b_1])\ \cup \\
&&\bigcup_{v\in U'\cap V'}\sigma_v([a_{\infty},a_0,b_3,b_4,b_2,c,b_{\infty}] \cup [a_1,b_1]);
\end{array}
$
\item 
$
\begin{array}[t]{lll}
J_5&=&\bigcup_{v\in U\cap V}\sigma_v([a_{\infty},b_3,a_4,b_1,a_2,b_{\infty}] \cup [a_0,b_0])\ \cup \\
&&\bigcup_{v\in U'\cap V}\sigma_v([a_{\infty},b_3,b_1,a_4,a_2,b_{\infty}] \cup [a_0,b_0])\ \cup \\
&&\bigcup_{v\in U\cap V'}\sigma_v([a_{\infty},b_3,a_4,b_1,b_0] \cup [a_0,a_2,c,b_{\infty}])\ \cup \\
&&\bigcup_{v\in U'\cap V'}\sigma_v([a_0,a_2,c,a_4,b_{\infty}] \cup [a_{\infty},b_1,b_3,b_0]);
\end{array}
$
\item 
$
\begin{array}[t]{lll}
J_6&=&\bigcup_{v\in U\cap V}\sigma_v([b_{\infty},c,a_2,a_3,b_4,b_1,b_0] \cup [a_0,a_{\infty}])\ \cup \\
&&\bigcup_{v\in U'\cap V}\sigma_v([a_{\infty},c,a_3,a_2,b_1,b_0] \cup [a_0,b_4,b_{\infty}])\ \cup \\
&&\bigcup_{v\in U\cap V'}\sigma_v([a_{\infty},a_3,b_4,b_1,a_2,b_{\infty}] \cup [a_0,b_0])\ \cup \\
&&\bigcup_{v\in U'\cap V'}\sigma_v([a_{\infty},a_3,a_2,b_1,b_4,b_{\infty}] \cup [a_0,b_0]);
\end{array}
$

\end{itemize}
It is easily checked that $\{J_1,J_2,J_3,J_4,J_5,J_6\}$ is a $2$-factorisation of $L$ and that $\{a_\infty^{v_1},a_\infty^{v_2},\ldots,a_\infty^{v_m}\}$ links $\{J_1,J_2,J_3,J_4,J_5,J_6\}$.
It remains to show that for $i=1,2,\ldots,m$, $\sigma_{v_i}(H')$ induces an $\{a_\infty^{v_{i}},b_\infty^{v_i}\}$-connector in $\{J_1,J_2,J_3,J_4,J_5,J_6\}$. 
There are four cases to consider: $v_i\in U\cap V$, $v_i\in U'\cap V$, $v_i\in U\cap V'$ and $v_i\in U'\cap V'$.

First suppose $v_i\in U\cap V$. In this case $a_\y^{v_i}$ and $b_\y^{v_i}$ are in the same component of $J_1$, $J_3$ and $J_5$,
and we have eight cases for the remaining $2$-factors. Namely, for each of $J_2$, $J_4$ and $J_6$, $a_\y^{v_i}$ and $b_\y^{v_i}$
are either in the same component or they are not. We number these eight cases as in the following table. 

\vspace{0.3cm}

\begin{center}
\begin{tabular}{c|c|c|c|}
&$J_2$&$J_4$&$J_6$\\
\hline
1.1& same&same&same\\
\hline
1.2& distinct&same&same\\
\hline
1.3& same&distinct&same\\
\hline
1.4& same&same&distinct\\
\hline
1.5& same&distinct&distinct\\
\hline
1.6& distinct&same&distinct\\
\hline
1.7& distinct&distinct&same\\
\hline
1.8& distinct&distinct&distinct\\
\hline
\end{tabular}
\end{center}

Depending on which of cases 1.1--1.8 that we are in, 
we can reallocate the edges of $\sigma_{v_i}(H)$ to the factors $J_1,J_2,\ldots,J_6$ 
as indicated in the following table to obtain a new 
$2$-factorisation of $L$ with the desired properties. If $J_x$ ($x\in\{1,2,\ldots,6\}$) is not listed for a particular case, then the edges of 
$\sigma_{v_i}(H)$ that are in $J_x$ are unchanged. 

\begin{center}
\begin{tabular}{|c|l|}
\hline
1.1& \\
\hline
1.2& $J_2:\sigma_{v_i}([a_1,b_2,a_0,b_4,b_1] \cup [a_{\infty},a_3,c,b_{\infty}])$\\
& $J_3:\sigma_{v_i}([a_{\infty},a_0,a_3,a_2,a_4,b_{\infty}] \cup [a_1,b_1])$\\
& $J_4:\sigma_{v_i}([a_1,a_0,b_3,b_4,b_2,b_{\infty}] \cup [a_{\infty},b_1])$\\
& $J_6:\sigma_{v_i}([b_{\infty},b_4,a_3,b_1,b_0] \cup [a_0,a_2,c,a_{\infty}])$\\
\hline
1.3a& $J_1:\sigma_{v_i}([a_{\infty},c,a_4,b_3,b_0,b_2,b_{\infty}] \cup [a_1,b_1])$\\
& $J_2:\sigma_{v_i}([a_1,a_0,b_4,b_2,c,b_{\infty}] \cup [a_{\infty},a_3,b_1])$\\
& $J_3:\sigma_{v_i}([a_{\infty},a_0,a_3,a_2,b_{\infty}] \cup [a_1,a_4,b_1])$\\
& $J_4:\sigma_{v_i}([a_1,b_2,a_0,b_3,b_4,b_1] \cup [a_{\infty},b_{\infty}])$\\
& $J_5:\sigma_{v_i}([a_{\infty},b_3,b_1,a_2,a_4,b_{\infty}] \cup [a_0,b_0])$\\
& $J_6:\sigma_{v_i}([a_0,a_2,c,a_3,b_4,b_{\infty}] \cup [a_{\infty},b_1,b_0])$\\
& if $a_\y^{v_i}$ and $b_\y^{v_i}$ are in the same component of $J_6\setminus\sigma_{v_i}(H')$.\\
1.3b& $J_1:\sigma_{v_i}([a_{\infty},b_3,b_0,b_2,c,a_4,b_{\infty}] \cup [a_1,b_1])$\\
& $J_3:\sigma_{v_i}([a_{\infty},a_3,a_0,a_2,b_{\infty}] \cup [a_1,a_4,b_1])$\\
& $J_4:\sigma_{v_i}([a_1,b_2,a_0,a_{\infty}] \cup [b_{\infty},b_4,b_3,b_1])$\\
& $J_5:\sigma_{v_i}([a_0,b_3,a_4,a_2,b_1,b_0] \cup [a_{\infty},b_{\infty}])$\\
& $J_6:\sigma_{v_i}([a_{\infty},b_1,b_4,a_3,a_2,c,b_{\infty}] \cup [a_0,b_0])$\\
& if $a_\y^{v_i}$ and $b_0^{v_i}$ are in the same component of $J_6\setminus\sigma_{v_i}(H')$.\\
\hline
1.4& $J_3:\sigma_{v_i}([a_{\infty},a_0,a_3,a_2,a_4,b_{\infty}] \cup [a_1,b_1])$\\
& $J_6:\sigma_{v_i}([a_{\infty},a_3,b_4,b_1,b_0] \cup [a_0,a_2,c,b_{\infty}])$\\
\hline
1.5& $J_1:\sigma_{v_i}([a_{\infty},b_3,b_0,b_2,c,a_4,b_{\infty}] \cup [a_1,b_1])$\\
& $J_3:\sigma_{v_i}([a_{\infty},a_3,a_0,a_2,b_{\infty}] \cup [a_1,a_4,b_1])$\\
& $J_4:\sigma_{v_i}([a_1,b_2,a_0,a_{\infty}] \cup [b_{\infty},b_4,b_3,b_1])$\\
& $J_5:\sigma_{v_i}([a_0,b_3,a_4,a_2,b_1,b_0] \cup [a_{\infty},b_{\infty}])$\\
& $J_6:\sigma_{v_i}([a_{\infty},b_1,b_4,a_3,a_2,c,b_{\infty}] \cup [a_0,b_0])$\\
\hline
1.6& $J_2:\sigma_{v_i}([b_{\infty},c,a_3,b_4,b_1] \cup [a_1,b_2,a_0,a_{\infty}])$\\
& $J_4:\sigma_{v_i}([a_1,a_0,b_3,b_4,b_2,b_{\infty}] \cup [a_{\infty},b_1])$\\
& $J_6:\sigma_{v_i}([a_{\infty},c,a_2,a_3,b_1,b_0] \cup [a_0,b_4,b_{\infty}])$\\
\hline
\end{tabular}
\end{center}

\begin{center}
\begin{tabular}{|c|l|}
\hline
1.7a& $J_1:\sigma_{v_i}([a_1,b_2,b_0,b_3,b_1] \cup [a_{\infty},c,a_4,b_{\infty}])$\\
& $J_2:\sigma_{v_i}([a_1,a_0,b_2,c,a_3,a_{\infty}] \cup [b_{\infty},b_4,b_1])$\\
& $J_3:\sigma_{v_i}([a_1,a_4,a_2,a_0,a_3,b_1] \cup [a_{\infty},b_{\infty}])$\\
& $J_4:\sigma_{v_i}([a_{\infty},a_0,b_3,b_4,b_2,b_{\infty}] \cup [a_1,b_1])$\\
& $J_6:\sigma_{v_i}([a_0,b_4,a_3,a_2,c,b_{\infty}] \cup [a_{\infty},b_1,b_0])$\\
& if $a_\y^{v_i}$ and $b_\y^{v_i}$ are in the same component of $J_6\setminus\sigma_{v_i}(H')$.\\
1.7b& $J_1:\sigma_{v_i}([a_{\infty},b_3,a_4,c,b_{\infty}] \cup [a_1,b_2,b_0,b_1])$\\
& $J_2:\sigma_{v_i}([a_{\infty},a_0,b_4,a_3,c,b_2,b_{\infty}] \cup [a_1,b_1])$\\
& $J_3:\sigma_{v_i}([a_{\infty},a_3,a_0,a_2,b_{\infty}] \cup [a_1,a_4,b_1])$\\
& $J_4:\sigma_{v_i}([a_1,a_0,b_2,b_4,b_3,b_1] \cup [a_{\infty},b_{\infty}])$\\
& $J_5:\sigma_{v_i}([a_{\infty},b_1,a_2,a_4,b_{\infty}] \cup [a_0,b_3,b_0])$\\
& $J_6:\sigma_{v_i}([a_{\infty},c,a_2,a_3,b_1,b_4,b_{\infty}] \cup [a_0,b_0])$\\
& if $a_\y^{v_i}$ and $b_0^{v_i}$ are in the same component of $J_6\setminus\sigma_{v_i}(H')$.\\
\hline
1.8& $J_1:\sigma_{v_i}([a_{\infty},b_3,a_4,c,b_{\infty}] \cup [a_1,b_2,b_0,b_1])$\\
& $J_2:\sigma_{v_i}([a_{\infty},a_0,b_4,a_3,c,b_2,b_{\infty}] \cup [a_1,b_1])$\\
& $J_3:\sigma_{v_i}([a_{\infty},a_3,a_0,a_2,b_{\infty}] \cup [a_1,a_4,b_1])$\\
& $J_4:\sigma_{v_i}([a_1,a_0,b_2,b_4,b_3,b_1] \cup [a_{\infty},b_{\infty}])$\\
& $J_5:\sigma_{v_i}([a_{\infty},b_1,a_2,a_4,b_{\infty}] \cup [a_0,b_3,b_0])$\\
& $J_6:\sigma_{v_i}([a_{\infty},c,a_2,a_3,b_1,b_4,b_{\infty}] \cup [a_0,b_0])$\\
\hline
\end{tabular}
\end{center}

Next suppose $v_i\in U'\cap V$. In this case $a_\y^{v_i}$ and $b_\y^{v_i}$ are in the same component of $J_2$, $J_4$ and $J_5$,
and we have eight cases for the remaining $2$-factors. Namely, for each of $J_1$, $J_3$ and $J_6$, $a_\y^{v_i}$ and $b_\y^{v_i}$
are either in the same component or they are not. We number these eight cases as in the following table. 

\vspace{0.3cm}

\begin{center}
\begin{tabular}{c|c|c|c|}
&$J_1$&$J_3$&$J_6$\\
\hline
2.1& same&same&same\\
\hline
2.2& distinct&same&same\\
\hline
2.3& same&distinct&same\\
\hline
2.4& same&same&distinct\\
\hline
2.5& same&distinct&distinct\\
\hline
2.6& distinct&same&distinct\\
\hline
2.7& distinct&distinct&same\\
\hline
2.8& distinct&distinct&distinct\\
\hline
\end{tabular}
\end{center}

Depending on which of cases 2.1--2.8 that we are in, 
we can reallocate the edges of $\sigma_{v_i}(H)$ to the factors $J_1,J_2,\ldots,J_6$ 
as indicated in the following table to obtain a new 
$2$-factorisation of $L$ with the desired properties. If $J_x$ ($x\in\{1,2,\ldots,6\}$) is not listed for a particular case, then the edges of 
$\sigma_{v_i}(H)$ that are in $J_x$ are unchanged. 

\begin{center}
\begin{tabular}{|c|l|}
\hline
2.1& \\
\hline
2.2& $J_1:\sigma_{v_i}([a_{\infty},b_3,b_0,b_2,b_{\infty}] \cup [a_1,a_4,b_1])$\\
& $J_5:\sigma_{v_i}([a_{\infty},b_1,b_3,a_4,a_2,b_{\infty}] \cup [a_0,b_0])$
\\
\hline
2.3 & $J_3:\sigma_{v_i}([a_{\infty},a_3,c,a_4,b_{\infty}] \cup [a_1,a_0,a_2,b_1])$\\
& $J_6:\sigma_{v_i}([a_{\infty},c,a_2,a_3,b_1,b_0] \cup [a_0,b_4,b_{\infty}])$
\\
\hline
2.4 & $J_2:\sigma_{v_i}([a_1,b_2,a_0,b_4,a_3,b_1] \cup [a_{\infty},b_{\infty}])$\\
& $J_3:\sigma_{v_i}([a_{\infty},a_3,c,a_4,b_1] \cup [a_1,a_0,a_2,b_{\infty}])$\\
& $J_5:\sigma_{v_i}([a_{\infty},b_3,b_1,a_2,a_4,b_{\infty}] \cup [a_0,b_0])$\\
& $J_6:\sigma_{v_i}([a_0,a_3,a_2,c,a_{\infty}] \cup [b_{\infty},b_4,b_1,b_0])$
\\
\hline
2.5 & $J_2:\sigma_{v_i}([a_1,b_2,a_0,b_4,a_3,b_1] \cup [a_{\infty},b_{\infty}])$\\
& $J_3:\sigma_{v_i}([a_1,a_0,a_3,a_2,b_1] \cup [a_{\infty},c,a_4,b_{\infty}])$\\
& $J_6:\sigma_{v_i}([a_0,a_2,c,a_3,a_{\infty}] \cup [b_{\infty},b_4,b_1,b_0])$
\\
\hline
2.6 & $J_1:\sigma_{v_i}([a_1,b_2,b_0,b_3,a_{\infty}] \cup [b_{\infty},a_4,b_1])$\\
& $J_2:\sigma_{v_i}([a_{\infty},a_3,a_0,b_2,b_4,b_{\infty}] \cup [a_1,b_1])$\\
& $J_3:\sigma_{v_i}([a_{\infty},a_0,a_2,a_3,b_1] \cup [a_1,a_4,c,b_{\infty}])$\\
& $J_4:\sigma_{v_i}([a_1,a_0,b_3,b_4,b_1] \cup [a_{\infty},c,b_2,b_{\infty}])$\\
& $J_5:\sigma_{v_i}([a_{\infty},b_1,b_3,a_4,a_2,b_{\infty}] \cup [a_0,b_0])$\\
& $J_6:\sigma_{v_i}([a_0,b_4,a_3,c,a_2,b_1,b_0] \cup [a_{\infty},b_{\infty}])$
\\
\hline
2.7 & $J_1:\sigma_{v_i}([a_{\infty},b_3,b_0,b_2,b_{\infty}] \cup [a_1,a_4,b_1])$\\
& $J_3:\sigma_{v_i}([a_{\infty},a_3,c,a_4,b_{\infty}] \cup [a_1,a_0,a_2,b_1])$\\
& $J_5:\sigma_{v_i}([a_{\infty},b_1,b_3,a_4,a_2,b_{\infty}] \cup [a_0,b_0])$\\
& $J_6:\sigma_{v_i}([a_{\infty},c,a_2,a_3,b_1,b_0] \cup [a_0,b_4,b_{\infty}])$
\\
\hline
2.8 & $J_1:\sigma_{v_i}([a_{\infty},b_3,b_0,b_2,b_{\infty}] \cup [a_1,a_4,b_1])$\\
& $J_2:\sigma_{v_i}([a_1,b_2,a_0,b_4,a_3,b_1] \cup [a_{\infty},b_{\infty}])$\\
& $J_3:\sigma_{v_i}([a_1,a_0,a_3,a_2,b_1] \cup [a_{\infty},c,a_4,b_{\infty}])$\\
& $J_5:\sigma_{v_i}([a_{\infty},b_1,b_3,a_4,a_2,b_{\infty}] \cup [a_0,b_0])$\\
& $J_6:\sigma_{v_i}([a_0,a_2,c,a_3,a_{\infty}] \cup [b_{\infty},b_4,b_1,b_0])$
\\
\hline
\end{tabular}
\end{center}

Next suppose $v_i\in U\cap V'$. In this case $a_\y^{v_i}$ and $b_\y^{v_i}$ are in the same component of $J_1$, $J_3$ and $J_6$,
and we have eight cases for the remaining $2$-factors. Namely, for each of $J_2$, $J_4$ and $J_5$, $a_\y^{v_i}$ and $b_\y^{v_i}$
are either in the same component or they are not. We number these eight cases as in the following table. 

\vspace{0.3cm}

\begin{center}
\begin{tabular}{c|c|c|c|}
&$J_2$&$J_4$&$J_5$\\
\hline
3.1& same&same&same\\
\hline
3.2& distinct&same&same\\
\hline
3.3& same&distinct&same\\
\hline
3.4& same&same&distinct\\
\hline
3.5& same&distinct&distinct\\
\hline
3.6& distinct&same&distinct\\
\hline
3.7& distinct&distinct&same\\
\hline
3.8& distinct&distinct&distinct\\
\hline
\end{tabular}
\end{center}

Depending on which of cases 3.1--3.8 that we are in, 
we can reallocate the edges of $\sigma_{v_i}(H)$ to the factors $J_1,J_2,\ldots,J_6$ 
as indicated in the following table to obtain a new 
$2$-factorisation of $L$ with the desired properties. If $J_x$ ($x\in\{1,2,\ldots,6\}$) is not listed for a particular case, then the edges of 
$\sigma_{v_i}(H)$ that are in $J_x$ are unchanged. 

\begin{center}
\begin{tabular}{|c|l|}
\hline
3.1& \\
\hline
3.2a& $J_1:\sigma_{v_i}([a_1,b_2,b_0,b_3,a_4,b_1] \cup [a_{\infty},c,b_{\infty}])$\\
& $J_2:\sigma_{v_i}([a_{\infty},a_0,b_2,c,a_3,b_4,b_{\infty}] \cup [a_1,b_1])$ \\
& $J_3:\sigma_{v_i}([a_1,a_4,a_2,a_0,a_3,b_1] \cup [a_{\infty},b_{\infty}])$ \\
& $J_4:\sigma_{v_i}([a_1,a_0,b_3,b_4,b_2,b_{\infty}] \cup [a_{\infty},b_1])$ \\
& $J_5:\sigma_{v_i}([a_{\infty},b_3,b_1,a_2,c,a_4,b_{\infty}] \cup [a_0,b_0])$ \\
& $J_6:\sigma_{v_i}([a_0,b_4,b_1,b_0] \cup [a_{\infty},a_3,a_2,b_{\infty}])$ \\
& if $a_\y^{v_i}$ and $a_0^{v_i}$ are in the same component of $J_5\setminus\sigma_{v_i}(H')$.\\
3.2b & $J_1:\sigma_{v_i}([a_{\infty},c,b_2,b_0,b_3,a_4,b_{\infty}] \cup [a_1,b_1])$\\
& $J_2:\sigma_{v_i}([a_1,a_0,b_2,b_4,b_1] \cup [a_{\infty},a_3,c,b_{\infty}])$ \\
& $J_3:\sigma_{v_i}([a_1,a_4,a_2,a_0,a_3,b_1] \cup [a_{\infty},b_{\infty}])$ \\
& $J_4:\sigma_{v_i}([a_{\infty},a_0,b_4,b_3,b_1] \cup [a_1,b_2,b_{\infty}])$ \\
& $J_5:\sigma_{v_i}([b_{\infty},a_2,c,a_4,b_1,b_0] \cup [a_0,b_3,a_{\infty}])$ \\
& $J_6:\sigma_{v_i}([a_{\infty},b_1,a_2,a_3,b_4,b_{\infty}] \cup [a_0,b_0])$ \\
& if $a_\y^{v_i}$ and $b_\y^{v_i}$ are in the same component of $J_5\setminus\sigma_{v_i}(H')$.\\
\hline
3.3a& $J_1:\sigma_{v_i}([a_{\infty},b_3,b_0,b_2,c,b_{\infty}] \cup [a_1,a_4,b_1])$\\
& $J_3:\sigma_{v_i}([a_{\infty},a_3,a_0,a_2,a_4,b_{\infty}] \cup [a_1,b_1])$ \\
& $J_4:\sigma_{v_i}([a_1,b_2,a_0,b_3,b_4,b_1] \cup [a_{\infty},b_{\infty}])$ \\
& $J_5:\sigma_{v_i}([b_{\infty},a_2,c,a_4,b_3,b_1,b_0] \cup [a_0,a_{\infty}])$ \\
& $J_6:\sigma_{v_i}([a_{\infty},b_1,a_2,a_3,b_4,b_{\infty}] \cup [a_0,b_0])$ \\
& if $a_\y^{v_i}$ and $b_\y^{v_i}$ are in the same component of $J_5\setminus\sigma_{v_i}(H')$.\\
3.3b& $J_1:\sigma_{v_i}([a_{\infty},b_3,b_0,b_2,c,b_{\infty}] \cup [a_1,a_4,b_1])$\\
& $J_4:\sigma_{v_i}([a_1,b_2,a_0,b_3,b_4,b_1] \cup [a_{\infty},b_{\infty}])$ \\
& $J_5:\sigma_{v_i}([a_{\infty},b_1,b_3,a_4,c,a_2,b_{\infty}] \cup [a_0,b_0])$ \\
& $J_6:\sigma_{v_i}([a_{\infty},a_3,b_4,b_{\infty}] \cup [a_0,a_2,b_1,b_0])$ \\
& if $a_\y^{v_i}$ and $a_0^{v_i}$ are in the same component of $J_5\setminus\sigma_{v_i}(H')$.\\
\hline
3.4& $J_1:\sigma_{v_i}([a_{\infty},b_3,b_0,b_2,c,b_{\infty}] \cup [a_1,a_4,b_1])$\\
& $J_5:\sigma_{v_i}([a_0,a_2,c,a_4,b_3,b_1,b_0] \cup [a_{\infty},b_{\infty}])$ \\
\hline
3.5& $J_1:\sigma_{v_i}([a_{\infty},b_3,b_0,b_2,c,b_{\infty}] \cup [a_1,a_4,b_1])$\\
& $J_4:\sigma_{v_i}([a_1,b_2,a_0,b_3,b_4,b_1] \cup [a_{\infty},b_{\infty}])$ \\
& $J_5:\sigma_{v_i}([a_{\infty},b_1,b_3,a_4,c,a_2,b_{\infty}] \cup [a_0,b_0])$ \\
& $J_6:\sigma_{v_i}([a_{\infty},a_3,b_4,b_{\infty}] \cup [a_0,a_2,b_1,b_0])$ \\
\hline
\end{tabular}
\end{center}

\begin{center}
\begin{tabular}{|c|l|}
\hline
3.6& $J_1:\sigma_{v_i}([a_1,b_2,b_0,b_3,a_4,b_1] \cup [a_{\infty},c,b_{\infty}])$\\
& $J_2:\sigma_{v_i}([a_{\infty},a_0,b_2,c,a_3,b_4,b_{\infty}] \cup [a_1,b_1])$ \\
& $J_3:\sigma_{v_i}([a_1,a_4,a_2,a_0,a_3,b_1] \cup [a_{\infty},b_{\infty}])$ \\
& $J_4:\sigma_{v_i}([a_1,a_0,b_3,b_4,b_2,b_{\infty}] \cup [a_{\infty},b_1])$ \\
& $J_5:\sigma_{v_i}([a_{\infty},b_3,b_1,a_2,c,a_4,b_{\infty}] \cup [a_0,b_0])$ \\
& $J_6:\sigma_{v_i}([a_0,b_4,b_1,b_0] \cup [a_{\infty},a_3,a_2,b_{\infty}])$ \\
\hline
3.7a& $J_1:\sigma_{v_i}([a_{\infty},b_3,a_4,c,b_{\infty}] \cup [a_1,b_2,b_0,b_1])$\\
& $J_2:\sigma_{v_i}([a_{\infty},a_3,c,b_2,a_0,b_4,b_{\infty}] \cup [a_1,b_1])$ \\
& $J_3:\sigma_{v_i}([a_1,a_4,a_2,a_0,a_3,b_1] \cup [a_{\infty},b_{\infty}])$ \\
& $J_4:\sigma_{v_i}([b_{\infty},b_2,b_4,b_3,b_1] \cup [a_1,a_0,a_{\infty}])$ \\
& $J_5:\sigma_{v_i}([a_{\infty},c,a_2,b_1,a_4,b_{\infty}] \cup [a_0,b_3,b_0])$ \\
& $J_6:\sigma_{v_i}([a_{\infty},b_1,b_4,a_3,a_2,b_{\infty}] \cup [a_0,b_0])$ \\
& if $a_\y^{v_i}$ and $a_0^{v_i}$ are in the same component of $J_5\setminus\sigma_{v_i}(H')$.\\
3.7b& $J_1:\sigma_{v_i}([a_1,b_2,b_0,b_3,a_4,b_1] \cup [a_{\infty},c,b_{\infty}])$\\
& $J_2:\sigma_{v_i}([a_{\infty},a_3,c,b_2,a_0,b_4,b_{\infty}] \cup [a_1,b_1])$ \\
& $J_3:\sigma_{v_i}([a_1,a_4,a_2,a_0,a_3,b_1] \cup [a_{\infty},b_{\infty}])$ \\
& $J_4:\sigma_{v_i}([b_{\infty},b_2,b_4,b_3,b_1] \cup [a_1,a_0,a_{\infty}])$ \\
& $J_5:\sigma_{v_i}([b_{\infty},a_4,c,a_2,b_1,b_0] \cup [a_0,b_3,a_{\infty}])$ \\
& $J_6:\sigma_{v_i}([a_{\infty},b_1,b_4,a_3,a_2,b_{\infty}] \cup [a_0,b_0])$ \\
& if $a_\y^{v_i}$ and $b_\y^{v_i}$ are in the same component of $J_5\setminus\sigma_{v_i}(H')$.\\
\hline
3.8& $J_1:\sigma_{v_i}([a_{\infty},b_3,a_4,c,b_{\infty}] \cup [a_1,b_2,b_0,b_1])$\\
& $J_2:\sigma_{v_i}([a_{\infty},a_3,c,b_2,a_0,b_4,b_{\infty}] \cup [a_1,b_1])$ \\
& $J_3:\sigma_{v_i}([a_1,a_4,a_2,a_0,a_3,b_1] \cup [a_{\infty},b_{\infty}])$ \\
& $J_4:\sigma_{v_i}([b_{\infty},b_2,b_4,b_3,b_1] \cup [a_1,a_0,a_{\infty}])$ \\
& $J_5:\sigma_{v_i}([a_{\infty},c,a_2,b_1,a_4,b_{\infty}] \cup [a_0,b_3,b_0])$ \\
& $J_6:\sigma_{v_i}([a_{\infty},b_1,b_4,a_3,a_2,b_{\infty}] \cup [a_0,b_0])$ \\
\hline
\end{tabular}
\end{center}

Finally suppose $v_i\in U'\cap V'$. In this case $a_\y^{v_i}$ and $b_\y^{v_i}$ are in the same component of $J_2$, $J_4$ and $J_6$,
and we have eight cases for the remaining $2$-factors. Namely, for each of $J_1$, $J_3$ and $J_5$, $a_\y^{v_i}$ and $b_\y^{v_i}$
are either in the same component or they are not. We number these eight cases as in the following table. 

\vspace{0.3cm}

\begin{center}
\begin{tabular}{c|c|c|c|}
&$J_1$&$J_3$&$J_5$\\
\hline
4.1& same&same&same\\
\hline
4.2& distinct&same&same\\
\hline
4.3& same&distinct&same\\
\hline
4.4& same&same&distinct\\
\hline
4.5& same&distinct&distinct\\
\hline
4.6& distinct&same&distinct\\
\hline
4.7& distinct&distinct&same\\
\hline
4.8& distinct&distinct&distinct\\
\hline
\end{tabular}
\end{center}

Depending on which of cases 4.1--4.8 that we are in, 
we can reallocate the edges of $\sigma_{v_i}(H)$ to the factors $J_1,J_2,\ldots,J_6$ 
as indicated in the following table to obtain a new 
$2$-factorisation of $L$ with the desired properties. If $J_x$ ($x\in\{1,2,\ldots,6\}$) is not listed for a particular case, then the edges of 
$\sigma_{v_i}(H)$ that are in $J_x$ are unchanged. 

\begin{center}
\begin{tabular}{|c|l|}
\hline
4.1& \\
\hline
4.2& $J_1:\sigma_{v_i}([a_1,a_4,b_3,b_0,b_2,b_{\infty}] \cup [a_{\infty},b_1])$\\
& $J_5:\sigma_{v_i}([a_0,a_2,c,a_4,b_{\infty}] \cup [a_{\infty},b_3,b_1,b_0])$ \\
\hline
4.3& $J_3:\sigma_{v_i}([a_1,a_0,a_3,c,a_2,a_4,b_{\infty}] \cup [a_{\infty},b_1])$ \\
& $J_5:\sigma_{v_i}([a_{\infty},c,a_4,b_1,b_3,b_0] \cup [a_0,a_2,b_{\infty}])$ \\
\hline
4.4& $J_2:\sigma_{v_i}([a_1,b_2,b_4,a_0,a_3,b_1] \cup [a_{\infty},b_{\infty}])$ \\
& $J_3:\sigma_{v_i}([a_1,a_0,a_2,c,a_3,a_{\infty}] \cup [b_{\infty},a_4,b_1])$ \\
& $J_4:\sigma_{v_i}([a_{\infty},c,b_2,a_0,b_3,b_4,b_{\infty}] \cup [a_1,b_1])$ \\
& $J_5:\sigma_{v_i}([b_{\infty},c,a_4,a_2,b_1,b_3,b_0] \cup [a_0,a_{\infty}])$ \\
& $J_6:\sigma_{v_i}([a_{\infty},b_1,b_4,a_3,a_2,b_{\infty}] \cup [a_0,b_0])$ \\
\hline
4.5& $J_1:\sigma_{v_i}([a_1,b_2,b_0,b_3,a_{\infty}] \cup [b_{\infty},a_4,b_1])$ \\
& $J_2:\sigma_{v_i}([a_{\infty},a_3,a_0,b_2,b_4,b_{\infty}] \cup [a_1,b_1])$ \\
& $J_3:\sigma_{v_i}([a_1,a_4,c,a_3,b_1] \cup [a_{\infty},a_0,a_2,b_{\infty}])$ \\
& $J_4:\sigma_{v_i}([a_1,a_0,b_3,b_4,b_1] \cup [a_{\infty},c,b_2,b_{\infty}])$ \\
& $J_5:\sigma_{v_i}([a_{\infty},b_1,b_3,a_4,a_2,c,b_{\infty}] \cup [a_0,b_0])$ \\
& $J_6:\sigma_{v_i}([a_0,b_4,a_3,a_2,b_1,b_0] \cup [a_{\infty},b_{\infty}])$ \\
\hline
4.6& $J_1:\sigma_{v_i}([a_1,b_2,b_0,b_3,a_4,b_1] \cup [a_{\infty},b_{\infty}])$ \\
& $J_2:\sigma_{v_i}([a_{\infty},a_3,a_0,b_2,b_4,b_{\infty}] \cup [a_1,b_1])$ \\
& $J_3:\sigma_{v_i}([a_1,a_4,a_2,a_0,a_{\infty}] \cup [b_{\infty},c,a_3,b_1])$ \\
& $J_4:\sigma_{v_i}([a_1,a_0,b_4,b_3,b_1] \cup [a_{\infty},c,b_2,b_{\infty}])$ \\
& $J_5:\sigma_{v_i}([b_{\infty},a_4,c,a_2,b_1,b_0] \cup [a_0,b_3,a_{\infty}])$ \\
& $J_6:\sigma_{v_i}([a_{\infty},b_1,b_4,a_3,a_2,b_{\infty}] \cup [a_0,b_0])$ \\
\hline
4.7& $J_1:\sigma_{v_i}([a_{\infty},b_3,b_0,b_2,b_{\infty}] \cup [a_1,a_4,b_1])$ \\
& $J_3:\sigma_{v_i}([a_1,a_0,a_3,c,a_2,a_4,b_{\infty}] \cup [a_{\infty},b_1])$ \\
& $J_5:\sigma_{v_i}([a_{\infty},c,a_4,b_3,b_1,b_0] \cup [a_0,a_2,b_{\infty}])$ \\
\hline
4.8& $J_1:\sigma_{v_i}([a_{\infty},b_3,b_0,b_2,b_{\infty}] \cup [a_1,a_4,b_1])$ \\
& $J_2:\sigma_{v_i}([a_1,b_2,b_4,a_0,a_3,b_1] \cup [a_{\infty},b_{\infty}])$ \\
& $J_3:\sigma_{v_i}([a_{\infty},a_3,c,a_4,b_{\infty}] \cup [a_1,a_0,a_2,b_1])$ \\
& $J_4:\sigma_{v_i}([a_{\infty},c,b_2,a_0,b_3,b_4,b_{\infty}] \cup [a_1,b_1])$ \\
& $J_5:\sigma_{v_i}([b_{\infty},c,a_2,a_4,b_3,b_1,b_0] \cup [a_0,a_{\infty}])$ \\
& $J_6:\sigma_{v_i}([a_{\infty},b_1,b_4,a_3,a_2,b_{\infty}] \cup [a_0,b_0])$ \\
\hline
\end{tabular}
\end{center}

\noindent{\bf The case n=7:}

For $i=1,2,\ldots,7$, let $X_1(i)$, $Y_1(i)$, $X'_1(i)$, $Y'_1(i)$ $X_2(i)$, $Y_2(i)$, $X'_2(i)$ and $Y'_2(i)$ be the subgraphs of $K_{A_n\cup B_n\cup\{c\}}$ given by the union of
the paths listed in the following tables, and let $H^i$ be the subgraph of $K_{A_n\cup B_n\cup\{c\}}$ with edge set
$E(X_1(i))\cup E(Y_1(i))$. It can be checked that $\{H^1,H^2,\ldots,H^7\}$ is a decomposition of $K_{A_n\cup B_n\cup\{c\}}$,
and applying Lemma \ref{4RegularHamiltonFragments} with $H=H^i$, $X_1=X_1(i)$, $Y_1=Y_1(i)$, $X'_1 = X'_1(i)$, $Y'_1=Y'_1(i)$, $X_2=X_2(i)$, $Y_2=Y_2(i)$, $X'_2= X'_2(i)$, and $Y'_2 = Y'_2(i)$ shows that each 
$H^i$ is a Hamilton fragment. The value of $t$ can be deduced from the ends of the given paths.

\vspace{0.3cm}

$
\begin{array}{|c|c|}
\hline
X_1(1)&[a_{\infty},a_0,b_2,b_3,b_4,c,b_1,b_{\infty}] \cup [a_5,b_5]\\
\hline
Y_1(1)&[a_{\infty},b_1,b_2,a_3,b_4,b_5] \cup [a_5,a_0,b_{\infty}]\\
\hline
X'_1(1)& [a_{\infty},b_1,b_2,b_3,b_4,b_5] \cup [a_5,a_0,b_{\infty}]\\
\hline
Y'_1(1) & [a_{\infty},a_0,b_2,a_3,b_4,c,b_1,b_{\infty}] \cup [a_5,b_5]\\
\hline
X_2(1)&[a_{\infty},b_1,c,b_4,b_3,b_2,a_0,b_{\infty}] \cup [a_5,b_5]\\
\hline
Y_2(1)&[b_{\infty},b_1,b_2,a_3,b_4,b_5] \cup [a_5,a_0,a_{\infty}]\\
\hline
X'_2(1) & [b_{\infty},b_1,b_2,b_3,b_4,b_5] \cup [a_5,a_0,a_{\infty}]\\
\hline
Y'_2(1) & [a_{\infty},b_1,c,b_4,a_3,b_2,a_0,b_{\infty}] \cup [a_5,b_5] \\
\hline
\end{array}
$

$\begin{array}{|c|c|}
\hline
X_1(2)& [a_{\infty},c,a_2,a_3,a_1,b_5,b_0,b_{\infty}] \cup [a_4,b_4]\\
\hline
Y_1(2)&[b_{\infty},a_3,b_0,a_1,b_4] \cup [a_4,b_5,a_2,a_{\infty}]\\
\hline
X'_1(2)& Y_1(2) \\
\hline
Y'_1(2) & X_1(2)\\
\hline
X_2(2)& [a_{\infty},c,a_2,b_5,a_1,a_3,b_0,b_{\infty}] \cup [a_4,b_4]\\
\hline
Y_2(2)& [a_4,b_5,b_0,a_1,b_4] \cup [a_{\infty},a_2,a_3,b_{\infty}]\\
\hline
X'_2(2) & X_2(2) \\
\hline
Y'_2(2) & Y_2(2) \\
\hline
\end{array}
$

$\begin{array}{|c|c|}
\hline
X_1(3)& [a_{\infty},b_2,c,a_4,b_1,b_0,a_5,b_{\infty}] \cup [a_3,b_3]\\
\hline
Y_1(3)& [a_3,b_1,a_5,b_2,b_0,a_{\infty}] \cup [b_{\infty},a_4,b_3]\\
\hline
X'_1(3)& Y_1(3) \\
\hline
Y'_1(3) & X_1(3)\\
\hline
X_2(3)& [a_{\infty},b_0,b_1,a_5,b_2,c,a_4,b_{\infty}] \cup [a_3,b_3] \\
\hline
Y_2(3)& [a_{\infty},b_2,b_0,a_5,b_{\infty}] \cup [a_3,b_1,a_4,b_3] \\
\hline
X'_2(3) & X_2(3) \\
\hline
Y'_2(3) & Y_2(3) \\
\hline
\end{array}
$

$\begin{array}{|c|c|}
\hline
X_1(4)& [a_{\infty},b_5,b_1,b_4,a_0,c,b_3,b_{\infty}] \cup [a_2,b_2] \\
\hline
Y_1(4)& [a_2,b_1,a_0,b_5,b_3,a_{\infty}] \cup [b_{\infty},b_4,b_2]\\
\hline
X'_1(4)& Y_1(4) \\
\hline
Y'_1(4) & X_1(4)\\
\hline
X_2(4)& [a_{\infty},b_3,c,a_0,b_5,b_1,b_4,b_{\infty}] \cup [a_2,b_2] \\
\hline
Y_2(4)& [a_2,b_1,a_0,b_4,b_2] \cup [a_{\infty},b_5,b_3,b_{\infty}]\\
\hline
X'_2(4) & X_2(4) \\
\hline
Y'_2(4) & Y_2(4) \\
\hline
\end{array}
$

$\begin{array}{|c|c|}
\hline
X_1(5)& [a_{\infty},a_5,b_3,a_2,b_4,b_0,c,b_{\infty}] \cup [a_1,b_1]\\
\hline
Y_1(5)& [b_{\infty},a_2,b_0,b_3,b_1] \cup [a_1,a_5,b_4,a_{\infty}]\\
\hline
X'_1(5)& Y_1(5) \\
\hline
Y'_1(5) & X_1(5)\\
\hline
X_2(5)& [a_{\infty},a_5,b_4,a_2,b_3,b_0,c,b_{\infty}] \cup [a_1,b_1]\\
\hline
Y_2(5)& [a_{\infty},b_4,b_0,a_2,b_{\infty}] \cup [a_1,a_5,b_3,b_1] \\
\hline
X'_2(5) & X_2(5) \\
\hline
Y'_2(5) & Y_2(5) \\
\hline
\end{array}
$

$\begin{array}{|c|c|}
\hline
X_1(6)& [a_{\infty},a_1,c,b_5,a_3,a_4,b_2,b_{\infty}] \cup [a_0,b_0]\\
\hline
Y_1(6)& [b_{\infty},b_5,b_2,a_1,a_4,b_0] \cup [a_0,a_3,a_{\infty}]\\
\hline
X'_1(6)& Y_1(6) \\
\hline
Y'_1(6) & X_1(6)\\
\hline
X_2(6)& [a_{\infty},a_3,b_5,c,a_1,a_4,b_2,b_{\infty}] \cup [a_0,b_0] \\
\hline
Y_2(6)& [a_{\infty},a_1,b_2,b_5,b_{\infty}] \cup [a_0,a_3,a_4,b_0]\\
\hline
X'_2(6) & X_2(6) \\
\hline
Y'_2(6) & Y_2(6) \\
\hline
\end{array}
$

$\begin{array}{|c|c|}
\hline
X_1(7)& [a_3,c,a_5,a_4,a_2,a_1,a_0,b_3] \cup [a_{\infty},b_{\infty}]\\
\hline
Y_1(7)&[a_3,a_5,a_2,a_0,a_4,a_{\infty}] \cup [b_{\infty},a_1,b_3] \\
\hline
X'_1(7)& Y_1(7) \\
\hline
Y'_1(7) & X_1(7)\\
\hline
X_2(7)& [a_3,c,a_5,a_4,a_2,a_0,a_1,b_3] \cup [a_{\infty},b_{\infty}]\\
\hline
Y_2(7)& [a_3,a_5,a_2,a_1,b_{\infty}] \cup [a_{\infty},a_4,a_0,b_3]\\
\hline
X'_2(7) & X_2(7) \\
\hline
Y'_2(7) & Y_2(7) \\
\hline
\end{array}
$

\noindent{\bf The case n=8:}

For $i=1,2,\ldots,8$, let $X_1(i)$, $Y_1(i)$, $X'_1(i)$, $Y'_1(i)$ $X_2(i)$, $Y_2(i)$, $X'_2(i)$ and $Y'_2(i)$ be the subgraphs of $K_{A_n\cup B_n\cup\{c\}}$ given by the union of
the paths listed in the following tables, and let $H^i$ be the subgraph of $K_{A_n\cup B_n\cup\{c\}}$ with edge set
$E(X_1(i))\cup E(Y_1(i))$. It can be checked that $\{H^1,H^2,\ldots,H^8\}$ is a decomposition of $K_{A_n\cup B_n\cup\{c\}}$,
and applying Lemma \ref{4RegularHamiltonFragments} with $H=H^i$, $X_1=X_1(i)$, $Y_1=Y_1(i)$, $X'_1 = X'_1(i)$, $Y'_1=Y'_1(i)$, $X_2=X_2(i)$, $Y_2=Y_2(i)$, $X'_2= X'_2(i)$, and $Y'_2 = Y'_2(i)$ shows that each 
$H^i$ is a Hamilton fragment. The value of $t$ can be deduced from the ends of the given paths.

$\begin{array}{|c|c|}
\hline
X_1(1)&[a_{\infty},a_0,c,b_4,b_5,b_3,b_2,b_1,b_{\infty}] \cup [a_6,b_6] \\
\hline
Y_1(1)& [a_{\infty},b_1,a_2,a_3,a_5,b_4,b_6] \cup [a_6,a_0,b_{\infty}]\\
\hline
X'_1(1)& [a_{\infty},b_1,b_2,b_3,b_5,b_4,b_6] \cup [a_6,a_0,b_{\infty}] \\
\hline
Y'_1(1) & [a_{\infty},a_0,c,b_4,a_5,a_3,a_2,b_1,b_{\infty}] \cup [a_6,b_6]\\
\hline
X_2(1)& [a_{\infty},b_1,b_2,b_3,b_5,b_4,c,a_0,b_{\infty}] \cup [a_6,b_6]\\
\hline
Y_2(1)& [b_{\infty},b_1,a_2,a_3,a_5,b_4,b_6] \cup [a_6,a_0,a_{\infty}]\\
\hline
X'_2(1) & [b_{\infty},b_1,b_2,b_3,b_5,b_4,b_6] \cup [a_6,a_0,a_{\infty}]\\
\hline
Y'_2(1) &[a_{\infty},b_1,a_2,a_3,a_5,b_4,c,a_0,b_{\infty}] \cup [a_6,b_6]\\
\hline
\end{array}
$

$\begin{array}{|c|c|}
\hline
X_1(2)& [a_{\infty},a_2,c,b_3,b_6,b_1,b_4,b_0,b_{\infty}] \cup [a_5,b_5]\\
\hline
Y_1(2)& [a_5,b_3,a_6,a_4,a_1,b_0,a_{\infty}] \cup [b_{\infty},a_2,b_5]\\
\hline
X'_1(2)& [a_5,b_3,b_6,b_1,b_4,b_0,a_{\infty}] \cup [b_{\infty},a_2,b_5] \\
\hline
Y'_1(2) & [a_{\infty},a_2,c,b_3,a_6,a_4,a_1,b_0,b_{\infty}] \cup [a_5,b_5]\\
\hline
X_2(2)& [a_{\infty},b_0,b_4,b_1,b_6,b_3,c,a_2,b_{\infty}] \cup [a_5,b_5]\\
\hline
Y_2(2)& [a_5,b_3,a_6,a_4,a_1,b_0,b_{\infty}] \cup [a_{\infty},a_2,b_5]\\
\hline
X'_2(2) & [a_5,b_3,b_6,b_1,b_4,b_0,b_{\infty}] \cup [a_{\infty},a_2,b_5] \\
\hline
Y'_2(2) &[a_{\infty},b_0,a_1,a_4,a_6,b_3,c,a_2,b_{\infty}] \cup [a_5,b_5]\\
\hline
\end{array}
$

$\begin{array}{|c|c|}
\hline
X_1(3)& [a_{\infty},c,b_0,b_6,b_2,a_3,b_1,a_5,b_{\infty}] \cup [a_4,b_4]\\
\hline
Y_1(3)& [a_4,b_1,b_0,a_5,b_6,a_3,a_{\infty}] \cup [b_{\infty},b_2,b_4]\\
\hline
X'_1(3)& Y_1(3) \\
\hline
Y'_1(3) & X_1(3)\\
\hline
X_2(3)& [a_{\infty},c,b_0,a_5,b_1,a_3,b_6,b_2,b_{\infty}] \cup [a_4,b_4]\\
\hline
Y_2(3)& [a_4,b_1,b_0,b_6,a_5,b_{\infty}] \cup [a_{\infty},a_3,b_2,b_4]\\
\hline
X'_2(3) & X_2(3) \\
\hline
Y'_2(3) & Y_2(3) \\
\hline
\end{array}
$

$\begin{array}{|c|c|}
\hline
X_1(4)& [a_{\infty},b_6,a_1,a_5,b_2,a_0,a_4,b_{\infty}] \cup [a_3,b_3]\\
\hline
Y_1(4)& [a_3,c,b_2,a_4,a_5,a_{\infty}] \cup [b_{\infty},b_6,a_0,a_1,b_3]\\
\hline
X'_1(4)& Y_1(4) \\
\hline
Y'_1(4) & X_1(4)\\
\hline
X_2(4)& [a_{\infty},b_6,a_1,a_0,b_2,a_5,a_4,b_{\infty}] \cup [a_3,b_3]\\
\hline
Y_2(4)& [a_3,c,b_2,a_4,a_0,b_6,b_{\infty}] \cup [a_{\infty},a_5,a_1,b_3]\\
\hline
X'_2(4) & X_2(4) \\
\hline
Y'_2(4) & Y_2(4) \\
\hline
\end{array}
$

$\begin{array}{|c|c|}
\hline
X_1(5)&[a_{\infty},b_5,c,b_1,a_6,b_0,a_4,b_3,b_{\infty}] \cup [a_2,b_2] \\
\hline
Y_1(5)& [b_{\infty},a_6,b_5,b_1,b_3,b_0,b_2] \cup [a_2,a_4,a_{\infty}]\\
\hline
X'_1(5)& Y_1(5) \\
\hline
Y'_1(5) & X_1(5)\\
\hline
X_2(5)& [a_{\infty},a_4,b_0,a_6,b_5,c,b_1,b_3,b_{\infty}] \cup [a_2,b_2]\\
\hline
Y_2(5)& [a_2,a_4,b_3,b_0,b_2] \cup [a_{\infty},b_5,b_1,a_6,b_{\infty}]\\
\hline
X'_2(5) & X_2(5) \\
\hline
Y'_2(5) & Y_2(5) \\
\hline
\end{array}
$

$\begin{array}{|c|c|}
\hline
X_1(6)& [a_{\infty},a_6,a_5,a_2,b_3,a_0,b_4,b_{\infty}] \cup [a_1,b_1]\\
\hline
Y_1(6)& [a_1,a_6,a_2,b_4,b_3,a_{\infty}] \cup [b_{\infty},c,a_5,a_0,b_1]\\
\hline
X'_1(6)& Y_1(6) \\
\hline
Y'_1(6) & X_1(6)\\
\hline
X_2(6)& [a_{\infty},a_6,a_5,a_0,b_3,a_2,b_4,b_{\infty}] \cup [a_1,b_1]\\
\hline
Y_2(6)& [a_1,a_6,a_2,a_5,c,b_{\infty}] \cup [a_{\infty},b_3,b_4,a_0,b_1]\\
\hline
X'_2(6) & X_2(6) \\
\hline
Y'_2(6) & Y_2(6) \\
\hline
\end{array}
$

$\begin{array}{|c|c|}
\hline
X_1(7)& [a_{\infty},b_2,b_5,a_3,b_4,a_6,c,a_1,b_{\infty}] \cup [a_0,b_0]\\
\hline
Y_1(7)& [a_0,a_3,a_6,b_2,a_1,b_4,a_{\infty}] \cup [b_{\infty},b_5,b_0]\\
\hline
X'_1(7)& Y_1(7) \\
\hline
Y'_1(7) & X_1(7)\\
\hline
X_2(7)& [a_{\infty},b_2,a_1,c,a_6,b_4,a_3,b_5,b_{\infty}] \cup [a_0,b_0]\\
\hline
Y_2(7)& [a_0,a_3,a_6,b_2,b_5,b_0] \cup [a_{\infty},b_4,a_1,b_{\infty}]\\
\hline
X'_2(7) & X_2(7) \\
\hline
Y'_2(7) & Y_2(7) \\
\hline
\end{array}
$

$\begin{array}{|c|c|}
\hline
X_1(8)& [a_0,b_5,b_6,a_4,a_3,a_1,a_2,b_0] \cup [a_{\infty},b_{\infty}]\\
\hline
Y_1(8)& [a_0,a_2,b_6,c,a_4,b_5,a_1,a_{\infty}] \cup [b_{\infty},a_3,b_0]\\
\hline
X'_1(8)& Y_1(8) \\
\hline
Y'_1(8) & X_1(8)\\
\hline
X_2(8)& [a_0,b_5,a_4,b_6,a_2,a_1,a_3,b_0] \cup [a_{\infty},b_{\infty}]\\
\hline
Y_2(8)& [a_{\infty},a_1,b_5,b_6,c,a_4,a_3,b_{\infty}] \cup [a_0,a_2,b_0]\\
\hline
X'_2(8) & X_2(8) \\
\hline
Y'_2(8) & Y_2(8) \\
\hline
\end{array}
$

\noindent{\bf The case n=9:}

For $i=1,2,\ldots,9$, let $X_1(i)$, $Y_1(i)$, $X'_1(i)$, $Y'_1(i)$ $X_2(i)$, $Y_2(i)$, $X'_2(i)$ and $Y'_2(i)$ be the subgraphs of $K_{A_n\cup B_n\cup\{c\}}$ given by the union of
the paths listed in the following tables, and let $H^i$ be the subgraph of $K_{A_n\cup B_n\cup\{c\}}$ with edge set
$E(X_1(i))\cup E(Y_1(i))$. It can be checked that $\{H^1,H^2,\ldots,H^9\}$ is a decomposition of $K_{A_n\cup B_n\cup\{c\}}$,
and applying Lemma \ref{4RegularHamiltonFragments} with $H=H^i$, $X_1=X_1(i)$, $Y_1=Y_1(i)$, $X'_1 = X'_1(i)$, $Y'_1=Y'_1(i)$, $X_2=X_2(i)$, $Y_2=Y_2(i)$, $X'_2= X'_2(i)$, and $Y'_2 = Y'_2(i)$ shows that each 
$H^i$ is a Hamilton fragment. The value of $t$ can be deduced from the ends of the given paths.

$\begin{array}{|c|c|}
\hline
X_1(1)& [a_{\infty},a_0,b_2,b_4,b_5,b_6,c,a_3,b_1,b_{\infty}]\cup [a_7,b_7]\\
\hline
Y_1(1)& [b_{\infty},a_0,b_4,a_5,b_6,b_7] \cup [a_7,a_3,b_2,b_1,a_{\infty}]\\
\hline
X'_1(1)& [b_{\infty},a_0,b_4,b_5,b_6,b_7] \cup [a_7,a_3,b_2,b_1,a_{\infty}] \\
\hline
Y'_1(1) & [a_{\infty},a_0,b_2,b_4,a_5,b_6,c,a_3,b_1,b_{\infty}]\cup [a_7,b_7]\\
\hline
X_2(1)& [a_{\infty},b_1,a_3,c,b_6,b_5,b_4,b_2,a_0,b_{\infty}] \cup [a_7,b_7]\\
\hline
Y_2(1)& [a_{\infty},a_0,b_4,a_5,b_6,b_7] \cup [a_7,a_3,b_2,b_1,b_{\infty}]\\
\hline
X'_2(1) & [a_{\infty},a_0,b_4,b_5,b_6,b_7] \cup [a_7,a_3,b_2,b_1,b_{\infty}] \\
\hline
Y'_2(1) & [a_{\infty},b_1,a_3,c,b_6,a_5,b_4,b_2,a_0,b_{\infty}] \cup [a_7,b_7] \\
\hline
\end{array}
$

$\begin{array}{|ll||ll|}
\hline
X_1(2): & [a_{\infty},b_2,b_5,a_7,c,b_3,a_4,a_1,b_0,b_{\infty}]\cup [a_6,b_6]& X'_1(2):& Y_1(2) \\
\hline
Y_1(2): & [b_{\infty},b_2,a_1,b_3,b_0,b_6]\cup [a_6,a_7,a_4,b_5,a_{\infty}] & 
Y'_1(2): & X_1(2)\\
\hline
X_2(2): & [a_{\infty},b_5,b_2,a_1,a_4,a_7,c,b_3,b_0,b_{\infty}]\cup [a_6,b_6] &
X'_2(2):& X_2(2)\\
\hline
Y_2(2): & [a_6,a_7,b_5,a_4,b_3,a_1,b_0,b_6]\cup [a_{\infty},b_2,b_{\infty}] &
Y'_2(2) :& Y_2(2)\\
\hline
\end{array}
$

$\begin{array}{|c|c|}
\hline
X_1(3)&[a_{\infty},c,b_4,b_7,a_6,a_2,a_0,b_1,b_3,b_{\infty}]\cup [a_5,b_5]\\
\hline
Y_1(3)& [b_{\infty},b_4,b_1,a_2,b_3,b_5]\cup [a_5,b_7,a_0,a_6,a_{\infty}]\\
\hline
X'_1(3)& Y_1(3) \\
\hline
Y'_1(3) & X_1(3)\\
\hline
X_2(3)& [a_{\infty},c,b_4,b_7,a_6,a_0,a_2,b_1,b_3,b_{\infty}]\cup [a_5,b_5]\\
\hline
Y_2(3)& [a_5,b_7,a_0,b_1,b_4,b_{\infty}]\cup [a_{\infty},a_6,a_2,b_3,b_5]\\
\hline
X'_2(3) & X_2(3) \\
\hline
Y'_2(3) & Y_2(3) \\
\hline
\end{array}
$

$\begin{array}{|c|c|}
\hline
X_1(4)& [a_{\infty},b_6,b_1,b_5,b_0,a_7,a_2,a_3,b_{\infty}]\cup [a_4,b_4]\\
\hline
Y_1(4)& [b_{\infty},b_5,a_2,b_6,a_3,b_0,b_4]\cup [a_4,c,b_1,a_7,a_{\infty}]\\
\hline
X'_1(4)& Y_1(4) \\
\hline
Y'_1(4) & X_1(4)\\
\hline
X_2(4)& [a_{\infty},a_7,b_1,b_6,a_2,b_5,b_0,a_3,b_{\infty}]\cup [a_4,b_4]\\
\hline
Y_2(4)& [a_{\infty},b_6,a_3,a_2,a_7,b_0,b_4]\cup [a_4,c,b_1,b_5,b_{\infty}]\\
\hline
X'_2(4) & X_2(4) \\
\hline
Y'_2(4) & Y_2(4) \\
\hline
\end{array}
$

$\begin{array}{|c|c|}
\hline
X_1(5)& [a_{\infty},a_1,c,a_0,a_4,b_7,a_2,a_5,a_6,b_{\infty}]\cup [a_3,b_3]\\
\hline
Y_1(5)& [a_3,a_0,a_5,a_1,b_7,a_{\infty}]\cup [b_{\infty},a_2,a_4,a_6,b_3]\\
\hline
X'_1(5)& Y_1(5) \\
\hline
Y'_1(5) & X_1(5)\\
\hline
X_2(5)&[a_{\infty},a_1,c,a_0,a_5,a_2,b_7,a_4,a_6,b_{\infty}]\cup [a_3,b_3] \\
\hline
Y_2(5)& [a_{\infty},b_7,a_1,a_5,a_6,b_3]\cup [a_3,a_0,a_4,a_2,b_{\infty}]\\
\hline
X'_2(5) & X_2(5) \\
\hline
Y'_2(5) & Y_2(5) \\
\hline
\end{array}
$

$\begin{array}{|c|c|}
\hline
X_1(6)& [a_{\infty},a_5,b_3,b_6,a_4,b_1,b_0,c,b_7,b_{\infty}]\cup [a_2,b_2]\\
\hline
Y_1(6)& [a_2,b_0,a_4,a_5,b_1,b_7,b_3,a_{\infty}]\cup [b_{\infty},b_6,b_2]\\
\hline
X'_1(6)& Y_1(6) \\
\hline
Y'_1(6) & X_1(6)\\
\hline
X_2(6)& [a_{\infty},a_5,a_4,b_1,b_0,c,b_7,b_3,b_6,b_{\infty}]\cup [a_2,b_2]\\
\hline
Y_2(6)& [a_{\infty},b_3,a_5,b_1,b_7,b_{\infty}]\cup [a_2,b_0,a_4,b_6,b_2]\\
\hline
X'_2(6) & X_2(6) \\
\hline
Y'_2(6) & Y_2(6) \\
\hline
\end{array}
$

$\begin{array}{|c|c|}
\hline
X_1(7)& [a_{\infty},a_4,a_3,b_7,b_2,b_0,a_6,b_5,c,b_{\infty}]\cup [a_1,b_1]\\
\hline
Y_1(7)& [a_1,a_3,b_5,b_7,b_0,a_{\infty}]\cup [b_{\infty},a_4,b_2,a_6,b_1]\\
\hline
X'_1(7)& Y_1(7) \\
\hline
Y'_1(7) & X_1(7)\\
\hline
X_2(7)& [a_{\infty},b_0,b_7,a_3,a_4,b_2,a_6,b_5,c,b_{\infty}]\cup [a_1,b_1]\\
\hline
Y_2(7)& [a_1,a_3,b_5,b_7,b_2,b_0,a_6,b_1]\cup [a_{\infty},a_4,b_{\infty}]\\
\hline
X'_2(7) & X_2(7) \\
\hline
Y'_2(7) & Y_2(7) \\
\hline
\end{array}
$

$\begin{array}{|c|c|}
\hline
X_1(8)& [a_{\infty},a_2,b_4,a_1,a_6,a_3,a_5,a_7,b_{\infty}]\cup [a_0,b_0]\\
\hline
Y_1(8)& [a_0,a_7,a_1,a_2,c,a_6,b_4,a_3,a_{\infty}]\cup [b_{\infty},a_5,b_0]\\
\hline
X'_1(8)& Y_1(8) \\
\hline
Y'_1(8) & X_1(8)\\
\hline
X_2(8)& [a_{\infty},a_3,a_6,b_4,a_2,a_1,a_7,a_5,b_{\infty}]\cup [a_0,b_0]\\
\hline
Y_2(8)& [a_{\infty},a_2,c,a_6,a_1,b_4,a_3,a_5,b_0]\cup [a_0,a_7,b_{\infty}]\\
\hline
X'_2(8) & X_2(8) \\
\hline
Y'_2(8) & Y_2(8) \\
\hline
\end{array}
$

$\begin{array}{|c|c|}
\hline
X_1(9)& [a_5,c,b_2,b_3,a_7,b_4,b_6,a_1,a_0,b_5] \cup [a_{\infty},b_{\infty}]\\
\hline
Y_1(9)& [a_5,b_2,a_7,b_6,a_0,b_3,b_4,a_{\infty}] \cup [b_{\infty},a_1,b_5]\\
\hline
X'_1(9)& Y_1(9) \\
\hline
Y'_1(9) & X_1(9)\\
\hline
X_2(9)& [a_5,c,b_2,a_7,b_3,b_4,b_6,a_0,a_1,b_5]\cup [a_{\infty},b_{\infty}]\\
\hline
Y_2(9)&[a_{\infty},b_4,a_7,b_6,a_1,b_{\infty}]\cup [a_5,b_2,b_3,a_0,b_5] \\
\hline
X'_2(9) & X_2(9) \\
\hline
Y'_2(9) & Y_2(9) \\
\hline
\end{array}
$
\vspace{0.3cm}

\noindent{\large\bf Acknowledgements}
The authors acknowledge the support of the Australian Research Council via grants
DP150100530, DP150100506, DP120100790, DP120103067 and DP130102987.

\end{document}